# Classification of irreducible holonomies of torsion-free affine connections

By Sergei Merkulov and Lorenz Schwachhöfer

## Contents

## 1. Introduction

One of the most informative characteristics of an affine connection on a smooth connected manifold $M$ is its holonomy group $\mathcal{H}$. At any particular point $p \in M$ one may represent $\mathcal{H}$ as the set of all linear automorphisms of the associated tangent space $T_pM$ which are induced by parallel translation along $p$-based loops. If $M$ is simply connected, then this set is known to be a Lie subgroup of $\mathrm{End}(T_pM)$. The notion of the holonomy group was introduced in the 1920s by Élie Cartan who used it to classify all Riemannian locally symmetric spaces [19]. This concept is widely used in many branches of mathematics and theoretical physics. One of the fundamental questions in differential geometry is the *irreducible holonomy problem*:

Which irreducible Lie subgroups $G \subset \mathrm{End}(V)$, where $V$ is a real or complex finite dimensional vector space, can occur as the holonomy group of a torsion-free affine connection?

The requirement of *torsion-freeness* is imposed because of a result of Hano and Ozeki [23] which says that *any* (closed) Lie subgroup $G \subset \mathrm{End}(V)$ can occur as the holonomy of some affine connection (with torsion, in general).

If a connection is (locally) symmetric, then its holonomy group, if irreducible, equals essentially the (local) isotropy group in the transvection group. Hence, the holonomy problem for these connections is equivalent to the classification problem of symmetric spaces which was solved long ago [19], [9]. Therefore, throughout this paper we shall use the term *holonomy group* for an irreducible Lie subgroup $G \subset \mathrm{End}(V)$ which can occur as the holonomy of a torsion-free affine connection which is *not* locally symmetric.

The next milestone toward the solution of the holonomy problem was set in 1955 by Berger [8]. He gave necessary conditions which a holonomy group must satisfy, and then classified the groups satisfying these conditions. His classification was given in two parts.

The first part ("metric part") consists of all holonomy groups which leave some symmetric bilinear form invariant. The connections with one of these holonomies are always Levi-Civita connections of a (pseudo-)Riemannian metric on $M$ and have been studied extensively in the intervening years. In fact, it is by now well-known which of the entries in this list do actually occur as holonomy groups and how the holonomy relates to both the geometry and the topology of the underlying manifold. We summarize the possible metric holonomies in Table 1 which represents a culmination of efforts of many people (see e.g. [4], [10], [14], [16], [36] and the references cited therein).

The second part ("non-metric part") was stated to contain all remaining holonomy groups, *up to a finite number of missing terms*. It is given in Table 2, modulo some slight modifications due to Bryant [16]. The holonomy groups which are missing from this list are called *exotic holonomies.* That the list of exotic holonomies is nonempty was established by Bryant [15]; several other (but still finitely many) exotic holonomies followed [16], [22], [21]. However, in [20], Chi et al. discovered an *infinite family* of exotic holonomies, thereby showing the incompleteness of the nonmetric part of Berger's list.

**Table 1**: Complete list of metric holonomies

| Group $G$ | Representation space | Group $G$ | Representation space |
|---|---|---|---|
| $\mathrm{SO}(p,q)$ | $\mathbb{R}^{p+q}$, $p+q \geqslant 3$ | $\mathrm{G}_2$ | $\mathbb{R}^7$ |
| $\mathrm{SO}(n,\mathbb{C})$ | $\mathbb{C}^n \simeq \mathbb{R}^{2n}$, $n \geqslant 3$ | $\mathrm{G}_2'$ | $\mathbb{R}^7$ |
| $\mathrm{U}(p,q)$ | $\mathbb{C}^{p+q} \simeq \mathbb{R}^{2(p+q)}$, $p+q \geqslant 2$ | $\mathrm{G}_2^{\mathbb{C}}$ | $\mathbb{R}^{14}$ |
| $\mathrm{SU}(p,q)$ | $\mathbb{C}^{p+q} \simeq \mathbb{R}^{2(p+q)}$, $\begin{cases} p+q \geqslant 2 \\ (p,q) \neq (1,1) \end{cases}$ | $\mathrm{Spin}(7)$ | $\mathbb{R}^8$ |
| $\mathrm{Sp}(p,q)\cdot\mathrm{Sp}(1)$ | $H^{p+q} \simeq \mathbb{R}^{4(p+q)}$, $p+q \geqslant 2$ | $\mathrm{Spin}(4,3)$ | $\mathbb{R}^8$ |
| $\mathrm{Sp}(p,q)$ | $H^{p+q} \simeq \mathbb{R}^{4(p+q)}$, $p+q \geqslant 2$ | $\mathrm{Spin}(7,\mathbb{C})$ | $\mathbb{R}^{16}$ |
| $\mathrm{Sp}(n,\mathbb{R})\cdot\mathrm{SL}(2,\mathbb{R})$ | $\mathbb{R}^2 \otimes \mathbb{R}^{2n} \simeq \mathbb{R}^{4n}$, $n \geqslant 2$ | | |
| $\mathrm{Sp}(n,\mathbb{C})\cdot\mathrm{SL}(2,\mathbb{C})$ | $\mathbb{C}^2 \otimes \mathbb{C}^{2n} \simeq \mathbb{R}^{8n}$, $n \geqslant 2$ | | |

In Table 3 we list all known exotic holonomies, including some new examples; that is, we have the following result.

Theorem A. (i) *The representations of the groups*

$$\mathrm{Sp}(3,\mathbb{R}),\ \mathrm{Sp}(3,\mathbb{C}),\ \mathrm{SL}(6,\mathbb{R}),\ \mathrm{SU}(1,5),\ \mathrm{SU}(3,3),\ \mathrm{SL}(6,\mathbb{C}),$$
$$\mathrm{Spin}(2,10),\ \mathrm{Spin}(6,6),\ \mathrm{Spin}(12,\mathbb{C})$$

*listed in Table* 3 *occur as holonomies of torsion-free affine connections.*

(ii) *Any torsion-free affine connection whose holonomy is one of those listed in* (i) *is analytic.*

(iii) *The moduli space of torsion-free affine connections whose holonomy is one of those listed in* (i) *is finite dimensional.*

These properties are due to the fact that all connections with these holonomies are constructed from certain quadratic deformations of linear Poisson structures; this method was developed for the construction of the infinite family of exotic holonomies [20].

**Table 2**: Berger's list of non-metric holonomies, modified by Bryant [16]

| Group $G$ | Representation $V$ | Restrictions |
|---|---|---|
| $T_{\mathbb{R}} \cdot SL(n,\mathbb{R})$ | $\mathbb{R}^n$ | $n \geqslant 2$ |
| | $\odot^2\mathbb{R}^n \simeq \mathbb{R}^{n(n+1)/2}$ | $n \geqslant 3$ |
| | $\Lambda^2\mathbb{R}^n \simeq \mathbb{R}^{n(n-1)/2}$ | $n \geqslant 5$ |
| $T_{\mathbb{C}} \cdot SL(n,\mathbb{C})$ | $\mathbb{C}^n \simeq \mathbb{R}^{2n}$ | $n \geqslant 1$ |
| | $\odot^2\mathbb{C}^n \simeq \mathbb{R}^{n(n+1)}$ | $n \geqslant 3$ |
| | $\Lambda^2\mathbb{C}^n \simeq \mathbb{R}^{n(n-1)}$ | $n \geqslant 5$ |
| $T_{\mathbb{R}} \cdot SL(n,\mathbb{C})$ | $\{A \in M_n(\mathbb{C}) : A = A^*\} \simeq \mathbb{R}^{n^2}$ | $n \geqslant 3$ |
| $T_{\mathbb{R}} \cdot SL(n,\mathbb{H})$ | $\mathbb{H}^n \simeq \mathbb{R}^{4n}$ | $n \geqslant 1$ |
| | $\{A \in M_n(\mathbb{H}) : A = -A^*\} \simeq \mathbb{R}^{n(2n+1)}$ | $n \geqslant 2$ |
| | $\{A \in M_n(\mathbb{H}) : A = A^*\} \simeq \mathbb{R}^{n(2n-1)}$ | $n \geqslant 3$ |
| $Sp(n,\mathbb{R})$ | $\mathbb{R}^{2n}$ | $n \geqslant 2$ |
| $Sp(n,\mathbb{C})$ | $\mathbb{C}^{2n} \simeq \mathbb{R}^{4n}$ | $n \geqslant 2$ |
| $\mathbb{R}^* \cdot SO(p,q)$ | $\mathbb{R}^{p+q}$ | $p+q \geqslant 3$ |
| $T_{\mathbb{C}} \cdot SO(n,\mathbb{C})$, $T_{\mathbb{C}} \neq 0$ | $\mathbb{C}^n \simeq \mathbb{R}^{2n}$ | $n \geqslant 3$ |
| $T_{\mathbb{R}} \cdot SL(m,\mathbb{R}) \cdot SL(n,\mathbb{R})$ | $\mathbb{R}^m \otimes \mathbb{R}^n \simeq \mathbb{R}^{mn}$ | $m > n \geqslant 2$ or $m \geqslant n > 2$ |
| $T_{\mathbb{C}} \cdot SL(m,\mathbb{C}) \cdot SL(n,\mathbb{C})$ | $\mathbb{C}^m \otimes \mathbb{C}^n \simeq \mathbb{R}^{2mn}$ | $m > n \geqslant 2$ or $m \geqslant n > 2$ |
| $T_{\mathbb{R}} \cdot SL(m,\mathbb{H}) \cdot SL(n,\mathbb{H})$ | $\mathbb{H}^m \otimes \mathbb{H}^n \simeq \mathbb{R}^{4mn}$ | $m > n \geqslant 1$ or $m \geqslant n > 1$ |
| Notation: $T_{\mathbb{F}}$ denotes any connected Lie subgroup of $\mathbb{F}^*$, $M_n(\mathbb{F})$ denotes the algebra of $n \times n$ matrices with entries in $\mathbb{F}$. | | |

Our main result is to show that these are the last entries missing from Berger's list, thus completing the holonomy classification.

Main Theorem. *If $G \subset \mathrm{End}(V)$ is an irreducible Lie subgroup which occurs as the holonomy of a torsion-free affine connection which is not locally symmetric, then $G$ is one of the entries in Tables* 1–3.

Moreover, due to the efforts of a number of people over the last 40 years (see the fundamental papers of Bryant [13]–[16], the most recent works on exotic holonomies [20], [21], the books [10], [36] and the references cited therein), all entries of Tables 1–3 are known to occur as holonomies, except for the 4-dimensional representations of $H_\lambda \cdot SU(2)$ and $H_\lambda \cdot SL(2,\mathbb{R})$ which were discovered as candidates for holonomy groups by Bryant [16].

There is no apparent pattern in Tables 1–3. A certain pattern does emerge, however, when one switches from the usual representation of holonomy as a pair (Lie group $G$, representation space $V$) to the Borel-Weil picture (homo-geneous manifold $X$, line bundle $L \to X$) which is defined as follows. To any irreducible non-abelian real matrix subgroup $G \subset \mathrm{End}(V)$, one first associates, by complexifying $G$ or complexifying both $G$ and $V$, an irreducible complex matrix subgroup $G_{\mathbb{C}} \subset \mathrm{End}(V_{\mathbb{C}})$. Then the projectivized $G_{\mathbb{C}}$-orbit $X \subset \mathbb{P}(V_{\mathbb{C}})$ of a maximal weight vector is well defined and is called the *sky of* $G$; define $L = \mathcal{O}(1) \mid_X$. Clearly, for any initial $G \subset \mathrm{End}(V)$ the associated sky $X$ is a compact complex homogeneous-rational manifold. Then the main theorem easily implies the following:

Corollary. *Let $G \subset \mathrm{End}(V)$ be an irreducible (real or complex) Lie subgroup. If $G$ occurs as the holonomy group of a torsion-free connection which is not locally symmetric, then its sky is biholomorphic to a compact hermitean-symmetric manifold.*

Remarkably, the implication Main Theorem $\Longrightarrow$ Corollary can be reversed. This is a central idea of our approach to the solution of the holonomy problem. In a sense, this is a holomorphic analogue of Simons' [37] approach to the classification of *compact* irreducible holonomies of Levi-Civita connections.

Using the above corollary as a working hypothesis we first were able to reproduce all the known holonomy groups in Tables 1–3 and second find new examples of Theorem A. This part of our work is explained in Sections 2-4 of the paper. Then using a combination of standard representation-theoretic methods (applicable in the picture $(G, V)$) with the methods of complex analysis and the Bott-Borel-Weil theorem (which are applicable in picture $(X, L)$) we showed that the working hypothesis gives rise to *all* possible holonomies, thereby completing the Berger classification (and giving, in particular, a new proof of the entire original Berger list). This second part of our work is explained in Sections 5 and 6 of the paper.

While working on the solution of the holonomy problem, we obtained the following classification result which is of interest on its own.

Theorem B. *Let $X$ be a compact complex homogeneous-rational manifold and $L$ an ample line bundle on $X$. Then*

$$\text{(i)}\quad \mathrm{H}^0(X, TX \otimes L^*) = \begin{cases} \mathbb{C} & \textit{for } (X, L) = (\mathbb{CP}_1, \mathcal{O}(2)), \\ \mathbb{C}^n & \textit{for } ,(X, L) = (\mathbb{CP}_n, \mathcal{O}(1)),\ n \geqslant 1 \\ 0 & \textit{otherwise.} \end{cases}$$

(ii) $\mathrm{H}^1(X, TX \otimes L^*) = 0$ *unless $(X, L)$ is one of the entries in Table* 4.

**Table 3:** List of exotic holonomies

| group $G$ | representation $V$ | restrictions/remarks |
|---|---|---|
| $T_{\mathbb{R}} \cdot \mathrm{Spin}(5,5)$ | $\mathbb{R}^{16}$ | |
| $T_{\mathbb{R}} \cdot \mathrm{Spin}(1,9)$ | $\mathbb{R}^{16}$ | |
| $T_{\mathbb{C}} \cdot \mathrm{Spin}(10,\mathbb{C})$ | $\mathbb{C}^{16} \simeq \mathbb{R}^{32}$ | |
| $T_{\mathbb{R}} \cdot E_6^1$ | $\mathbb{R}^{27}$ | |
| $T_{\mathbb{R}} \cdot E_6^4$ | $\mathbb{R}^{27}$ | |
| $T_{\mathbb{C}} \cdot E_6^{\mathbb{C}}$ | $\mathbb{C}^{27} \simeq \mathbb{R}^{54}$ | |
| $T_{\mathbb{R}} \cdot \mathrm{SL}(2,\mathbb{R})$ | $\odot^3\mathbb{R}^2 \simeq \mathbb{R}^4$ | |
| $\mathrm{SL}(2,\mathbb{C})$ | $\odot^3\mathbb{C}^2 \simeq \mathbb{R}^8$ | |
| $\mathbb{C}^* \cdot \mathrm{SL}(2,\mathbb{C})$ | $\odot^3\mathbb{C}^2 \simeq \mathbb{R}^8$ | |
| $\mathbb{R}^* \cdot \mathrm{Sp}(2,\mathbb{R})$ | $\mathbb{R}^4$ | |
| $\mathbb{C}^* \cdot \mathrm{Sp}(2,\mathbb{C})$ | $\mathbb{C}^4 \simeq \mathbb{R}^8$ | |
| $\mathbb{R}^* \cdot \mathrm{SO}(2) \cdot \mathrm{SL}(2,\mathbb{R})$ | $\mathbb{R}^2 \otimes \mathbb{R}^2 \simeq \mathbb{R}^4$ | |
| $\mathbb{C}^* \cdot \mathrm{SU}(2)$ | $\mathbb{C}^2 \simeq \mathbb{R}^4$ | |
| $H_\lambda \cdot \mathrm{SU}(2)$ | $\mathbb{C}^2 \simeq \mathbb{R}^4$ | existence unknown |
| $H_\lambda \cdot \mathrm{SU}(1,1)$ | $\mathbb{C}^2 \simeq \mathbb{R}^4$ | existence unknown |
| $\mathrm{SL}(2,\mathbb{R}) \cdot \mathrm{SO}(p,q)$ | $\mathbb{R}^2 \otimes \mathbb{R}^{p+q} \simeq \mathbb{R}^{2(p+q)}$ | $p+q \geqslant 3$ |
| $\mathrm{Sp}(1) \cdot \mathrm{SO}(n,\mathbb{H})$ | $\mathbb{H}^n \simeq \mathbb{R}^{4n}$ | $n \geqslant 2$ |
| $\mathrm{SL}(2,\mathbb{C}) \cdot \mathrm{SO}(n,\mathbb{C})$ | $\mathbb{C}^2 \otimes \mathbb{C}^n \simeq \mathbb{R}^{4n}$ | $n \geqslant 3$ |
| $E_7^5$ | $\mathbb{R}^{56}$ | |
| $E_7^7$ | $\mathbb{R}^{56}$ | |
| $E_7^{\mathbb{C}}$ | $\mathbb{R}^{112} \simeq \mathbb{C}^{56}$ | |
| $\mathrm{Sp}(3,\mathbb{R})$ | $\mathbb{R}^{14} \subset \Lambda^3\mathbb{R}^6$ | |
| $\mathrm{Sp}(3,\mathbb{C})$ | $\mathbb{R}^{28} \simeq \mathbb{C}^{14} \subset \Lambda^3\mathbb{C}^6$ | |
| $\mathrm{SL}(6,\mathbb{R})$ | $\mathbb{R}^{20} \simeq \Lambda^3\mathbb{R}^6$ | |
| $\mathrm{SU}(1,5)$ | $\mathbb{R}^{20}$ | |
| $\mathrm{SU}(3,3)$ | $\mathbb{R}^{20}$ | |
| $\mathrm{SL}(6,\mathbb{C})$ | $\mathbb{R}^{40} \simeq \Lambda^3\mathbb{C}^6$ | |
| $\mathrm{Spin}(2,10)$ | $\mathbb{R}^{32}$ | |
| $\mathrm{Spin}(6,6)$ | $\mathbb{R}^{32}$ | |
| $\mathrm{Spin}(12,\mathbb{C})$ | $\mathbb{C}^{32} \simeq \mathbb{R}^{64}$ | |

Notation: $T_{\mathbb{F}}$ denotes any connected Lie subgroup of $\mathbb{F}^*$,
$H_\lambda = \left\{ e^{(2\pi i+\lambda)t} \mid t \in \mathbb{R} \right\} \subseteq \mathbb{C}^*, \quad \lambda > 0.$

**Table 4**: The list of all $(X, L)$ with $\mathrm{H}^1(X, TX \otimes L^*) \neq 0$

| $\mathrm{Aut}^0 X$ | $(X, L)$ | $\mathrm{H}^1(X, TX \otimes L^*)$ |
|---|---|---|
| $A_l \quad (l \geqslant 1)$ | 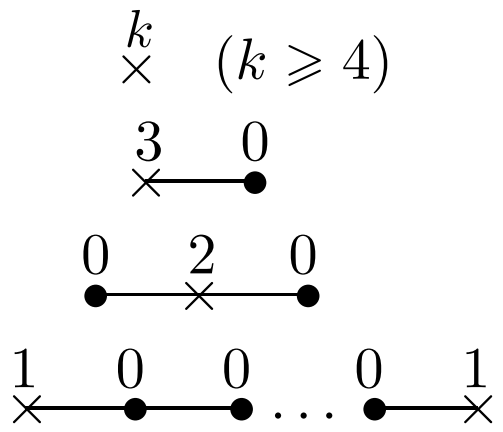 | $\odot^{k-4}\mathbb{C}^2$<br>$\mathbb{C}$<br>$\mathbb{C}$<br>$\mathbb{C}$ |
| $B_l \quad (l \geqslant 3)$ | 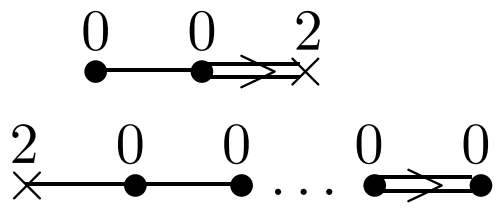 | $\mathbb{C}$<br>$\mathbb{C}$ |
| $C_l \ (l \geqslant 2)$ | 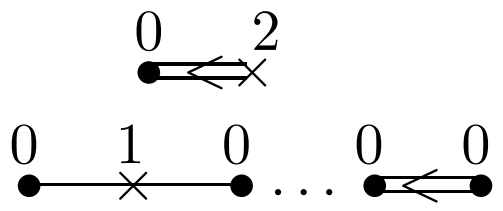 | $\mathbb{C}$<br>$\mathbb{C}$ |
| $D_l \quad (l \geqslant 4)$ | 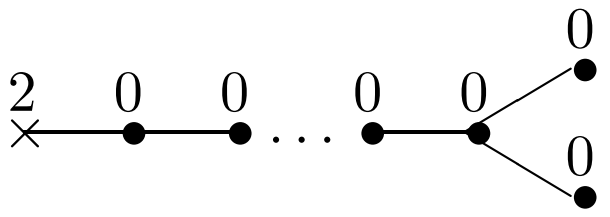 | $\mathbb{C}$ |
| $F_4$ | 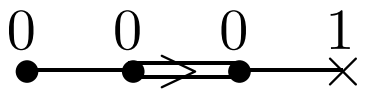 | $\mathbb{C}$ |
| $A_1 \times A_1$ | 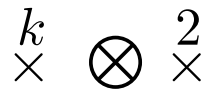 $(k \geqslant 2)$ | $\odot^{k-2}\mathbb{C}^2$ |
| $A_1 \times A_l \ (l \geqslant 1)$ | 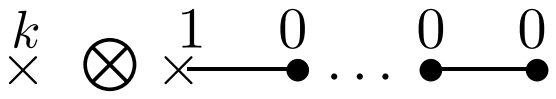 $(k \geqslant 2)$ | $\odot^{k-2}\mathbb{C}^2 \otimes \mathbb{C}^{l+1}$ |
| Notation: $\mathrm{Aut}^0 X$ denotes the universal covering of the component, the identity of the group of all automorphisms of $X$ | | |

This theorem is a "nonvanishing" counterpart to the vanishing theorem of Kobayashi and Ochiai [27] which says that if $X$ is a compact complex rational manifold and $L \to X$ a line bundle such that $\det(TX) \otimes L^*$ is ample, then $\mathrm{H}^i(X, TX \otimes L^*) = 0$ for all $i \geqslant 2$.

Another by-product result of our work on the holonomy problem unveils a striking supersymmetry property of *all* complex symplectic representations and their real forms which occur as holonomies of torsion-free affine connections. With any symplectic Lie algebra $\mathfrak{g} \subset \mathfrak{sp}(V)$, $\dim \mathfrak{g} = m$, $\dim V = 2n$, one may associate naturally an $(m|2n)$-dimensional supermanifold $\mathcal{M}_{\mathfrak{g}} \simeq \mathfrak{g} \oplus \Pi V$, $\Pi$ being the parity change functor [29], which comes equipped with a SUSY (supersymmetry)-structure induced by the canonical projection $\odot^2 V^* \to \mathfrak{g}$.

The structure sheaf $\mathcal{O}_{\mathcal{M}}$ has *supersymmetry* automorphisms which mix purely even and purely odd sectors; in physics literature such automorphisms are often called SUSY-translations.

THEOREM C. *Given an irreducible symplectic representation* $\mathfrak{g} \subset \mathfrak{sp}(V)$, *the structure sheaf of the associated supermanifold* $\mathcal{M}_{\mathfrak{g}}$ *contains* $\mathfrak{g}$*-invariant and* SUSY*-invariant polynomials of order two in odd coordinates and of order at least two in even coordinates if and only if* $\mathfrak{g}$ *occurs as the holonomy of a torsion-free affine connection which is not locally symmetric, i.e.* $\mathfrak{g}$ *is an entry of the following table*:

| Group $G$ | Representation space | Group $G$ | Representation space |
|---|---|---|---|
| $\mathrm{Sp}(n,\mathbb{R})$ | $\mathbb{R}^{2n}$ | $\mathrm{E}_7^5$ | $\mathbb{R}^{56}$ |
| $\mathrm{Sp}(n,\mathbb{C})$ | $\mathbb{C}^{2n}$ | $\mathrm{E}_7^7$ | $\mathbb{R}^{56}$ |
| $\mathrm{SL}(2,\mathbb{R})$ | $\mathbb{R}^4 \simeq \odot^3\mathbb{R}^2$ | $\mathrm{E}_7^{\mathbb{C}}$ | $\mathbb{C}^{56}$ |
| $\mathrm{SL}(2,\mathbb{C})$ | $\mathbb{C}^4 \simeq \odot^3\mathbb{C}^2$ | $\mathrm{Spin}(2,10)$ | $\mathbb{R}^{32}$ |
| $\mathrm{SL}(2,\mathbb{R})\cdot\mathrm{SO}(p,q)$ | $\mathbb{R}^{2(p+q)}$, $p+q \geqslant 3$ | $\mathrm{Spin}(6,6)$ | $\mathbb{R}^{32}$ |
| $\mathrm{SL}(2,\mathbb{C})\cdot\mathrm{SO}(n,\mathbb{C})$ | $\mathbb{C}^{2n}$, $n \geqslant 3$ | $\mathrm{Spin}(12,\mathbb{C})$ | $\mathbb{C}^{32}$ |
| $\mathrm{Sp}(1)\cdot\mathrm{SO}(n,\mathbb{H})$ | $\mathbb{H}^n \simeq \mathbb{R}^{4n}$, $n \geqslant 2$ | $\mathrm{Sp}(3,\mathbb{R})$ | $\mathbb{R}^{14} \subset \Lambda^3\mathbb{R}^6$ |
| $\mathrm{SL}(6,\mathbb{R})$ | $\mathbb{R}^{20} \simeq \Lambda^3\mathbb{R}^6$ | $\mathrm{Sp}(3,\mathbb{C})$ | $\mathbb{C}^{14} \subset \Lambda^3\mathbb{C}^6$ |
| $\mathrm{SU}(1,5)$ | $\mathbb{R}^{20}$ | | |
| $\mathrm{SU}(3,3)$ | $\mathbb{R}^{20}$ | | |
| $\mathrm{SL}(6,\mathbb{C})$ | $\mathbb{C}^{20} \simeq \Lambda^3\mathbb{C}^6$ | | |

A few words about our notation. A short exact sequence of sheaves

$$0 \longrightarrow A_1 \longrightarrow C \longrightarrow A_2 \longrightarrow 0$$

is written sometimes shortly as $C = A_1 + A_2$. One may check that $A_1 + (A_2 + A_3) = (A_1 + A_2) + A_3$ so that the notation $\sum_{i=1}^n A_i$ is non-ambigously defined (and should not be confused with $\oplus_{i=1}^n A_i$).

We also identify vector bundles with the associated locally free sheaves.

## 2. Preliminary facts and results

2.1. *Holonomy groups.* Let $M$ be a smooth connected and simply connected $n$-manifold and $\nabla$ an affine connection on $M$, i.e. a linear connection in the tangent bundle $TM$. Fix a point $x \in M$ and let

$$\mathcal{L}_x = \{\gamma : [0,1] \to M \mid \gamma(0) = \gamma(1) = x\}$$

be the set of piecewise smooth loops based at $x$. For $\gamma \in \mathcal{L}_x$, denote by $P_\gamma : T_xM \longrightarrow T_xM$ the linear automorphism induced by the $\nabla$-parallel translations along $\gamma$.

The (restricted) *holonomy of* $\nabla$ *at* $x \in M$ is defined as a subset $H_x := \{P_\gamma \mid \gamma \in \mathcal{L}_x\} \subset \mathrm{GL}(T_xM)$. Its basic properties are: (i) $H_x$ is a connected Lie subgroup of $\mathrm{GL}(T_xM)$; (ii) if one fixes an isomorphism $i : T_xM \simeq V$, where $V$ is any fixed vector space with $\dim V = \dim M$ (typically, $V = \mathbb{R}^n$), then the conjugacy class of $i(H_x) \subset \mathrm{GL}(V)$ does not depend on the choice of $x \in M$ (see, e.g., [10]). The *holonomy group* of $\nabla$ is defined as any linear subgroup $G \subset \mathrm{GL}(V)$ in the conjugacy class of $i(H_x)$ for some $x \in M$. The Lie algebra $\mathfrak{g} \subset \mathfrak{gl}(V)$ of $G$ is called the *holonomy algebra* of $\nabla$.

2.2. *Berger's criteria.* Let $V$ be a vector space and $\mathfrak{g} \subset \mathfrak{gl}(V)$ a Lie subalgebra. We shall now describe a *necessary condition* for $\mathfrak{g}$ to be the holonomy algebra of a torsion-free connection which is due to Berger [8]. We define the *space of formal curvature maps*

$$K(\mathfrak{g}) := \{R \in \Lambda^2 V^* \otimes \mathfrak{g} \mid R(x,y)z + R(y,z)x + R(z,x)y = 0 \text{ for all } x,y,z \in V\},$$

and the *space of formal curvature derivatives*

$$\begin{aligned} K^1(\mathfrak{g}) \quad := \quad & \{\phi \in V^* \otimes K(\mathfrak{g}) \mid \phi(x)(y,z) + \phi(y)(z,x) \\ & + \phi(z)(x,y) = 0 \text{ for all } x,y,z \in V\}. \end{aligned}$$

We also let $\underline{\mathfrak{g}} := \{R(x,y) \mid R \in K(\mathfrak{g}), x, y \in V\} \subset \mathfrak{g}$.

The equation defining the elements of $K(\mathfrak{g})$ ($K^1(\mathfrak{g})$ respectively) is called the *first* (*second*) *Bianchi identity*. The significance of these spaces is that for a torsion-free connection $\nabla$ on $M$, the curvature $R$ and its covariant derivative $\nabla R$ at a point $x \in M$ lie in $K(\mathfrak{hol}_x)$ and $K^1(\mathfrak{hol}_x)$, respectively, where $\mathfrak{hol}_x \subset \mathfrak{gl}(T_xM)$ is the holonomy algebra of $\nabla$ at $x$.

Note that $K(\mathfrak{g})$ and $K^1(\mathfrak{g})$ are defined by the exact sequences

$$0 \to K(\mathfrak{g}) \longrightarrow \Lambda^2 V^* \otimes \mathfrak{g} \longrightarrow \Lambda^3 V^* \otimes V \tag{1}$$

and

$$0 \to K^1(\mathfrak{g}) \longrightarrow V^* \otimes K(\mathfrak{g}) \longrightarrow \Lambda^3 V^* \otimes V, \tag{2}$$

where in each case, the last map is given by the composition of the natural inclusion and the skew-symmetrization map, i.e. $\Lambda^2 V^* \otimes \mathfrak{g} \hookrightarrow \Lambda^2 V^* \otimes V^* \otimes V \to \Lambda^3 V^* \otimes V$ in the first and $V^* \otimes K(\mathfrak{g}) \hookrightarrow V^* \otimes \Lambda^2 V^* \otimes \mathfrak{g} \to \Lambda^3 V^* \otimes \mathfrak{g}$ in the second case.

*Definition* 2.1. A Lie subalgebra $\mathfrak{g} \subset \mathfrak{gl}(V)$ is called a Berger algebra if

$$1.\ \underline{\mathfrak{g}} = \mathfrak{g}, \qquad \text{and} \qquad 2.\ K^1(\mathfrak{g}) \neq 0.$$

A Lie subgroup $G \subset \mathrm{GL}(V)$ is called a Berger group if its Lie algebra $\mathfrak{g} \subset \mathfrak{gl}(V)$ is a Berger algebra.

In the literature, the two criteria in this definition are usually referred to as *Berger's first and second criterion.* The following simple observation of Berger which is an immediate consequence of the Ambrose-Singer Holonomy Theorem [5] explains the relevance of this definition.

PROPOSITION 2.2 [8]. *Let $G \subset \mathrm{GL}(V)$ be an irreducible Lie subgroup which occurs as the holonomy group of a torsion-free affine connection on some manifold $M$ which is* not *locally symmetric. Then $G$ must be a Berger group, and its Lie algebra $\mathfrak{g} \subset \mathfrak{gl}(V)$ must be a Berger algebra.*

We shall often utilize the following simple result.

LEMMA 2.3. *If $\mathfrak{g} \subset \mathfrak{gl}(V)$ is a Berger algebra, then $K(\mathfrak{g})$ is a nontrivial $\mathfrak{g}$-module.*

*Proof.* Without loss of generality, we may assume that $\dim V > 2$. Suppose $K(\mathfrak{g})$ is a trivial $\mathfrak{g}$-module. Then $K^1(\mathfrak{g}) \subset V^* \otimes K(\mathfrak{g})$ is a submodule and thus, since $V$ is irreducible, we have $K^1(\mathfrak{g}) = V^* \otimes W$ for some subspace $W \subset K(\mathfrak{g})$. Suppose there is a $0 \neq R \in W$. Pick independent elements $x, y, z \in V$ such that $R(x, y) \neq 0$, and define $\phi : V \to W$ such that $\phi(x) = \phi(y) = 0$ and $\phi(z) = R$. Then it follows that $\phi \notin K^1(\mathfrak{g})$ which is a contradiction. Therefore, $W = 0$, i.e. $K^1(\mathfrak{g}) = 0$, and thus $\mathfrak{g}$ is not Berger. □

2.3. *Real Berger algebras.* In this subsection we shall use the following notation: If $W$ is a *complex* vector space, then the Lie algebras of real and complex endomorphisms of $W$ are denoted by $\mathrm{End}_{\mathbb{R}}(W)$ and $\mathrm{End}_{\mathbb{C}}(W)$, respectively.

Let $V$ be a finite dimensional real vector space, and let $\mathfrak{h} \subset \mathrm{End}_{\mathbb{R}}(V)$ be a real Lie subalgebra. We denote their complexifications by $V_{\mathbb{C}} := V \otimes_{\mathbb{R}} \mathbb{C}$ and $\mathfrak{h}_{\mathbb{C}} := \mathfrak{h} \otimes_{\mathbb{R}} \mathbb{C}$. Then obviously, $\mathfrak{h}_{\mathbb{C}} \subset \mathrm{End}_{\mathbb{C}}(V_{\mathbb{C}})$, and by complexifying the exact sequences (1) and (2), we obtain

$$K(\mathfrak{h}_{\mathbb{C}}) = K(\mathfrak{h}) \otimes_{\mathbb{R}} \mathbb{C} \quad \text{and} \quad K^1(\mathfrak{h}_{\mathbb{C}}) = K^1(\mathfrak{h}) \otimes_{\mathbb{R}} \mathbb{C}.$$

In particular, $\mathfrak{h} \subset \mathrm{End}_{\mathbb{R}}(V)$ is a (symmetric, respectively, nonsymmetric) Berger algebra if and only if $\mathfrak{h}_{\mathbb{C}} \subset \mathrm{End}_{\mathbb{C}}(V_{\mathbb{C}})$ is.

Let us now assume that $\mathfrak{h} \subset \mathrm{End}_{\mathbb{R}}(V)$ is *irreducible.* Then there are two cases to be distinguished.

First, suppose that $\mathfrak{h}$ is of *real type*, i.e. there is no complex structure on $V$ which commutes with the elements of $\mathfrak{h}$. This happens if and only if $\mathfrak{h}_{\mathbb{C}} \subset \mathrm{End}_{\mathbb{C}}(V_{\mathbb{C}})$ is also irreducible.

Second, suppose that $\mathfrak{h}$ is *not* of real type, i.e. there is a complex structure $J$ on $V$ which commutes with the elements of $\mathfrak{h}$. That is, $\mathfrak{h} \subset \mathrm{End}_{\mathbb{C}}(V)$ with respect to this complex structure $J$. In this case, $V_{\mathbb{C}} = W \oplus \overline{W}$ decomposes into two irreducible $\mathfrak{h}_{\mathbb{C}}$-submodules of equal dimension given by

$$W = \{x + iJx \mid x \in V\} \text{ and } \overline{W} = \{x - iJx \mid x \in V\}.$$

Let $\mathfrak{h}_1 := \{A \in \mathfrak{h} \mid JA \in \mathfrak{h}\}$. Then $\mathfrak{h}_1 \lhd \mathfrak{h}$, and $J$ induces a complex Lie algebra structure on $\mathfrak{h}_1$; $(\mathfrak{h}_1)_{\mathbb{C}}$ can be written as the direct sum of complex Lie algebras $(\mathfrak{h}_1)_{\mathbb{C}} = \mathfrak{h}_1^+ \oplus \mathfrak{h}_1^-$ with

$$\mathfrak{h}_1^+ = \{A + iJA \mid A \in \mathfrak{h}_1\} \text{ and } \mathfrak{h}_1^- = \{A - iJA \mid A \in \mathfrak{h}_1\}.$$

Let $R \in K(\mathfrak{h}_{\mathbb{C}})$. Then for $u, v \in W$ and $\overline{w} \in \overline{W}$ the first Bianchi identity implies that $R(u,v)\overline{w} = 0$. Since this is true for all $\overline{w} \in \overline{W}$, it follows that $R(u,v) \in \mathfrak{h}_1^+$. Likewise, for $\overline{u}, \overline{v} \in \overline{W}$, we have $R(\overline{u}, \overline{v}) \in \mathfrak{h}_1^-$.

Next, for any $R \in K(\mathfrak{h}_{\mathbb{C}})$ the first Bianchi identity also implies that $R(\overline{u}, v)w = R(\overline{u}, w)v$ for all $\overline{u} \in \overline{W}$, $v, w \in W$. Thus, we have a map

$$\overline{W} \longrightarrow (\mathfrak{h}_{\mathbb{C}}|_W)^{(1)}, \quad \overline{u} \longmapsto R(\overline{u}, _).$$

If $(\mathfrak{h}_{\mathbb{C}}|_W)^{(1)} = 0$ then this implies that $R(W, \overline{W}) = 0$, and hence by the above, $R(V_{\mathbb{C}}, V_{\mathbb{C}}) \subset \mathfrak{h}_1^+ \oplus \mathfrak{h}_1^- = (\mathfrak{h}_1)_{\mathbb{C}}$ for all $R \in K(\mathfrak{h}_{\mathbb{C}})$; that is, $\underline{\mathfrak{h}}_{\mathbb{C}} \subset (\mathfrak{h}_1)_{\mathbb{C}}$. Hence $\mathfrak{h}_{\mathbb{C}}$ is not Berger unless $\mathfrak{h}_1 = \mathfrak{h}$; i.e. $\mathfrak{h}$ is a complex Lie algebra which acts irreducibly on the complex vector space $V$.

We define a map $\imath : \mathfrak{h}_{\mathbb{C}} \to \mathrm{End}_{\mathbb{C}}(V)$ by

$$\imath(A + iB) := A + JB. \tag{3}$$

In fact, it is easy to see that $\imath(\mathfrak{h}_{\mathbb{C}}) \subset \mathrm{End}_{\mathbb{C}}(V)$ is congruent to $(\mathfrak{h}_{\mathbb{C}})|_W \subset \mathrm{End}_{\mathbb{C}}(W)$, and hence $(\mathfrak{h}_{\mathbb{C}}|_W)^{(1)} = 0$ if and only if $(\imath(\mathfrak{h}_{\mathbb{C}}))^{(1)} = 0$. Thus, we obtain the following.

Proposition 2.4. *Let $V$ be a finite dimensional real vector space, and let $\mathfrak{h} \subset \mathrm{End}_{\mathbb{R}}(V)$ be an irreducible real subalgebra with complexification $\mathfrak{h}_{\mathbb{C}} \subset \mathrm{End}_{\mathbb{C}}(V_{\mathbb{C}})$.*

1. *If $\mathfrak{h}$ is of real type, i.e. if there is no complex structure on $V$ which commutes with the elements of $\mathfrak{h}$, then $\mathfrak{h}$ is a Berger algebra if and only if $\mathfrak{h}_{\mathbb{C}} \subset \mathrm{End}_{\mathbb{C}}(V_{\mathbb{C}})$ is an irreducible Berger algebra.*

2. *If $\mathfrak{h}$ is not of real type, i.e. if there is a complex structure $J$ on $V$ which commutes with the elements of $\mathfrak{h}$, and if the subalgebra $\imath(\mathfrak{h}_{\mathbb{C}}) \subset \mathrm{End}_{\mathbb{C}}(V)$ given by (3) satisfies $(\imath(\mathfrak{h}_{\mathbb{C}}))^{(1)} = 0$, then $\mathfrak{h}$ is a Berger algebra if and only if $J\mathfrak{h} = \mathfrak{h}$ and $\mathfrak{h} \subset \mathrm{End}_{\mathbb{C}}(V)$ is a complex irreducible Lie subalgebra.*

Thus, in order to classify all Berger algebras we need to classify all irreducible *complex* Berger subalgebras $\mathfrak{h}_{\mathbb{C}} \subset \mathrm{End}_{\mathbb{C}}(V_{\mathbb{C}})$, add all their real forms of real type, and finally, to investigate the real algebras $\mathfrak{h} \subset \mathrm{End}_{\mathbb{C}}(V)$ for which $\imath(\mathfrak{h}_{\mathbb{C}}) \subset \mathrm{End}_{\mathbb{C}}(V)$ is one of the entries of Table 5. The latter task has been completed by Bryant [16]; hence, we shall mainly concern ourselves with the investigation of complex Berger algebras.

2.4. *Conformal representations.* In this subsection, we shall prove the following result.

PROPOSITION 2.5. *Let $\mathfrak{g} \subset \mathfrak{so}(V)$ be a proper irreducible Lie subalgebra on some finite dimensional euclidean vector space $V$ with $\dim V \geqslant 5$. If $\mathrm{Hom}_{\mathfrak{g}}(\odot^2 \mathrm{Ad}(\mathfrak{g}), \Lambda^2 V) = 0$, i.e. if the $\mathfrak{g}$-modules $\odot^2 \mathrm{Ad}(\mathfrak{g})$ and $\Lambda^2 V$ have no irreducible summand in common, then $K(\mathfrak{g} \oplus \mathbb{C}\,\mathrm{Id}) = K(\mathfrak{g})$. In particular, $\mathfrak{g} \oplus \mathbb{C}\,\mathrm{Id}$ is not a Berger algebra.*

*Proof.* We identify $\mathfrak{so}(V)$ and $\Lambda^2 V$ by the equation $(x \wedge y) \cdot z := (x, z)y - (y, z)x$. It is well known that $K(\mathfrak{so}(V) \oplus \mathbb{C}\,\mathrm{Id}) = K(\mathfrak{so}(V)) \oplus K^c(V)$, where the subspace $K^c(V)$ is isomorphic to $\mathfrak{so}(V)$ via

$$R_A := (A, _)\,\mathrm{Id} + \mathrm{ad}(A) \quad \text{for all } A \in \mathfrak{so}(V).$$

Here, ( , ) denotes a suitable multiple of the Killing form on $\mathfrak{g}$. Moreover, it is known that $K(\mathfrak{so}(V)) \subset \odot^2(\mathrm{Ad}(\mathfrak{so}(V)))$, and therefore, $K(\mathfrak{g}) \subset \odot^2(\mathrm{Ad}(\mathfrak{g}))$.

By hypothesis, $K(\mathfrak{g})$ and $K^c(V) \cong \Lambda^2 V$ have no irreducible summand in common, and hence

$$K(\mathfrak{g} \oplus \mathbb{C}\,\mathrm{Id}) = K(\mathfrak{g}) \oplus [K^c(V) \cap K(\mathfrak{g} \oplus \mathbb{C}\,\mathrm{Id})].$$

Let

$$\mathfrak{h} := \{A \in \mathfrak{so}(V) \mid R_A \in K^c(V) \cap K(\mathfrak{g} \oplus \mathbb{C}\,\mathrm{Id})\} = \{A \in \mathfrak{so}(V) \mid [A, \mathfrak{so}(V)] \subset \mathfrak{g}\}.$$

It is then straightforward to show that $\mathfrak{h} \trianglelefteq \mathfrak{so}(V)$. But $\dim(V) \geqslant 5$, thus $\mathfrak{so}(V)$ is simple, and obviously, $\mathfrak{h} \neq \mathfrak{so}(V)$ if $\mathfrak{g}$ is a proper subalgebra. Thus, $\mathfrak{h} = 0$, and this completes the proof. □

COROLLARY 2.6. *Let $\mathfrak{g} \subset \mathfrak{gl}(V)$ be one of the following subalgebras:*

$$\mathfrak{g}_2^{\mathbb{C}} \subset \mathfrak{gl}(\mathbb{C}^7), \quad \mathfrak{spin}(7, \mathbb{C}) \subset \mathfrak{gl}(\mathbb{C}^8), \quad \mathfrak{spin}(9, \mathbb{C}) \subset \mathfrak{gl}(\mathbb{C}^{16}).$$

*Then $\mathfrak{g} \oplus \mathbb{C}\,\mathrm{Id} \subset \mathfrak{gl}(V)$ is not a Berger algebra.*

*Proof.* It is well-known that all these representations are metric. The decomposition of the $\mathfrak{g}$-modules $\odot^2(\mathrm{Ad}(\mathfrak{g}))$ and $\Lambda^2 V$ is straightforward and yields that Proposition 2.5 applies in each case. □

2.5. *Review of representation theory* [7], [24]. Let $\mathfrak{g}$ be a semisimple complex Lie algebra and $G$ the associated simply connected Lie group. Fix a maximally Abelian self-normalizing subalgebra $\mathfrak{t} \subset \mathfrak{g}$ (any two such subalgebras, called Cartan subalgebras, are conjugate under the adjoint action of $G$). If $\rho : \mathfrak{g} \to \mathfrak{gl}(V)$ is a representation of $\mathfrak{g}$ in a complex vector space $V$, then with any $\omega \in \mathfrak{t}^* \equiv \mathrm{Hom}_{\mathbb{C}}(\mathfrak{t}, \mathbb{C})$ one may associate the *weight space* of $V$ by $V_\omega = \{v \in V : \rho(h)v = \omega(h)v \text{ for all } h \in \mathfrak{t}\}$. An element $\omega \in \mathfrak{t}^*$ is called a weight of $V$ if $V_\omega \neq 0$.

In the particular case when $V = \mathfrak{g}$ and $\rho : \mathfrak{g} \to \mathfrak{gl}(\mathfrak{g})$ is the adjoint representation of $\mathfrak{g}$ on itself, the nonzero weights of $\mathfrak{g}$ are called the *roots* of $\mathfrak{g}$. Thus

$$\mathfrak{g} = \mathfrak{t} \oplus \sum_{\alpha \in \Phi} \mathfrak{g}_\alpha$$

where $\Phi$ is the set of all roots of $\mathfrak{g}$ and all sums are direct. A subset $\Delta = \{\alpha_1, \ldots, \alpha_r\} \subset \Phi$ with the property that every $\omega \in \Phi$ may be expressed as a linear combination $\omega = \sum_{i=1}^r a_i \alpha_i$ with all $a_i$ being nonnegative or all nonpositive integers is called a *system of simple roots* of $\mathfrak{g}$. Such a $\Delta$ exists, and any two such sets are conjugate under the adjoint action of $G$. If we fix $\Delta$, then $\Phi = \Phi^+ \cup \Phi^-$, where $\Phi^+ = \{\omega \in \Phi : \omega = \sum_{i=1}^r a_i \alpha_i \text{ with } a_i \geqslant 0\}$ is the set of *positive roots* and $\Phi^- = \{\omega \in \Phi : \omega = \sum_{i=1}^r a_i \alpha_i \text{ with } a_i \leqslant 0\}$ is the set of *negative roots* (both with respect to $\Delta$).

For any root $\alpha \in \Phi^+$ there is a unique element $H_\alpha$ in $[\mathfrak{g}_\alpha, \mathfrak{g}_{-\alpha}] \subset \mathfrak{t}$ such that $\alpha(H_\alpha) = 2$. If $\Delta = \{\alpha_1, \ldots, \alpha_r\}$ is the set of simple roots, then the associated set $\{H_{\alpha_1}, \ldots, H_{\alpha_r}\}$ forms a basis of $\mathfrak{g}$. Its dual basis $\{\omega_{\alpha_1}, \ldots, \omega_{\alpha_r}\}$ of $\mathfrak{t}^*$ is called the set of *fundamental weights*. One may use it to define the following three important subsets of $\mathfrak{t}^*$: the set of (*integral*) *weights* $\Lambda = \{\lambda \in \mathfrak{t}^* : \omega = \sum_{i=1}^r \lambda_i \omega_i \text{ with } \lambda_i \in \mathbb{Z}\}$; the set of *dominant weights* $\Lambda^+ = \{\lambda \in \Lambda : \lambda = \sum_{i=1}^r \lambda_i \omega_i \text{ with } \lambda_i \geqslant 0\}$; and the set of *strongly dominant weights* $\Lambda^{++} = \{\lambda \in \Lambda^+ : \lambda = \sum_{i=1}^r \lambda_i \omega_i \text{ with } \lambda_i > 0\}$. Note that $\lambda_i = \lambda(H_{\alpha_i})$. The minimal integral element $\omega_1 + \omega_2 + \ldots + \omega_r$ in $\Lambda^{++}$ is denoted by $\eta$. Any integral weight $\lambda$ of $\mathfrak{g}$ can be graphically represented by inscribing the integer $\lambda_i$ over the $i^{\text{th}}$ node of the Dynkin diagram for $\mathfrak{g}$. For example, the fundamental weight $\omega_1$ of $\mathfrak{sl}(3, \mathbb{C})$ is $\overset{1}{\bullet}\!\!-\!\!\!-\!\!\!-\!\!\overset{0}{\bullet}$ .

Let $\lambda \in \Lambda$ be an integral weight. It is called *singular* if $\lambda(H_\alpha) = 0$ for some $\alpha \in \Phi^+$, and *regular* otherwise. The *index* of $\lambda$ is defined to be the number of positive roots $\alpha$ for which $\lambda(H_\alpha) < 0$ holds; it is denoted by $\mathrm{ind}(\lambda)$.

If $\rho : \mathfrak{g} \to \mathfrak{gl}(V)$ is an irreducible representation of $\mathfrak{g}$, then there exists a unique weight $\omega(V) \in \Lambda^+$ of $V$, called the *highest weight* of $V$ (relative to $\mathfrak{t}$ and $\Delta$) such that $\dim V_\omega = 1$ and $\rho(\mathfrak{g}_\alpha) V_\omega = 0$ for all $\alpha \in \Delta$. This establishes a one-to-one correspondence, $V \Leftrightarrow \omega(V)$, between finite-dimensional irreducible $\mathfrak{g}$-modules and dominant weights, and allows us to use the graphical description of $\omega(V)$ to represent $\rho : \mathfrak{g} \to \mathfrak{gl}(V)$. For example, the standard representation of $\mathfrak{sl}(3, C)$ in $\mathbb{C}^3$ is denoted by 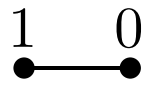 .

If $\mathfrak{g}$ is simple, then the adjoint representation $\rho : \mathfrak{g} \to \mathfrak{gl}(\mathfrak{g})$ is irreducible. The associated highest weight of $V = \mathfrak{g}$ is a root $\mu \in \Lambda^+$ which is called the *maximal root* of $\mathfrak{g}$. The following is the list of all maximal roots [7], [24]:

(4)

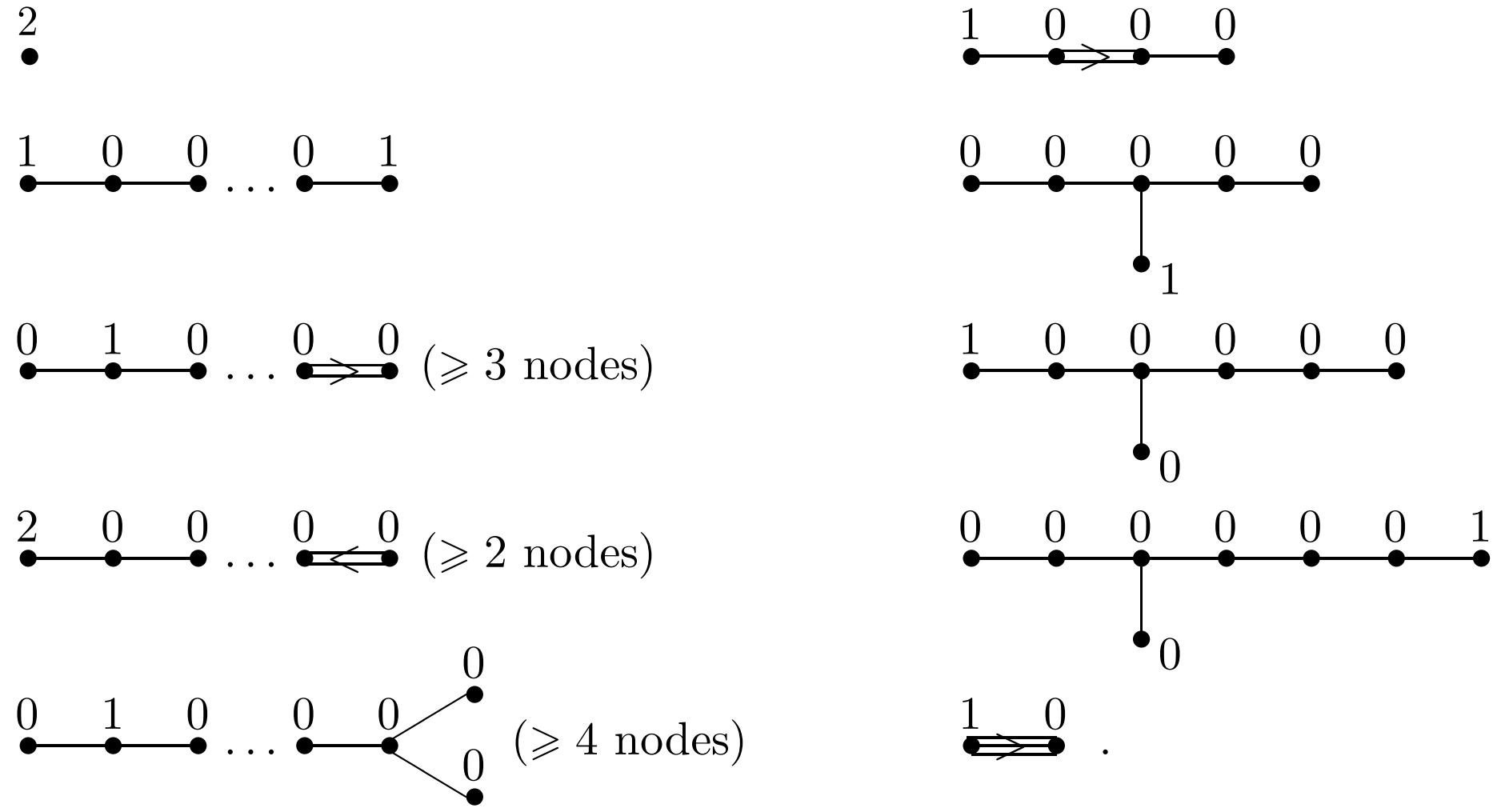

For any simple root $\alpha_i \in \Delta$, denote by $\sigma_i$ the reflection in the hyperplane perpendicular to $\alpha_i$. The *Weyl group* $W$ of $\mathfrak{g}$ is the group generated by all the simple reflections $\sigma_i$. If $\mathfrak{g}$ is simple then $W$ acts irreducibly on $\mathfrak{t}$, and moreover, there is a $W$-invariant inner product $( \ , \ )$ on $\mathfrak{t}$ (the so-called *Killing form*). Given a weight $\lambda \in \Lambda$ and a root $\alpha$ of $\mathfrak{g}$, the *Cartan number*, which is known to be an integer, is given by

$$\langle \lambda, \alpha \rangle = \frac{2(\lambda, \alpha)}{(\alpha, \alpha)} \in \mathbb{Z}.$$

Note that $\langle \ , \ \rangle$ is linear in the first entry only. Its significance is the following. If $\lambda$ occurs as the weight of an irreducible representation of $\mathfrak{g}$ and $\langle \lambda, \alpha \rangle > 0$ ($\langle \lambda, \alpha \rangle < 0$, respectively) then $\lambda - k\alpha$ ($\lambda + k\alpha$, respectively) is also a weight of that representation for $k = 1, \dots |\langle \lambda, \alpha \rangle|$. For the adjoint representation, it is well-known that $|\langle \alpha, \beta \rangle| \leqslant 3$, and $|\langle \alpha, \beta \rangle| = 3$ occurs if and only if $\mathfrak{g}$ contains $\mathfrak{g}_2$ as a direct summand.

The action of the simple reflection $\sigma_i$ on a weight $\lambda \in \Lambda$ can be described by the following rule [7]: to compute $\sigma_i(\lambda)$, let $c = \lambda(H_{\alpha_i})$ be the coefficient of the node associated to $\alpha_i$; add $c$ to the adjacent coefficients, with multiplicity if there is a multiple edge directed towards the adjacent node, and then replace

$c$ by $-c$. For example

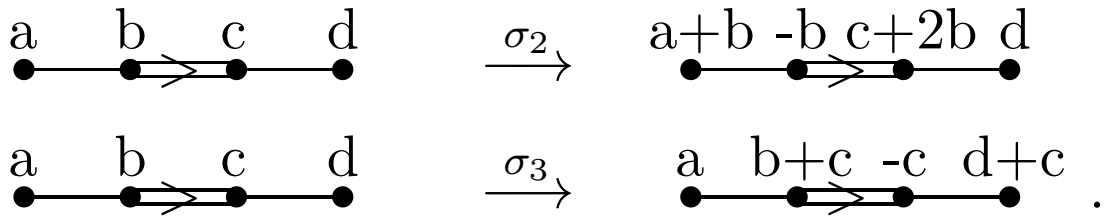

For any $w \in W$, there exists a minimal integer $l(w)$ such that $w$ can be expressed as a composition of $l(w)$ simple reflections. This integer is called the *length* of $w$.

2.6. *Homogeneous manifolds and vector bundles.* A maximal solvable subalgebra of a semisimple Lie algebra $\mathfrak{g}$ is called a *Borel subalgebra.* A subalgebra $\mathfrak{p} \subset \mathfrak{g}$ is called *parabolic* if it contains a Borel subalgebra. Every Borel subalgebra is $G$-conjugate to the standard one,

$$\mathfrak{b} := \mathfrak{t} \oplus \mathfrak{n}$$

where $\mathfrak{n} := \sum_{\alpha\in\Phi^+} \mathfrak{g}_\alpha$. There is a standard form $\mathfrak{p}$ for any parabolic subalgebra as well. Let $\Delta_\mathfrak{p}$ be a subset of $\Delta$ and let $\Phi^+_\mathfrak{p} = \text{span}\,\{\Delta_\mathfrak{p}\} \cap \Phi^+$. Then

$$\mathfrak{p} = \mathfrak{t} + \mathfrak{n} + \sum_{\alpha\in\Phi^+_\mathfrak{p}} \mathfrak{g}_{-\alpha}$$

is the standard parabolic subalgebra of $\mathfrak{g}$. A useful notation for a standard parabolic $\mathfrak{p} \subset \mathfrak{g}$ (and for the associated subgroup $P \subset G$) is to cross all nodes in the Dynkin diagram for $\mathfrak{g}$ which correspond to simple roots of $\mathfrak{g}$ in $\Delta \setminus \Delta_\mathfrak{p}$.

It is well-known that any compact complex homogeneous-rational manifold $X$ is isomorphic to the quotient space $G/P$, where $G$ is a simply connected Lie group and $P \subset G$ is a parabolic subgroup. It is then very useful to denote $X$ by the same Dynkin diagram as $\mathfrak{p}$, the Lie algebra of $P$. For example, the odd dimensional quadric $Q_{2n-1}$ is denoted by 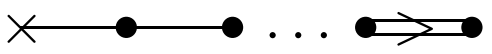 .

The number of crossed nodes in the Dynkin diagram for $X$ is called the *rank* of $X$ and is denoted by rank$X$. This number is independent of the representation of $X$ as a quotient $G/P$.

A vector bundle $E \to X = G/P$ is called *$G$-homogeneous* if there is a holomorphic representation $\rho : P \to \mathrm{GL}(V)$ such that $E = G \times_\rho V$; i.e., $E$ is the quotient $G \times V/P$, where every $p \in P$ acts on $G \times V$ as follows

$$\begin{array}{ccc} G \times V & \longrightarrow & G \times V \\ (g, v) & \longrightarrow & (g \cdot p, \rho(p^{-1})v). \end{array}$$

If $\rho : P \to \mathrm{GL}(V)$ is irreducible, then $E$ is said to be *irreducible* as well.

The finite-dimensional irreducible representations of $P$ are in one-to-one correspondence with integral weights $\lambda \in \Lambda$ whose Dynkin diagram has non-negative coefficients over the uncrossed nodes for $\mathfrak{p}$. A useful notation for an irreducible homogeneous vector bundle $E \to X$ is to combine the Dynkin diagram for the associated integral weight $\lambda$ with the Dynkin diagram for $\mathfrak{p}$ into one picture. For example, if $X = \times\!\!-\!\!\!-\!\!\bullet$ is the projective plane $\mathbb{CP}_2$, then $\mathcal{O}(-1) = \overset{-1}{\times}\!\!-\!\!\!-\!\!\overset{0}{\bullet}$ and $TX = \overset{1}{\times}\!\!-\!\!\!-\!\!\overset{1}{\bullet}$ .

The cohomology ring $\mathrm{H}^*(X, E)$ of an irreducible homogeneous vector bundle $E \to X$ with integral weight $\lambda \in \Lambda$ can be computed, according to Bott [11], as follows:

(i) If $\lambda + \eta$ is singular, then $\mathrm{H}^*(X, E) = 0$;

(ii) If $\lambda + \eta$ is regular and if $\mathrm{ind}(\lambda + \eta) = p$, then there is a unique element $\sigma_\lambda$ (of length $p$) in the Weyl group of $\Phi$ such that $\sigma_\lambda(\lambda + \eta) \in \Lambda^{++}$. Then $\mathrm{H}^*(X, E) = \mathrm{H}^p(X, E)$ and $\mathrm{H}^p(X, E)$ is an irreducible $\mathfrak{g}$-module whose highest weight is $\sigma_\lambda(\lambda + \eta) - \eta$.

2.7. *The tangent bundle of a compact homogeneous manifold.* Let $X = G/P$ be a compact complex homogeneous-rational manifold. The associated number $|\Delta \backslash \Delta_{\mathfrak{p}}|$ of crossed nodes does not depend on the representation of $X$ as the quotient $G/P$ of a semisimple complex Lie group $G$ by a parabolic subgroup $P$; i.e., if $G'$ is another semisimple complex Lie group acting transitively on $X$, $X = G'/P'$, then $|\Delta \setminus \Delta_{\mathfrak{p}}| = |\Delta' \setminus \Delta'_{\mathfrak{p}'}|$; this number is called the *rank* of $X$.

The tangent bundle $TX$ is the homogeneous vector bundle on $X$ associated with the natural representation of $\mathfrak{p}$ on $\mathfrak{g}/\mathfrak{p}$. Since $\mathfrak{p}$ is never semisimple, $TX$ is *not*, in general, decomposable into a direct sum of irreducible homogeneous vector bundles. Nevertheless, there is a simple recipe to exhibit the weight structure of $TX$ which works a follows. First define the map

$$\pi : \Phi^+ \setminus \Phi^+_{\mathfrak{p}} \longrightarrow \mathbb{Z}^{\mathrm{rank} X}$$

by taking coefficients of positive roots from $\Phi^+ \setminus \Phi^+_{\mathfrak{p}}$ over crossed nodes. Let $\Pi = \mathrm{Im}\pi$. Define the height of a positive root $\alpha \in \Phi^+$ to be the sum of all coefficients of $\alpha$ in the basis $\Delta$. Define the section of $\pi$

$$\begin{array}{rccc} \alpha : & \Pi & \longrightarrow & \Phi^+ \setminus \Phi^+_{\mathfrak{p}} \\ & n & \longrightarrow & \alpha(n) \end{array}$$

by the condition that $\alpha(n)$ is the root in the fibre $\pi^{-1}(n)$ which has the maximal height. By considering root strings, one easily obtains the following decomposition [7]:

$$TX = \sum_{i \geqslant 1} \left[ \bigoplus_{|n|=i} J(\alpha(n)) \right]. \tag{5}$$

where $J(\lambda)$ denotes the irreducible homogeneous bundle on $X$ with the highest weight $\lambda$ and the meaning of $\sum$ is as explained at the end of the introduction.

If $\mathfrak{q} \subset \mathfrak{p}$ are parabolic subalgebras of a semisimple Lie algebra $\mathfrak{g}$, then there is a natural fibration

$$X = G/Q \longrightarrow X_0 = G/P$$

whose fibre $X_v$ can be identified by the following rule [7]: delete from the Dynkin diagram for $\mathfrak{q}$ all crossed nodes (and incident) edges shared with $\mathfrak{p}$ and then delete all connected components with no crossed through nodes.

For example, the fibration ×—×—• $\longrightarrow$ ×—•—• has the fibre $\mathbb{P}^2 =$ ×—• .

If $\nu : X \to X_0$ is the fibration of the above type, then

$$TX = T_v + \nu^*(TX_0) \tag{6}$$

where $T_v$ is the homogeneous vector bundle on $X$ whose restriction to the fibre $X_v$ of $\nu$ is isomorphic to $TX_v$. Any ample line bundle $L \to X$ can be uniquely factored as $L_v \otimes \nu^*(L_0)$, where $L_0$ is an ample line bundle on $X_0$ and the line bundle $L_v \to X$ restricted to the fibers of $\nu$ is an ample line bundle on $X_v$.

Then

$$\begin{aligned} \odot^2 TX \otimes L^{*2} \;=\;& \odot^2 T_v \otimes L^{*2} + (T_v \otimes L_v^{*2}) \otimes \nu^*(TX_0 \otimes L_0^{*2}) \\ &+ L_v^{*2} \otimes \nu^*(\odot^2 TX_0 \otimes L_0^{*2}) \end{aligned} \tag{7}$$

and

$$\begin{aligned} \odot^3 TX \otimes L^{*2} \;=\;& \odot^3 T_v \otimes L^{*2} \;+\; (\odot^2 T_v \otimes L_v^{*2}) \otimes \nu^*(TX_0 \otimes L_0^{*2}) \\ &+ (T_v \otimes L_v^{*2}) \otimes \nu^*(\odot^2 TX_0 \otimes L_0^{*2}) \\ &+ L_v^{*2} \otimes \nu^*(\odot^3 TX_0 \otimes L_0^{*2}). \end{aligned} \tag{8}$$

If $\mathrm{rank} X \geqslant 2$, then there is always a fibration $X \to X_0$ of the above type whose fibre $X_v$ is either a projective space $\mathbb{P}_n$ or a Grassmanian $G(k; \mathbb{C}^{n+1})$, $k \geqslant 2$. We shall find very useful the following simple statement.

Lemma 2.7. (i) *If* $X_v = \mathbb{P}_n$ *with* $n \geqslant 2$, *then*

$$\mathrm{H}^1(X, \odot^2 TX \otimes L^{*2}) \;\subset\; \mathrm{H}^1(X, \odot^2 T_v \otimes L^{*2}), \tag{9}$$

$$\begin{aligned} \mathrm{H}^1(X, \odot^3 TX \otimes L^{*2}) \;\subset\;& \mathrm{H}^1\left(X, \odot^3 T_v \otimes L^{*2}\right) \\ &\oplus \mathrm{H}^1\left(X_0, \nu^0_*(\odot^2 T_v \otimes L_v^{*2}) \otimes TX_0 \otimes L_0^{*2}\right). \end{aligned} \tag{10}$$

*In particular, if* $L_v|_{X_v} = \mathcal{O}(k)$ *with* $k \geqslant 2$, *then*

$$\mathrm{H}^1(X, \odot^2 TX \otimes L^{*2}) = \mathrm{H}^1(X, \odot^3 TX \otimes L^{*2}) = 0.$$

(ii) *If* $X_v = G(k, \mathbb{C}^{n+1})$ *with* $(k, n) \neq (3, 6)$ *or* $(4, 8)$, *then*

$$\mathrm{H}^1(X, \odot^2 TX \otimes L^{*2}) = \mathrm{H}^1(X, \odot^3 TX \otimes L^{*2}) = 0. \tag{11}$$

*Proof.* The statements easily follow from (7), (8), the Leray spectral sequence and the fact that for $X_v = \mathbb{P}_n$ with $n \geqslant 2$ one has, for any ample line bundle $L_v \to X_v$,

$$\mathrm{H}^1(X_v, \odot^2 TX_v \otimes L_v^{*2}) = \mathrm{H}^i(X_v, L_v^{*2}) = \mathrm{H}^i(X_v, TX_v \otimes L_v^{*2}) = 0, \quad i = 0, 1,$$

while for $(X_v, L_v) = (\mathbb{P}_n, \mathcal{O}(k))$, $n, k \geqslant 2$ and for $X_v = G(k, \mathbb{C}^{n+1})$ with $(k, n) \neq (3, 6)$ or $(4, 8)$ one has, in addition to the above, the vanishing of the following cohomology groups:

$$\mathrm{H}^0(X_v, \odot^2 TX_v \otimes L_v^{*2}) = \mathrm{H}^i(X_v, \odot^3 TX_v \otimes L_v^{*2}) = 0, \quad i = 0, 1. \qquad \square$$

## 3. Spencer cohomology and complex homogeneous manifolds

3.1. *Spencer cohomology.* Let $V$ be a vector space and $\mathfrak{g}$ a Lie subalgebra of $\mathfrak{gl}(V) := V \otimes V^*$. Define recursively the $\mathfrak{g}$-modules

$$\begin{aligned} \mathfrak{g}^{(-1)} &= V \\ \mathfrak{g}^{(0)} &= \mathfrak{g} \\ \mathfrak{g}^{(k)} &= [\mathfrak{g}^{(k-1)} \otimes V^*] \cap [V \otimes \odot^{k+1} V^*] \ , \ k = 1, 2, \dots , \end{aligned}$$

and define the map

$$\partial : \mathfrak{g}^{(k)} \otimes \Lambda^{l-1} V^* \longrightarrow \mathfrak{g}^{(k-1)} \otimes \Lambda^l V^*$$

as the antisymmetrisation over the last $l$ indices. Here and elsewhere the symbols $\odot^k$ and $\Lambda^k$ stand for $k^{\text{th}}$ order symmetric and antisymmetric powers, respectively.

Since $\partial^2 = 0$, there is a complex

$$\mathfrak{g}^{(k)} \otimes \Lambda^{l-1} V^* \overset{\partial}{\longrightarrow} \mathfrak{g}^{(k-1)} \otimes \Lambda^l V^* \overset{\partial}{\longrightarrow} \mathfrak{g}^{(k-2)} \otimes \Lambda^{l+1} V^*$$

whose cohomology at the center term is denoted by $\mathrm{H}^{k,l}(\mathfrak{g})$ and is called the $(k, l)$ *Spencer cohomology group.* In particular,

$$\begin{aligned} \mathrm{H}^{k,1}(\mathfrak{g}) &= 0, \\ \mathrm{H}^{k,2}(\mathfrak{g}) &= \frac{\mathrm{Ker} : \mathfrak{g}^{(k-1)} \otimes \Lambda^2 V^* \overset{\partial}{\to} \mathfrak{g}^{(k-2)} \otimes \Lambda^3 V^*}{\mathrm{Image} : \mathfrak{g}^{(k)} \otimes V^* \overset{\partial}{\to} \mathfrak{g}^{(k-1)} \otimes \Lambda^2 V^*}. \end{aligned} \tag{12}$$

The $\mathfrak{g}$-module $\mathrm{H}^{k,2}(\mathfrak{g})$ has the following geometric meaning: if $G$ is a matrix Lie group whose Lie algebra is $\mathfrak{g}$ and $\mathcal{G} \to M$ is a $G$-structure on a manifold $M$ which is infinitesimally flat to $k^{\text{th}}$ order, then the obstruction for $\mathcal{G}$ to be infinitesimally flat to $(k+1)^{\text{th}}$ order is given by a section of the associated vector bundle $\mathcal{G} \times_G \mathrm{H}^{k,2}(\mathfrak{g})$.

The $\mathfrak{g}$-module $\mathfrak{g}^{(1)}$ has a clear geometric interpretation as well. If a $G$-structure $\mathcal{G} \to M$ is infinitesimally flat to 1st order (which is equivalent to saying that $\mathcal{G}$ admits a *torsion-free* affine connection), then the set of all torsion-free affine connections in $\mathcal{G}$ is an affine space modelled on the vector space $\mathrm{H}^0(M, \mathcal{G} \times_G \mathfrak{g}^{(1)})$. In particular, if $G \subset \mathrm{GL}(V)$ is such that $\mathfrak{g}^{(1)} = 0$, then any $G$-structure admits at most one torsion-free affine connection. If $K(\mathfrak{g})$ denotes the $\mathfrak{g}$-module of formal curvature tensors of $\mathfrak{g}$ defined in subsection 2.2. then

$$\mathrm{H}^{1,2}(\mathfrak{g}) = \frac{K(\mathfrak{g})}{\partial(\mathfrak{g}^{(1)} \otimes V^*)};$$

i.e., the cohomology group $\mathrm{H}^{1,2}(\mathfrak{g})$ represents the part of $K(\mathfrak{g})$ which is invariant under $\mathfrak{g}^{(1)}$-valued shifts in a formal torsion-free affine connection with holonomy in $\mathfrak{g}$. For example, if $(G, V) = (\mathrm{CO}(n, \mathbb{R}), \mathbb{R}^n)$, then $\mathfrak{g}^{(1)} = V^*$ and $\mathrm{H}^{1,2}(\mathfrak{g})$ is the vector space of formal Weyl tensors.
Note that $K(\mathfrak{g})$ fits into the exact sequence

$$\mathfrak{g}^{(1)} \otimes V^* \longrightarrow K(\mathfrak{g}) \longrightarrow H^{1,2}(\mathfrak{g}) \longrightarrow 0.$$

In particular, if both $\mathfrak{g}^{(1)}$ and $H^{1,2}(\mathfrak{g})$ vanish, then $K(\mathfrak{g}) = 0$ and hence, $\mathfrak{g}$ is *not* Berger.

If $\mathfrak{g}^{(1)} = 0$, then $\mathrm{H}^{1,2}(\mathfrak{g})$ is exactly $K(\mathfrak{g})$, the $\mathfrak{g}$-module which plays a key role in the theory of torsion-free affine connections with holonomy in $\mathfrak{g}$. The case $\mathfrak{g}^{(1)} = 0$ is generic — there are very few irreducibly acting Lie subgroups $\mathfrak{g} \subset \mathfrak{gl}(V)$ which have $\mathfrak{g}^{(1)} \neq 0$. For future reference we list in Table 5 all complex irreducible Lie subgroups $G \subset \mathrm{GL}(V)$ with $\mathfrak{g}^{(1)} \neq 0$ which is due to Cartan [17] and Kobayashi and Nagano [26]. We emphasize that our classification result also yields the classification of those irreducibly acting subalgebras $\mathfrak{g} \subset \mathfrak{gl}(V)$ with $\mathrm{H}^{1,2}(\mathfrak{g}) \neq 0$.

**Table 5:** The list of all irreducible complex matrix Lie groups $G$ with $\mathfrak{g}^{(1)}\neq 0$

| group $G$ | representation $V$ | $\mathfrak{g}^{(1)}$ |
|---|---|---|
| $\mathrm{SL}(n,\mathbb{C})$ | $V=\mathbb{C}^n$, $n\geqslant 2$ | $(V\otimes\odot^2V^*)_0$ |
| $\mathrm{GL}(n,\mathbb{C})$ | $V=\mathbb{C}^n$, $n\geqslant 1$ | $V\otimes\odot^2V^*$ |
| $\mathrm{GL}(n,\mathbb{C})$ | $V\simeq\odot^2\mathbb{C}^n$, $n\geqslant 2$ | $V^*$ |
| $\mathrm{GL}(n,\mathbb{C})$ | $V\simeq\Lambda^2\mathbb{C}^n$, $n\geqslant 5$ | $V^*$ |
| $\mathrm{GL}(m,\mathbb{C})\cdot\mathrm{GL}(n,\mathbb{C})$ | $V\simeq\mathbb{C}^m\otimes\mathbb{C}^n$, $m,n\geqslant 2$ | $V^*$ |
| $\mathrm{Sp}(n,\mathbb{C})$ | $V=\mathbb{C}^{2n}$, $n\geqslant 2$ | $\odot^3V^*$ |
| $\mathbb{C}^*\cdot\mathrm{Sp}(n,\mathbb{C})$ | $V=\mathbb{C}^{2n}$, $n\geqslant 2$ | $\odot^3V^*$ |
| $\mathrm{CO}(n,\mathbb{C})$ | $V=\mathbb{C}^n$, $n\geqslant 5$ | $V^*$ |
| $\mathbb{C}^*\cdot\mathrm{Spin}(10,\mathbb{C})$ | $V=\mathbb{C}^{16}$ | $V^*$ |
| $\mathbb{C}^*\cdot\mathrm{E}_6^{\mathbb{C}}$ | $V=\mathbb{C}^{27}$ | $V^*$ |

Unless otherwise explicitly stated, $\mathbb{Z}_{\mathbb{C}}$ denotes in what follows either a trivial group or the multiplicative group $\mathbb{C}^*$ and $z_{\mathbb{C}}$ denotes the Lie algebra of $\mathbb{Z}_{\mathbb{C}}$.

3.2. *Twistor formulae for Spencer cohomology groups.* Let $V$ be a finite dimensional complex vector space and $G\subset\mathrm{GL}(V)$ an irreducible representation of a reductive complex Lie group in $V$. Then $G$ also acts irreducibly in $V^*$ via the dual representation. Let $\tilde{X}$ be the $G$-orbit of a highest weight vector in $V^*\setminus 0$. Then the quotient $X:=\tilde{X}/\mathbb{C}^*$ is a compact complex homogeneous-rational manifold called the *sky* of $\mathfrak{g}$. There is a commutative diagram

$$\begin{array}{ccc}\tilde{X} & \hookrightarrow & V^*\setminus 0\\ \downarrow & & \downarrow\\ X & \hookrightarrow & \mathbb{P}(V^*)\end{array}.$$

In fact, $X=G_s/P$, where $G_s$ is the semisimple part of $G$ and $P$ is the parabolic subgroup of $G_s$ leaving a highest weight vector in $V^*$ invariant up to a scale. Let $L$ be the restriction of the hyperplane section bundle $\mathcal{O}(1)$ on $\mathbb{P}(V^*)$ to the submanifold $X$. Clearly, $L$ is an ample homogeneous line bundle on $X$.

In summary, there is a natural map

$$(G,V)\longrightarrow(X,L)$$

which associates with an irreducibly acting reductive Lie group $G\subset\mathrm{GL}(V)$ a pair $(X,L)$ consisting of a compact complex homogeneous-rational manifold $X$ and an ample line bundle $L$ on $X$. We call $(X, L)$ the *Borel-Weil data* associated with $(G, V)$.

Can this map be reversed? According to Borel-Weil, the representation space $V$ can be reconstructed very easily:

$$V = \mathrm{H}^0(X, L).$$

What about $G$? According to Onishchik, with a few (but notable) exceptions, $G$ can be reconstructed as well.

*Fact* 3.1 [2]. Assume that $G$ is simple. The Lie algebra of $G$ is isomorphic to $H^0(X, TX)$ unless one of the following holds:

(i) $G$ is the representation of $\mathrm{Sp}(n, \mathbb{C})$ in $\mathbb{C}^{2n}$ in which case $H^0(X, TX) \simeq \mathfrak{sl}(2n, \mathbb{C})$;

(ii) $G$ is the representation of $\mathrm{G}_2$ in $\mathbb{C}^7$ in which case $H^0(X, TX) \simeq \mathfrak{so}(7, \mathbb{C})$;

(iii) $G$ is the fundamental spinor representation of $\mathrm{Spin}(2n+1, \mathbb{C})$ in which case $H^0(X, TX) \simeq \mathfrak{spin}(2n+2, \mathbb{C})$.

Another proof of this fact is given in [39]. Therefore, if $G \subset \mathrm{GL}(V)$ is semisimple then, with few exceptions, $G$ can be reconstructed from $(X, L)$. However, it is often undesirable to restrict oneself to semisimple groups only (especially in the context of the holonomy classification problem). There is a natural central extension of the Lie algebra $\mathrm{H}^0(X, TX)$:

*Fact* 3.2. For any $(X, L)$, $\mathfrak{g} := H^0(X, L \otimes (J^1L)^*)$ is a reductive Lie algebra canonically represented in $H^0(X, L)$. This fact is easy to explain —

$\mathrm{H}^0(X, L \otimes (J^1L)^*)$ is exactly the Lie algebra of the Lie group $G$ of all global biholomorphisms of the line bundle $L$ which commute with the projection $L \to X$.

In summary, with a given irreducible representation $G \subset \mathrm{GL}(V)$ there is canonically associated a pair $(X, L)$ consisting of a compact complex homogeneous-rational manifold $X$ and a very ample line bundle on $X$ such that much of the original information about $G$ can be restored from $(X, L)$. For our purposes the crucial observation is that the $\mathfrak{g}$-modules $\mathfrak{g}^{(k)}$ and $\mathrm{H}^{k,2}(\mathfrak{g})$ also admit a simple description in terms of $(X, L)$.

THEOREM 3.3. *For a compact complex homogeneous manifold $X$ and a very ample line bundle $L$ on $X$, there is an isomorphism*

$$\mathfrak{g}^{(k)} = \mathrm{H}^0\left(X, L \otimes \odot^{k+1} N^*\right), \quad k = 0, 1, 2, \ldots,$$

*and an exact sequence of* $\mathfrak{g}$*-modules,*

$$\begin{aligned}0 \longrightarrow \mathrm{H}^{k,2}(\mathfrak{g}) \longrightarrow & \mathrm{H}^1\left(X, L\otimes \odot^{k+2}N^*\right)\\ &\longrightarrow \mathrm{H}^1\left(X, L\otimes\odot^{k+1}N^*\right)\otimes V^*, \quad k=1,2,\dots,\end{aligned}$$

*where* $\mathfrak{g} := \mathrm{H}^0(X, L\otimes N^*)$, $N := J^1L$, *and* $\mathrm{H}^{k,2}(\mathfrak{g})$ *are the Spencer cohomology groups associated with the canonical representation of* $\mathfrak{g}$ *in the vector space* $V := \mathrm{H}^0(X, L)$.

*Proof.* Since $L$ is very ample, there is a natural "evaluation" epimorphism

$$V\otimes \mathcal{O}_X \to J^1L \to 0$$

whose dualization gives rise to the canonical monomorphism $0\to N^*\to V^*\otimes \mathcal{O}_X$. Then one may construct the following sequences of locally free sheaves,

$$0 \longrightarrow L\otimes\odot^{k+1}N^* \longrightarrow L\otimes \odot^k N^*\otimes V^* \longrightarrow L\otimes\odot^{k-1}N^*\otimes \Lambda^2V^* \tag{13}$$

and

$$\begin{aligned}0 \longrightarrow L\otimes\odot^{k+2}N^* \longrightarrow L\otimes \odot^{k+1} N^*\otimes V^* \longrightarrow L\otimes\odot^{k}(N^*)\otimes \Lambda^2V^*\\ \longrightarrow L\otimes\odot^{k-1}N^*\otimes\Lambda^3V^*,\end{aligned} \tag{14}$$

and verify that both are exact. (Hint: For any vector space $W$ one has $W\otimes \Lambda^2 W \bmod \Lambda^3 W \simeq W\otimes\odot^2 W \bmod \odot^3 W$.)

Then computing $\mathrm{H}^0(X,\dots)$ of (13) and using the inductive definition of $\mathfrak{g}^{(k)}$ one easily obtains the first statement of the theorem.

The second statement follows from (14) and the definition (12) of $\mathrm{H}^{k,2}(\mathfrak{g})$. Indeed, define $E_k$ by the exact sequence

$$0\longrightarrow L\otimes\odot^{k+2}N^*\longrightarrow L\otimes\odot^{k+1}N^*\otimes V^*\longrightarrow E_k\longrightarrow 0.$$

The associated long exact sequence implies the following *exact* sequence of vector spaces,

$$\begin{aligned}0 \longrightarrow \mathrm{H}^0(X,E_k)/\partial[\mathfrak{g}^{(k)}\otimes V^*] \longrightarrow \mathrm{H}^1\left(X, L\otimes\odot^{k+2}N^*\right)\\ \longrightarrow \mathrm{H}^1\left(X, L\otimes\odot^{k+1}N^*\right)\otimes V^*.\end{aligned}$$

On the other hand, the exact sequence

$$0\longrightarrow E_k\longrightarrow L\otimes\odot^k N^*\otimes\Lambda^2V^*\longrightarrow L\otimes\odot^{k-1}N^*\otimes\Lambda^3V^*$$

implies

$$\mathrm{H}^0(X,E_k) = \ker : \mathfrak{g}^{(k-1)}\otimes\Lambda^2V^* \stackrel{\partial}{\longrightarrow} \mathfrak{g}^{(k-2)}\otimes\Lambda^3V^*,$$

which in turn implies

$$\mathrm{H}^{k,2}(\mathfrak{g}) = \mathrm{H}^0(X,E_k)/\partial[\mathfrak{g}^{(k)}\otimes V^*].$$

This completes the proof of the second part of the theorem. □

COROLLARY 3.4. *If a complex irreducible matrix group $G \subset \mathrm{GL}(V)$ is a Berger group, then the associated Borel-Weil data $(X, L)$ satisfy at least one of the following inequalities*:

$$\dim \mathrm{H}^0\left(X, L \otimes \odot^2 N^*\right) > 0, \tag{15}$$

$$\dim \tilde{\mathrm{H}}^1\left(X, L \otimes \odot^3 N^*\right) > 1, \tag{16}$$

*where $N = J^1 L$, and where*

$$\tilde{\mathrm{H}}^1\left(X, L \otimes \odot^3 N^*\right) = \mathrm{Ker} : \mathrm{H}^1(X, L \otimes \odot^3 N^* \longrightarrow \mathrm{H}^1\left(X, L \otimes \odot^2 N^*\right) \otimes (\mathrm{H}^0(X, L))^*.$$

The short exact sequence $N^* = L^* + TX \otimes L^*$ implies

$$L \otimes \odot^3 N^* = \odot^2 N^* + \odot^3 TX \otimes L^{*2} \tag{17}$$

$$= L^{*2} + TX \otimes L^{*2} + \odot^2 TX \otimes L^{*2} + \odot^3 TX \otimes L^{*2}. \tag{18}$$

LEMMA 3.5. *If $\dim X > 1$ and $X$ is not a quadric*, *then*

$$\mathrm{H}^1(X, L \otimes \odot^3 N^*) \subset \mathrm{H}^1(X, \odot^2 TX \otimes L^{*2}) \oplus \mathrm{H}^1(X, \odot^3 TX \otimes L^{*2}), \tag{19}$$

$$\mathrm{H}^1(X, L \otimes \odot^2 N^*) \subset \mathrm{H}^1(X, TX \otimes L^*) \oplus \mathrm{H}^1(X, \odot^2 TX \otimes L^*). \tag{20}$$

*Proof.* If $\dim X > 1$, then $\mathrm{H}^1(X, L^{*2}) = 0$ by the Kodaira vanishing theorem. If $X$ is not a quadric, then, by Theorem B, for any ample line bundle $L$ on $X$ one has $\mathrm{H}^1(X, TX \otimes L^{*2}) = 0$. Then long exact sequences of the short exact sequences (18) immediately imply (19). The second estimation can be checked analogously. □

In conclusion, every irreducibly acting holonomy group $\mathfrak{g} \subset \mathrm{End}(V)$ has an associated compact homogeneous-rational manifold $X$, the sky of $\mathfrak{g}$. The class of all such manifolds contains a very special subclass — the family of compact hermitean-symmetric manifolds. Therefore, there is a natural question: Which holonomy groups have the members of this subfamily as their skies? The answer to this question is given in the next section and provides us with a surprising insight into the structure of the completed Berger list.

## 4. Holonomy groups associated with

**compact Hermitean-symmetric manifolds**

The list of irreducible compact hermitean symmetric manifolds is very short:

(i) Projective spaces $\mathbb{P}_n =$ 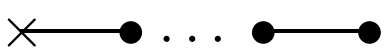 ($n$ nodes, $n \geqslant 1$). Note also that there is another representation for odd-dimensional projective spaces; cf. Fact 3.1.

$\mathbb{P}_{2n+1} =$ 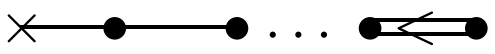 .

(ii) Grassmanians $G(k, \mathbb{C}^{n+1}) =$ 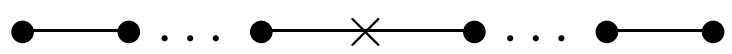 ($n$ nodes, $k^{\text{th}}$ node is crossed, $2 \leqslant k \leqslant (n+1)/2, n \geqslant 3$).

(iii) Quadrics $Q_{2n-1} =$ 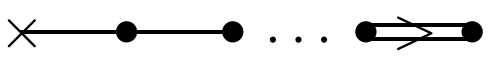 ($n$ nodes, $n \geqslant 2$) and

$Q_{2n-2} =$ 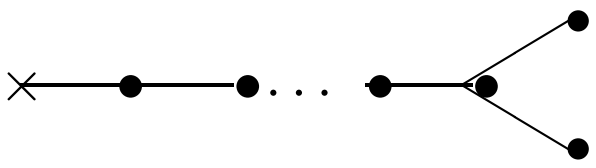 ($n$ nodes, $n \geqslant 4$).

Note that $Q_4 = G(2, \mathbb{C}^4)$, and 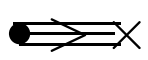 is biholomorphic to $Q_5$; cf. Fact 3.1.

(iv) Manifolds $Y_n =$ 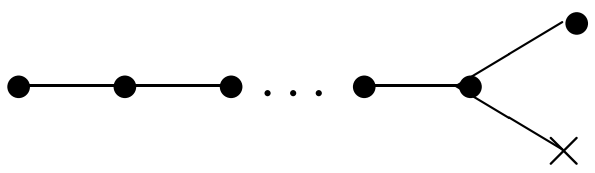 ($n$ nodes in both cases, $n \geqslant 5$).

Note that they are biholomorphic to 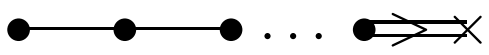 ($n-1$ nodes); cf. Fact 3.1.

(v) 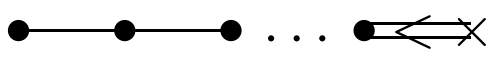 ($n$ nodes, $n \geqslant 3$),

(vi) 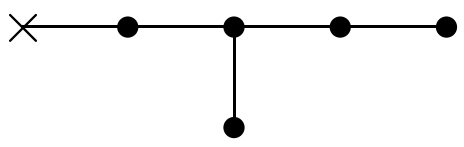 ,

(vii) 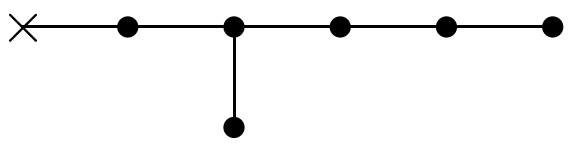 .

Theorem 4.1. *Let $\mathfrak{g} \subset \mathfrak{gl}(V)$ be an irreducible Berger algebra. Then its sky is biholomorphic to an irreducible compact complex Hermitean manifold if and only if $\mathfrak{g}$ is an entry of Table* 6.

**Table 6**

| No. | group $G$ | representation $V$ | restrictions |
|---|---|---|---|
| 1 | $\mathbb{Z}_{\mathbb{C}} \cdot \mathrm{SL}(n, \mathbb{C})$ | $\mathbb{C}^n$ | $n \geqslant 1$ |
| 2 | $\mathrm{GL}(n, \mathbb{C})$ | $\odot^2\mathbb{C}^n$ | $n \geqslant 3$ |
| 3 | $\mathrm{GL}(n, \mathbb{C})$ | $\Lambda^2\mathbb{C}^n$ | $n \geqslant 3$ |
| 4 | $\mathbb{Z}_{\mathbb{C}} \cdot \mathrm{SL}(2, \mathbb{C})$ | $\odot^3\mathbb{C}^2$ | |
| 5 | $\mathrm{SL}(6, \mathbb{C})$ | $\Lambda^3\mathbb{C}^6$ | |
| 6 | $\mathbb{Z}_{\mathbb{C}} \cdot \mathrm{SO}(n, \mathbb{C})$ | $\mathbb{C}^n$ | $n \geqslant 3$ |
| 7 | $\mathbb{Z}_{\mathbb{C}} \cdot \mathrm{Sp}(n, \mathbb{C})$ | $\mathbb{C}^{2n}$ | $n \geqslant 2$ |
| 8 | $\mathrm{Sp}(3, \mathbb{C})$ | $\mathbb{C}^{14} \subset \Lambda^3\mathbb{C}^6$ | |
| 9 | $\mathrm{G}_2^{\mathbb{C}}$ | $\mathbb{C}^7$ | |
| 10 | $\mathbb{Z}_{\mathbb{C}} \cdot \mathrm{E}_6^{\mathbb{C}}$ | $\mathbb{C}^{27}$ | |
| 11 | $\mathrm{E}_7^{\mathbb{C}}$ | $\mathbb{C}^{56}$ | |
| 12 | $\mathrm{Spin}(7, \mathbb{C})$ | $\mathbb{C}^8$ | |
| 13 | $\mathbb{Z}_{\mathbb{C}} \cdot \mathrm{Spin}(10, \mathbb{C})$ | $\mathbb{C}^{16}$ | |
| 14 | $\mathrm{Spin}(12, \mathbb{C})$ | $\mathbb{R}^{64} \simeq \mathbb{C}^{32}$ | |
| Notation: $\mathbb{Z}_{\mathbb{C}}$ denotes either the trivial group or $\mathbb{C}^*$. | | | |

*Proof.* All the entries in Table 6 are known to be Berger algebras (for the new exotic entries, numbers 5, 8, 11 and 12, this will be shown in detail in §7) and to have hermitean symmetric skies. So Theorem 4.1 will be proved if we show that there are no other irreducible matrix Lie algebras with these properties.

*Case* (i). Let $L$ be an ample line bundle on $X = \mathbb{P}_n$, i.e. $L = \mathcal{O}(s)$ for some $s \in \mathbb{N}$. Using the isomorphism $J^1L = \mathbb{C}^{n+1} \otimes \mathcal{O}(s-1)$, we show easily

$$\mathrm{H}^0(X, L \otimes \odot^2 N^*) \quad = \quad \begin{cases} \mathbb{C}^{n+1} \otimes \odot^2\mathbb{C}^{n+1} & \text{for } s = 1, \\ \odot^2\mathbb{C}^{n+1} & \text{for } s = 2, \\ 0 & \text{otherwise,} \end{cases}$$

$$\tilde{\mathrm{H}}^1(X, L\otimes \odot^3 N^*) \;=\; \begin{cases} \mathbb{C}^2\otimes\odot^3\mathbb{C}^2 & \text{for } n=1, s=3,\\ \mathbb{C} & \text{for } n=1, s=4,\\ 0 & \text{otherwise,}\end{cases}$$

which together with Fact 3.1 imply that the only irreducible representations which have $\mathbb{P}_n$ as the sky and satisfy the criteria (15) and (16) are the entries No.1, 2, 4 and 7 in Table 6.

*Case* (ii). Any ample line bundle $L$ on $X = G(k, \mathbb{C}^{n+1})$ is of the form

$$L = \overset{0}{\bullet}\!-\!\overset{0}{\bullet}\ \ldots\ \overset{0}{\bullet}\!-\!\overset{s}{\times}\!-\!\overset{0}{\bullet}\ \ldots\ \overset{0}{\bullet}\!-\!\overset{0}{\bullet}$$

for some $s\in\mathbb{N}$. Since

$$TX = \overset{1}{\bullet}\!-\!\overset{0}{\bullet}\ \ldots\ \overset{0}{\bullet}\!-\!\overset{0}{\times}\!-\!\overset{0}{\bullet}\ \ldots\ \overset{0}{\bullet}\!-\!\overset{1}{\bullet}$$

one has

$$\odot^2TX = \overset{2}{\bullet}\!-\!\overset{0}{\bullet}\ \ldots\ \overset{0}{\bullet}\!-\!\overset{0}{\times}\!-\!\overset{0}{\bullet}\ \ldots\ \overset{0}{\bullet}\!-\!\overset{2}{\bullet} \quad\oplus\quad \overset{0}{\bullet}\!-\!\overset{1}{\bullet}\ \ldots\ \overset{0}{\bullet}\!-\!\overset{0}{\times}\!-\!\overset{0}{\bullet}\ \ldots\ \overset{1}{\bullet}\!-\!\overset{0}{\bullet}$$

and

$$\begin{aligned}\odot^3TX \;=\; & \overset{3}{\bullet}\!-\!\overset{0}{\bullet}\ \ldots\ \overset{0}{\bullet}\!-\!\overset{0}{\times}\!-\!\overset{0}{\bullet}\ \ldots\ \overset{0}{\bullet}\!-\!\overset{3}{\bullet}\\ & \oplus\ \overset{1}{\bullet}\!-\!\overset{1}{\bullet}\ \ldots\ \overset{0}{\bullet}\!-\!\overset{0}{\times}\!-\!\overset{0}{\bullet}\ \ldots\ \overset{1}{\bullet}\!-\!\overset{1}{\bullet}\ \oplus\ \overset{0}{\bullet}\!-\!\overset{0}{\bullet}\!-\!\overset{1}{\bullet}\ \ldots\ \overset{0}{\bullet}\!-\!\overset{0}{\times}\!-\!\overset{0}{\bullet}\ \ldots\ \overset{1}{\bullet}\!-\!\overset{0}{\bullet}\!-\!\overset{0}{\bullet}\end{aligned}$$

where the $k^{\mathrm{th}}$ node is crossed, $k\geqslant 2$, and when $k=2$ the third summand in the latter formula is assumed to be zero.

Therefore,

$$\begin{aligned}\mathrm{H}^0(X, L\otimes\odot^2N^*) &= \mathrm{H}^0(X, \odot^2TX\otimes L^*)\\ &= \begin{cases}\Lambda^2\mathbb{C}^{n+1} & \text{for } \ X = G(2,\mathbb{C}^{n+1}),\ s=1,\\ 0 & \text{otherwise,}\end{cases}\end{aligned}$$

in accordance with the Cartan-Kobayashi-Nagano classification (see Table 5).

By the Bott-Borel-Weil theorem,

$$\begin{aligned}&\mathrm{H}^1\big(\overset{a}{\bullet}\ \ldots\ \overset{b}{\bullet}\!-\!\!-\!\overset{-s}{\bullet}\!-\!\!-\!\overset{c}{\bullet}\ \ldots\ \overset{d}{\bullet}\big)\\ &= \begin{cases}\overset{a}{\bullet}\ \ldots\ \overset{1+b-s}{\bullet}\!-\!\!-\!\overset{s-2}{\bullet}\!-\!\!-\!\overset{1+c-s}{\bullet}\ \ldots\ \overset{d}{\bullet} & \text{if } s\geqslant 2, c\geqslant s-2, d\geqslant s-1,\\ 0 & \text{otherwise.}\end{cases}\end{aligned}$$

Then one easily finds that

$$\mathrm{H}^1(X, \odot^2TX\otimes L^{*2}) = \begin{cases}\mathbb{C} & \text{for } \ (X,L) = \overset{0}{\bullet}\!-\!\overset{0}{\bullet}\!-\!\overset{1}{\times}\!-\!\overset{0}{\bullet}\!-\!\overset{0}{\bullet}\\ 0 & \text{otherwise,}\end{cases}$$

and

$$\mathrm{H}^1(X, \odot^3 TX \otimes L^{*2}) = \begin{cases} \overset{1}{\bullet}-\overset{0}{\bullet}-\overset{0}{\bullet}-\overset{0}{\bullet}-\overset{1}{\bullet} & \text{for } (X, L) = \overset{0}{\bullet}-\overset{0}{\bullet}-\overset{1}{\times}-\overset{0}{\bullet}-\overset{0}{\bullet} \\ \mathbb{C} & \text{for } (X, L) = \overset{0}{\bullet}-\overset{0}{\bullet}-\overset{0}{\bullet}-\overset{1}{\times}-\overset{0}{\bullet}-\overset{0}{\bullet}-\overset{0}{\bullet} \\ 0 & \text{otherwise.} \end{cases}$$

In addition to entry No.3 in Table 6 only the following two representations

$$\mathfrak{g} = \begin{cases} \overset{0}{\bullet}-\overset{0}{\bullet}-\overset{1}{\bullet}-\overset{0}{\bullet}-\overset{0}{\bullet} \\ \overset{0}{\bullet}-\overset{0}{\bullet}-\overset{0}{\bullet}-\overset{1}{\bullet}-\overset{0}{\bullet}-\overset{0}{\bullet}-\overset{0}{\bullet} \end{cases}$$

may have, by the estimation (19), $K(\mathfrak{g} \oplus \mathbb{C})$ nonvanishing. Moreover, the long exact sequence of (17) allows us to compute $K(\mathfrak{g})$ precisely. Indeed, since $\mathrm{H}^0(X, \odot^3 TX \otimes L^{*2}) = \mathrm{H}^1(X, \odot^2 TX \otimes L^{*2}) = \mathbb{C}$ for the first representation and $\mathrm{H}^0(X, \odot^3 TX \otimes L^{*2}) = \mathrm{H}^1(X, \odot^2 TX \otimes L^{*2}) = 0$ for the second one, the following piece of the long exact sequence

$$\begin{aligned} 0 \quad \to \quad & \mathrm{H}^0(X, \odot^3 TX \otimes L^{*2}) \to \mathrm{H}^1(X, \odot^2 TX \otimes L^{*2}) \\ \to \quad & \mathrm{H}^1(X, L \otimes \odot^3 N^*) \to \mathrm{H}^1(X, \odot^3 TX \otimes L^{*2}) \to 0 \end{aligned}$$

implies

$$K(\mathfrak{g} \oplus \mathbb{C}) = K(\mathfrak{g}) = \mathrm{H}^1(X, \odot^3 TX \otimes L^{*2}).$$

In fact the statement $K(\mathfrak{g} \oplus \mathbb{C}) = K(\mathfrak{g})$ above requires an additional simple calculation involving the first Bianchi identity, which we postpone to Section 7.

Therefore, the only irreducible representations which have a Grassmanian $G(k, \mathbb{C}^{n+1})$ as the sky and satisfy the criteria (15) and (16) are entries No. 3 and 5 in Table 6.

*Case* (iii). Let $i : Q_n \hookrightarrow V$ be the standard embedding, where $V$ is an $(n+1)$-dimensional vector space. Denote by [1] the 1-dimensional subspace in $\odot^2 V^*$ spanned by the conformal metric $[g_{ab}]$, and by $[-1] = [1]^*$.

Any ample line bundle $L$ on $Q_n$ is of the form $i^*(\mathcal{O}(s))$ for some $s \in \mathbb{N}$. The long exact sequence of the extension

$$0 \longrightarrow N^* \longrightarrow V^* \otimes (\mathcal{O}(1-s)) \longrightarrow i^*(\mathcal{O}(2-s))[-1] \longrightarrow 0$$

immediately implies

$$\mathrm{H}^1(X, L \otimes \odot^3 N^*) = \begin{cases} \text{space of Weyl tensors} & \text{for } s = 1, \\ \odot^2 V^*[-1] = [-1] \oplus \odot^2_0 V^* & \text{for } s = 2, \\ 0 & \text{otherwise,} \end{cases}$$

$$\mathrm{H}^1(X, L \otimes \odot^2 N^*) = \begin{cases} \Lambda^2 V^* & \text{for } s = 2, \\ V^*[-1] & \text{for } s = 3, \\ 0 & \text{otherwise.} \end{cases}$$

Then in the case $s = 2$ we have

$$\begin{array}{rcccc}\tilde{\mathrm{H}}^1(X, L \otimes \odot^3 N^*) & = & \ker : \odot^2 V^*[-1] & \longrightarrow & \Lambda^2 V^* \otimes \odot^2_0 V^* \\ & & \phi_{ab} & \longrightarrow & \phi_{ac} g_{bd} + \phi_{bc} g_{ad} - \phi_{ad} g_{bc} - \phi_{bd} g_{ac}\end{array}$$

implying $\tilde{\mathrm{H}}^1(X, L \otimes \odot^3 N^*) = [-1] \simeq \mathbb{C}$.

Taking into account Fact 3.1 we conclude that it is only the classical representations No. 6, exceptional ones No. 9 and No. 11 and the central extensions of the latter two which have the quadric as their sky. Their central extensions have already been ruled out as Berger algebras in Corollary 2.6.

*Case* (iv). Any ample line bundle $L \to X$ is of the form

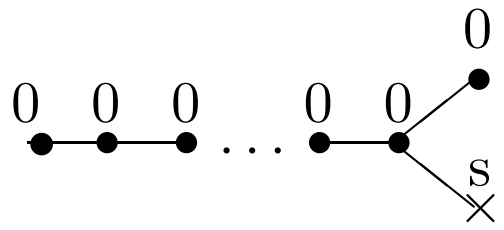

for some positive integer $s$. Since

$TX =$ 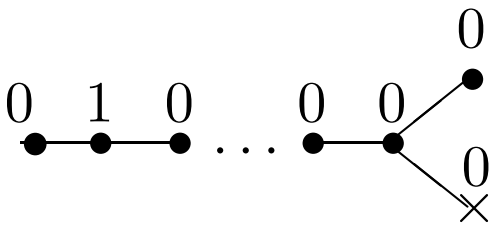

one has

$\odot^2 TX =$ 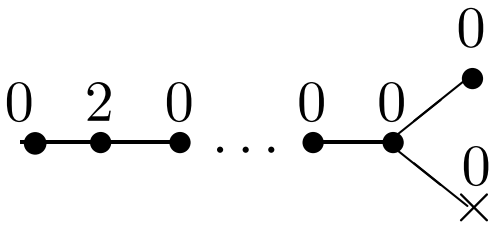

$\oplus$ 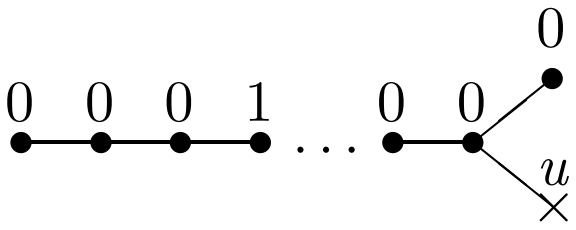

and

$\odot^3 TX =$ 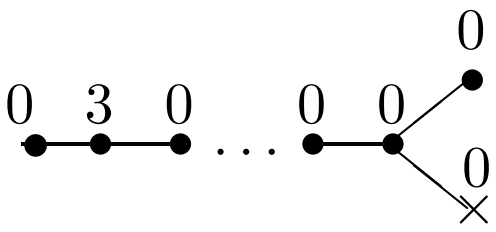

$\oplus$ 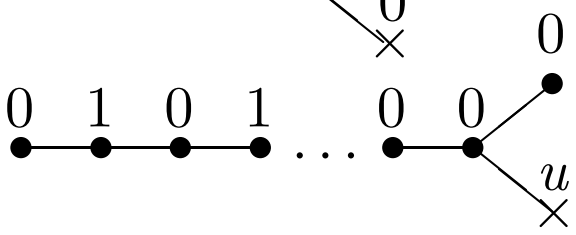

$\oplus$ 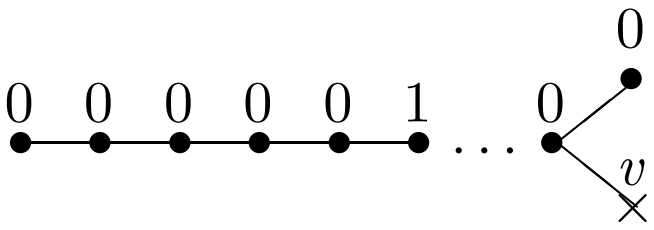

,

where

$$u = \begin{cases} 1 & \text{for the number of nodes n=5,} \\ 0 & \text{otherwise,} \end{cases} \quad , \quad v = \begin{cases} 1 & \text{for n=7,} \\ 0 & \text{otherwise.} \end{cases}$$

Therefore,

$$\mathrm{H}^0(X, L \otimes \odot^2 N^*) = \mathrm{H}^0(X, \odot^2 TX \otimes L^*)$$

$=$ 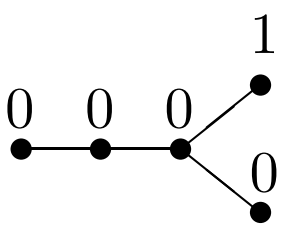

for $n = 5,\ s = 1$,

$0$ otherwise,

in accordance with the Cartan-Kobayashi-Nagano classification (see Table 5).

By the Bott-Borel-Weil theorem

$\mathrm{H}^1($ 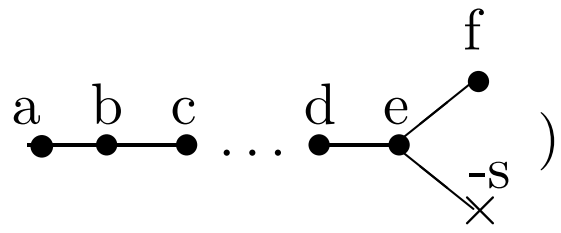

$)$

$=$ 
- 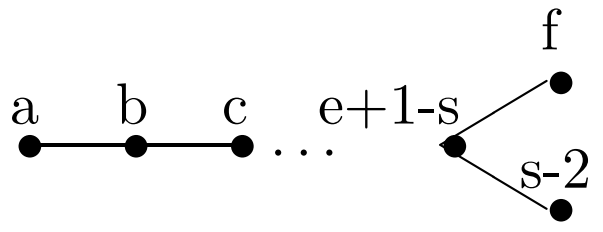

for $s \geqslant 2, e \geqslant s-1$,
- $0$ otherwise,

which immediately implies

$\mathrm{H}^1(X, \odot^2 TX \otimes L^{*2}) =$
- $\mathbb{C}$ for $(X, L) =$ 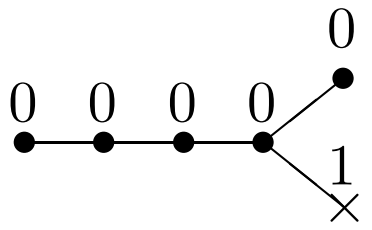

- $0$ otherwise,

and

$\mathrm{H}^1(X, \odot^3 TX \otimes L^{*2})$

$=$
- 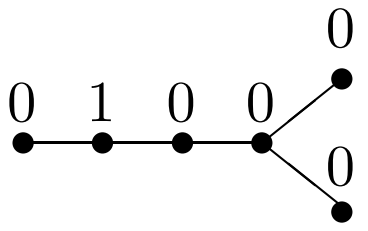

for $(X, L) =$ 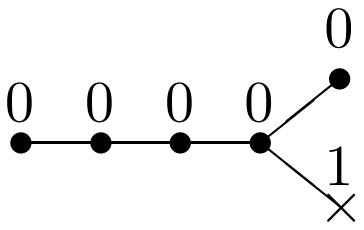

- $\mathbb{C}$ for $(X, L) =$ 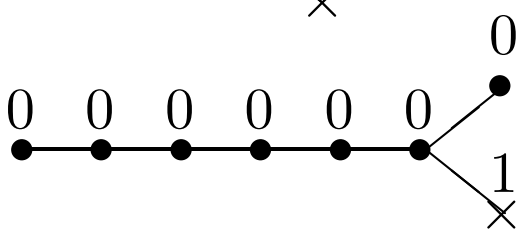

- $0$ otherwise.

In addition to entry No. 13 in Table 6, only the following two representations of $\mathfrak{so}(2n, \mathbb{C})$

$\mathfrak{g} =$ 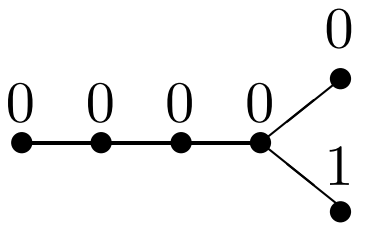

and $\mathfrak{g} =$ 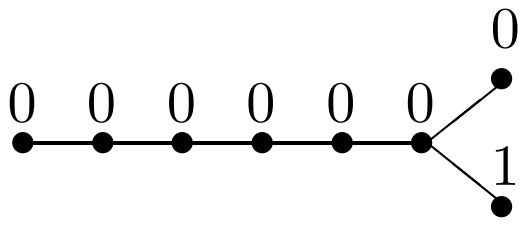

may have, by estimation (19), $K(\mathfrak{g} \oplus \mathbb{C})$ nonvanishing. Analyzing the long exact sequence of (17) in exactly the same way as in the case (ii) above, one obtains

$$K(\mathfrak{g} \oplus \mathbb{C}) = K(\mathfrak{g}) = \mathrm{H}^1(X, \odot^3 TX \otimes L^{*2}).$$

Therefore, we are left with only three candidates to holonomies having the sky biholomorphic to type (iv): the entries No. 13, 14 in Table 6 which are known to satisfy the Berger criteria (that $K^1(\mathfrak{g}) \neq 0$ for No. 14 will be shown in detail in §8) and the representations

$\mathfrak{g}_1 =$ 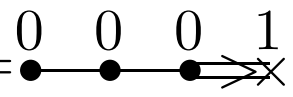

and $\mathfrak{g}_2 =$ 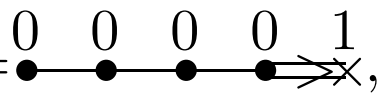

,

which, by Fact 3.1, may, in principle, have $K(\mathfrak{g}_i \oplus \mathbb{C}) \neq 0$. But the Borel-Weil approach can not be used to compute these modules. However, $\mathfrak{g}_1$ corresponds to the spinor representation of $\mathfrak{spin}(9)$; hence $\mathfrak{g}_1 \oplus \mathbb{C}\mathrm{Id}$ is not a Berger algebra by Corollary 2.6, and $\mathfrak{g}_1$ is known not to be a Berger algebra, either [4].

The second case, which corresponds to the spinor representation of $\mathfrak{spin}(11)$, will be ruled out as a Berger algebra in Lemma 7.3.

*Case* (v). Any ample line bundle $L \to X$ is of the form

$$\overset{0}{\bullet}\!-\!\overset{0}{\bullet}\!-\!\overset{0}{\bullet}\ldots\overset{0}{\bullet}\!\Leftarrow\!\overset{s}{\times}$$

for some positive integer $s$.

Since

$$TX = \overset{2}{\bullet}\!-\!\overset{0}{\bullet}\!-\!\overset{0}{\bullet}\ldots\overset{0}{\bullet}\!\Leftarrow\!\overset{0}{\times},$$

one has

$$\odot^2 TX = \overset{4}{\bullet}\!-\!\overset{0}{\bullet}\!-\!\overset{0}{\bullet}\ldots\overset{0}{\bullet}\!\Leftarrow\!\overset{0}{\times} \;\oplus\; \overset{0}{\bullet}\!-\!\overset{2}{\bullet}\!-\!\overset{0}{\bullet}\ldots\overset{0}{\bullet}\!\Leftarrow\!\overset{0}{\times}$$

and

$$\odot^3 TX = \overset{6}{\bullet}\!-\!\overset{0}{\bullet}\!-\!\overset{0}{\bullet}\ldots\overset{0}{\bullet}\!\Leftarrow\!\overset{0}{\times} \;\oplus\; \overset{2}{\bullet}\!-\!\overset{2}{\bullet}\!-\!\overset{0}{\bullet}\ldots\overset{0}{\bullet}\!\Leftarrow\!\overset{0}{\times} \;\oplus\; \overset{0}{\bullet}\!-\!\overset{0}{\bullet}\!-\!\overset{2}{\bullet}\ldots\overset{0}{\bullet}\!\Leftarrow\!\overset{0}{\times}.$$

By the Bott-Borel-Weil theorem

$$\mathrm{H}^1(\overset{a}{\bullet}\!-\!\overset{b}{\bullet}\!-\!\overset{c}{\bullet}\ldots\overset{d}{\bullet}\!\Leftarrow\!\overset{-s}{\times}) = \begin{cases} \overset{a}{\bullet}\!-\!\overset{b}{\bullet}\!-\!\overset{c}{\bullet}\ldots\overset{2+d-2s}{\bullet}\!\Leftarrow\!\overset{s-2}{\bullet} & \text{for } s \geqslant 2, d \geqslant 2s-2, \\ 0 & \text{otherwise,} \end{cases}$$

which implies

$$\mathrm{H}^1(X, \odot^2 TX \otimes L^{*2}) = \begin{cases} \mathbb{C} & \text{for } (X,L) = \overset{0}{\bullet}\!-\!\overset{0}{\bullet}\!\Leftarrow\!\overset{1}{\times} \\ 0 & \text{otherwise,} \end{cases}$$

and

$$\mathrm{H}^1(X, \odot^3 TX \otimes L^{*2}) = \begin{cases} \overset{2}{\bullet}\!-\!\overset{0}{\bullet}\!\Leftarrow\!\overset{0}{\bullet} & \text{for } (X,L) = \overset{0}{\bullet}\!-\!\overset{0}{\bullet}\!\Leftarrow\!\overset{1}{\times} \\ \mathbb{C} & \text{for } (X,L) = \overset{0}{\bullet}\!-\!\overset{0}{\bullet}\!-\!\overset{0}{\bullet}\!-\!\overset{0}{\bullet}\!\Leftarrow\!\overset{1}{\times} \\ 0 & \text{otherwise.} \end{cases}$$

By the estimation (19), only representations

$$\mathfrak{g} = \overset{0}{\bullet}\!-\!\overset{0}{\bullet}\!\Leftarrow\!\overset{1}{\bullet} \quad \text{and} \quad \mathfrak{g} = \overset{0}{\bullet}\!-\!\overset{0}{\bullet}\!-\!\overset{0}{\bullet}\!-\!\overset{0}{\bullet}\!\Leftarrow\!\overset{1}{\bullet}$$

may have $K(\mathfrak{g}) \neq 0$. Moreover, a straightforward analysis of the long exact sequence of (17) gives a precise expression for $K(\mathfrak{g})$:

$$K(\mathfrak{g} \oplus \mathbb{C}) = K(\mathfrak{g}) = \mathrm{H}^1(X, \odot^3 TX \otimes L^{*2}).$$

Therefore, the only irreducible representation which has a submanifold of the type (v) as the sky and satisfies the criteria (15) and (16) is entry 8 in Table 6.

*Case* (vi). Any ample line bundle $L \to X$ is of the form

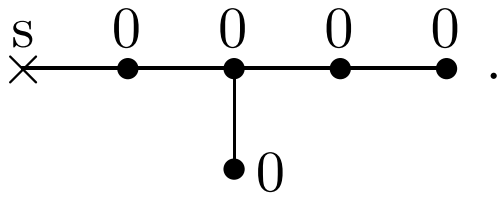

for some positive integer $s$.

Since

$TX =$ 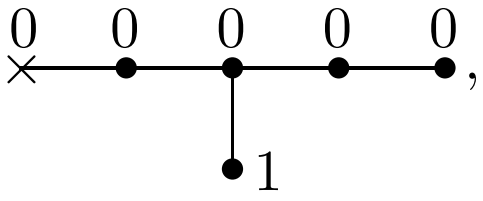

one has

$\odot^2 TX =$ 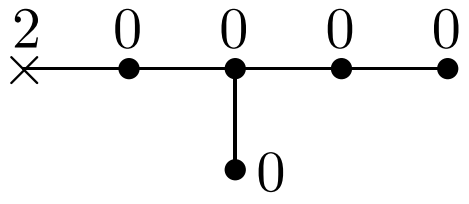

$\oplus$ 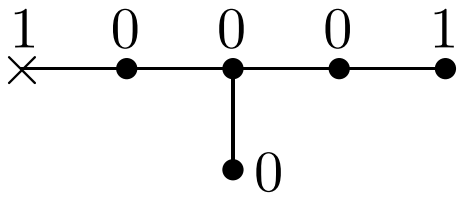

and

$\odot^3 TX =$ 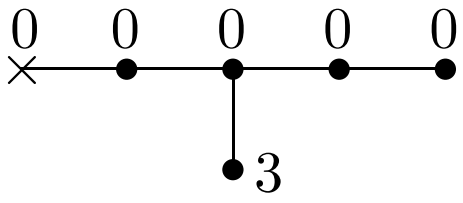

$\oplus$ 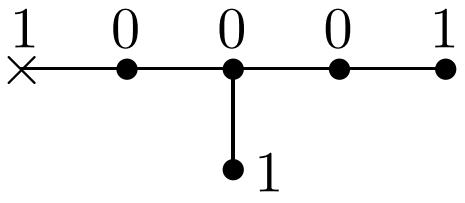

Then, by the Bott-Borel-Weil theorem,

$$\mathrm{H}^1(X, \odot^2 TX \otimes L^{*2}) = \mathrm{H}^1(X, \odot^3 TX \otimes L^{*2}) = 0$$

for any $s$, implying that all irreducible representations $\mathfrak{g}$ of $\mathfrak{e} \oplus \mathbb{C}$ have $H^{1,2}(\mathfrak{g}) = 0$, while, in agreement with Table 5, one has

$\mathrm{H}^0(X, \odot^2 TX \otimes L^*) = \mathrm{H}^0(X, L \otimes \odot^* N^{*2}) =$ 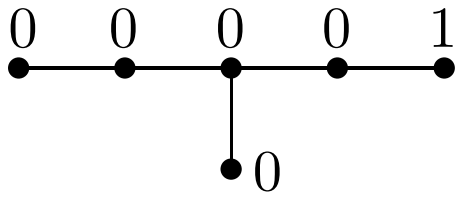

for $s = 1$; $0$ otherwise.

Therefore, the only irreducible representation which has a submanifold of type (vi) as the sky and satisfies the criteria (15) and (16) is entry 10 in Table 6.

*Case* (vii). Any ample line bundle $L \to X$ is of the form

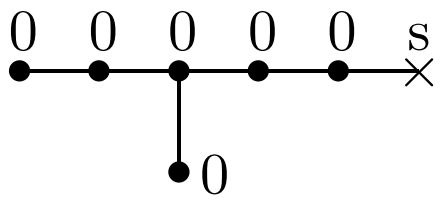

for some positive integer $s$.

Since

$TX =$ 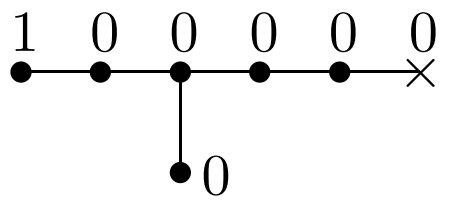

,

one has

$\odot^2 TX =$ 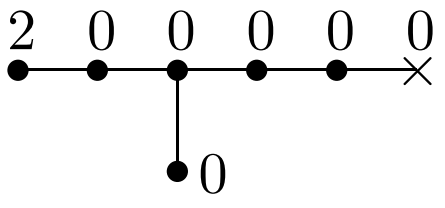

$\oplus$ 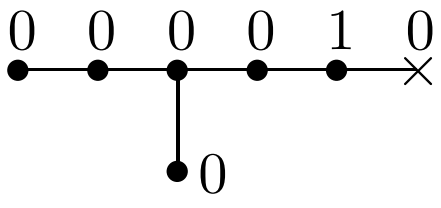

and

$\odot^3 TX =$ 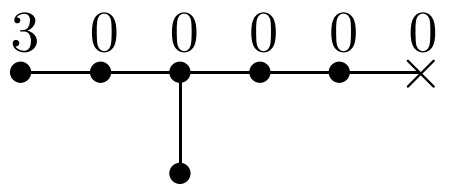

$\oplus$ 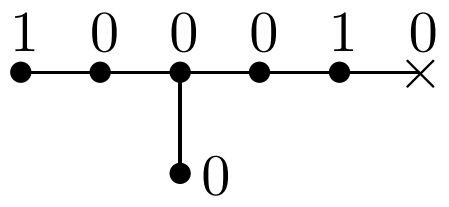

$\oplus$ 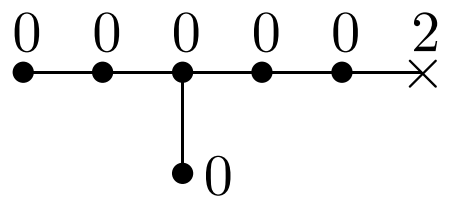

.

Then, by the Bott-Borel-Weil theorem,

$$\mathrm{H}^1(X, \odot^2 TX \otimes L^{*2}) = \begin{cases} \mathbb{C} & \text{for } s = 1 \\ 0 & \text{otherwise} \end{cases}$$

and

$$\mathrm{H}^1(X, \odot^3 TX \otimes L^{*2}) = \begin{cases} \mathrm{Ad}(\mathfrak{e}_7) & \text{for } s = 1 \\ 0 & \text{otherwise.} \end{cases}$$

By estimation (19), only the representation

$\mathfrak{g} =$ 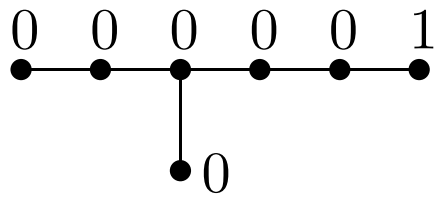

may have $K(\mathfrak{g}) \neq 0$. Again, the long exact sequence of (17) gives a precise expression for $K(\mathfrak{g})$:

$$K(\mathfrak{g} \oplus \mathbb{C}) = K(\mathfrak{g}) = \mathrm{H}^1(X, \odot^3 TX \otimes L^{*2}) = \mathrm{Ad}(\mathfrak{e}_7).$$

The statement $K(\mathfrak{g} \oplus \mathbb{C}) = K(\mathfrak{g})$ requires an additional calculation involving the first Bianchi identities. The necessary details are given in Section 7.

Therefore, the only irreducible representation which has a submanifold of type (v) as the sky and satisfies the criteria (15) and (16) is entry 11 in Table 6.

□

## 5. Classification of Segre holonomies

5.1. *Segre holonomies.* The main result of this section is the following:

THEOREM 5.1. *Let $G$ be the irreducible holonomy of a torsion-free affine connection which is not locally symmetric. If the semisimple part of $G$ is not simple, then $G$ is one of the groups listed in Table 7.*

**Table 7**: List of Segre holonomies

| group $G$ | representation $V$ | restrictions |
|---|---|---|
| $T_{\mathbb{R}} \cdot SL(m,\mathbb{R}) \cdot SL(n,\mathbb{R})$ | $\mathbb{R}^m \otimes \mathbb{R}^n \simeq \mathbb{R}^{mn}$ | $m \geqslant n \geqslant 2$ |
| $T_{\mathbb{C}} \cdot SL(m,\mathbb{C}) \cdot SL(n,\mathbb{C})$ | $\mathbb{C}^m \otimes \mathbb{C}^n \simeq \mathbb{R}^{2mn}$ | $m \geqslant n \geqslant 2$ |
| $T_{\mathbb{R}} \cdot SL(m,\mathbb{H}) \cdot SL(n,\mathbb{H})$ | $\mathbb{H} \otimes \mathbb{H} \simeq \mathbb{R}^{4mn}$ | $m \geqslant n \geqslant 1$ |
| $Sp(1) \cdot Sp(p,q)$ | $\mathbb{H}^{p+q} \simeq \mathbb{R}^{4(p+q)}$ | $p+q \geqslant 2$ |
| $SL(2,\mathbb{R}) \cdot Sp(n,\mathbb{R})$ | $\mathbb{R}^2 \otimes \mathbb{R}^{2n} \simeq \mathbb{R}^{4n}$ | $n \geqslant 2$ |
| $SL(2,\mathbb{C}) \cdot Sp(n,\mathbb{C})$ | $\mathbb{C}^2 \otimes \mathbb{C}^{2n} \simeq \mathbb{R}^{8n}$ | $n \geqslant 2$ |
| $SL(2,\mathbb{R}) \cdot SO(p,q)$ | $\mathbb{R}^2 \otimes \mathbb{R}^{p+q} \simeq \mathbb{R}^{2(p+q)}$ | $p+q \geqslant 3$ |
| $Sp(1) \cdot SO(n,\mathbb{H})$ | $\mathbb{H}^n \simeq \mathbb{R}^{4n}$ | $n \geqslant 2$ |
| $SL(2,\mathbb{C}) \cdot SO(n,\mathbb{C})$ | $\mathbb{C}^2 \otimes \mathbb{C}^n \simeq \mathbb{R}^{4n}$ | $n \geqslant 3$ |
| Notation: $T_{\mathbb{F}}$ denotes any connected Lie subgroup of $\mathbb{F}^*$. | | |

The class of holonomy groups (and the associated geometric structures) studied by the above theorem appear in the literature under different names. For example, the authors of [3], [29] call these *almost Grassmanian* and the authors of [6] call these *paraconformal*. In this paper we follow the terminology of Bryant [16] who suggested calling them *Segre* holonomies. Some applications of Segre structures to high energy physics are discussed in [28], [31].

5.2. *Cohomology on reducible homogeneous-rational manifolds.* From now on we assume that $\mathbb{X} = X_1 \times X_2$ is a direct product of two compact complex homogeneous-rational manifolds $X_1$ and $X_2$ and that $\mathbb{L}$ is an ample holomorphic line bundle on $X$. Denoting by $\pi_1 : \mathbb{X} \to X_1$ and $\pi_2 : \mathbb{X} \to X_2$ the natural projections, we may write $\mathbb{L} = \pi_1^*(L_1) \otimes \pi_2^*(L_2)$ for some uniquely specified ample line bundles $L_1$ and $L_2$ on $X_1$ and $X_2$ respectively. We denote $\mathbb{N} := J^1\mathbb{L}$ and $N_i := J^1 L_i$, $i = 1, 2$.

Since

$$0 \longrightarrow \Omega X_i \otimes L_i \longrightarrow N_i \longrightarrow L_i \longrightarrow 0,$$

one has

$$0 \quad \longrightarrow \quad \pi_1^*(\Omega^1 X_1) \otimes \mathbb{L} + \pi_2^*(\Omega^1 X_2) \otimes \mathbb{L}$$

$$\longrightarrow \quad \pi_1^*(N_1)\otimes\pi_2^*(L_2)+\pi_2^*(N_2)\otimes\pi_1^*(L_1)\longrightarrow \mathbb{L}+\mathbb{L}\longrightarrow 0.$$

The latter extension combined with

$$0\longrightarrow \pi_1^*(\Omega^1X_1)\otimes\mathbb{L}+\pi_2^*(\Omega^1X_2)\otimes\mathbb{L}\longrightarrow\mathbb{N}\longrightarrow\mathbb{L}\longrightarrow 0.$$

implies

$$0\longrightarrow\mathbb{N}\longrightarrow\pi_1^*(N_1)\otimes\pi_2^*(L_2)+\pi_1^*(L_1)\otimes\pi_2^*(N_2)\longrightarrow\mathbb{L}\longrightarrow 0,$$

or

$$0\longrightarrow\mathbb{L}^*\longrightarrow\pi_1^*(N_1^*)\otimes\pi_2^*(L_2^*)+\pi_1^*(N_1^*)\otimes\pi_2^*(L_2^*)\longrightarrow\mathbb{N}^*\longrightarrow 0, \tag{21}$$

which in turn implies the following two exact sequences

(22)

$$0\longrightarrow \begin{array}{c}\pi_1^*(N_1^*)\otimes\pi_2^*(L_2^*)\\ +\\ \pi_1^*(L_1^*)\otimes\pi_2^*(N_2^*)\end{array}\longrightarrow\begin{array}{c}\pi_1^*(L_1\otimes\odot^2N_1^*)\otimes\pi_2^*(L_2^*)\\ +\\ \pi_1^*(N_1^*)\otimes\pi_2^*(N_2^*)\\ +\\ \pi_1^*(L_1^*)\otimes\pi_2^*(L_2\otimes\odot^2N_2^*)\end{array}\longrightarrow\mathbb{L}\otimes\odot^2\mathbb{N}^*\longrightarrow 0$$

and

(23)

$$0\longrightarrow \begin{array}{c}\pi_1^*(\odot^2N_1^*)\otimes\pi_2^*(L_2^*)^2\\ +\\ \pi_1^*(L_1\otimes N_1^*)\otimes\pi_2^*(L_2^*\otimes N_2^*)\\ +\\ \pi_1^*(L_1^*)^2\otimes\pi_2^*(\odot^2N_2^*)\end{array}\longrightarrow\begin{array}{c}\pi_1^*(L_1\otimes\odot^3N_1^*)\otimes\pi_2^*(L_2^*)^2\\ +\\ \pi_1^*(\odot^2N_1^*)\otimes\pi_2^*(L_2^*\otimes N_2^*)\\ +\\ \pi_1^*(L_1^*\otimes N_1^*)\otimes\pi_2^*(\odot N_2^*)\\ +\\ \pi_1^*(L_1^*)^2\otimes\pi_2^*(L_2\otimes\odot^3N_2^*)\end{array}\longrightarrow\mathbb{L}\otimes\odot^3\mathbb{N}^*\longrightarrow 0.$$

Proposition 5.2. *Let $X$ be a compact complex homogeneous-rational manifold and $L$ an ample line bundle on $X$. Then*

$$\mathrm{H}^0(X,TX\otimes L^*)=\begin{cases}\mathbb{C} & \textit{for } (X,L)=(\mathbb{CP}_1,\mathcal{O}(2))\\ \mathbb{C}^{n+1} & \textit{for } (X,L)=(\mathbb{CP}_n,\mathcal{O}(1)),\ n\geqslant 1\\ 0 & \textit{otherwise.}\end{cases}$$

*Proof.* If $\dim X=1$, i.e. $X=\mathbb{CP}^1$, then the statement follows from the isomorphism $TX\simeq\mathcal{O}(2)$.

Asssume now that $\dim X\geqslant 2$. Then, by the Kodaira vanishing theorem, $\mathrm{H}^1(X,L^*)=0$ for any ample line bundle $L$ on $X$. Applying the Künneth formula to the long exact sequence of (22) with $\mathbb{X}=X\times X$ and $\mathbb{L}=\pi_1^*(L)\otimes\pi_2^*(L)$, one easily obtains

$$\begin{aligned}\mathrm{H}^0(\mathbb{X},\mathbb{L}\otimes\odot^2\mathbb{N}^*) &= \mathrm{H}^0(X,N^*)\otimes\mathrm{H}^0(X,N^*)\\ &= \mathrm{H}^0(X,TX\otimes L^*)\otimes\mathrm{H}^0(X,TX\otimes L^*).\end{aligned}$$

On the other hand, by Theorem 3.3,

$$\mathrm{H}^0(\mathbb{X}, \mathbb{L} \otimes \odot^2 \mathbb{N}^*) = \mathfrak{g}^{(1)},$$

where $\mathfrak{g}$ is the irreducible representation of

$$\mathrm{H}^0(\mathbb{X}, \mathbb{L} \otimes \mathbb{N}^*) \simeq \mathbb{C} \oplus \mathrm{H}^0(X, TX) \oplus \mathrm{H}^0(X, TX)$$

in $V \otimes V$ with $V = \mathrm{H}^0(X, L)$. Table 5 implies that such a $\mathfrak{g}^{(1)}$ can be nonzero if and only if $\mathrm{H}^0(X, TX) \simeq \mathfrak{sl}(n+1, \mathbb{C})$ irreducibly represented in $\mathbb{C}^{n+1}$, i.e. $X = \mathbb{CP}_n$. Then the isomorphism $\mathrm{H}^0(X, L) = \mathbb{C}^{n+1}$ implies $L = \mathcal{O}(1)$. Therefore, $\mathrm{H}^0(\mathbb{X}, \mathbb{L} \otimes \odot^2 \mathbb{N}^*)$ with $\mathbb{X} = X \times X$ and $\mathbb{L} = \pi_1^*(L) \otimes \pi_2^*(L)$ vanishes unless $(X, L) = (\mathbb{CP}_n, \mathcal{O}(1))$ which implies that $\mathrm{H}^0(X, TX \otimes L^*)$ vanishes unless $(X, L) = (\mathbb{CP}_n, \mathcal{O}(1))$. Finally, the extension

$$0 \longrightarrow \mathcal{O}(-1) \longrightarrow \mathbb{C}^{n+1} \otimes \mathcal{O}_{\mathbb{CP}_n} \longrightarrow T\mathbb{CP}_n(-1) \longrightarrow 0$$

implies $\mathrm{H}^0(\mathbb{CP}_n, T\mathbb{CP}_n(-1)) = \mathbb{C}^{n+1}$ which completes the proof of Proposition 5.2 when $\dim X \geqslant 2$. □

Corollary 5.3. *Let $X$ be a compact complex homogeneous-rational manifold, $L$ an ample line bundle on $X$ and $N = J^1 L$. Then, for any $k \geqslant 1$,*

$$\mathrm{H}^0\left(X, \odot^k N^*\right) = \begin{cases} \odot^k \mathbb{C}^{n+1} & \text{for } (X, L) = (\mathbb{CP}_n, \mathcal{O}(1)),\ n \geqslant 1 \\ 0 & \text{otherwise.} \end{cases}$$

*Proof.* The statement is true for $(X, L) = (\mathbb{CP}_n, \mathcal{O}(1))$ since $J^1\mathcal{O}(1) = \mathbb{C}^{n+1} \otimes \mathcal{O}_X$.

The case $k = 1$ of the required statement follows immediately from Propostion 5.2 and the extension

$$0 \longrightarrow L^* \longrightarrow N^* \longrightarrow TX \otimes L^* \longrightarrow 0. \tag{24}$$

The latter also implies

$$0 \longrightarrow \odot^k N^* \longrightarrow L \otimes \odot^{k+1} N^* \longrightarrow \odot^{k+1} TX \otimes L^{*k} \longrightarrow 0$$

which in turn implies

$$\mathrm{H}^0(X, \odot^k N^*) \subset \mathrm{H}^0(X, L \otimes \odot^{k+1} N^*).$$

According to Cartan [17] (see also [38] for another proof), the only irreducible complex Lie subalgebras $\mathfrak{g} \subset \mathfrak{gl}(V)$ which have $\mathfrak{g}^{(k)} \neq 0$ for $k \geqslant 3$ are $\mathfrak{gl}(m, \mathbb{C})$, $\mathfrak{sl}(m, \mathbb{C})$, $\mathfrak{sp}(m/2, \mathbb{C})$ and $\mathfrak{sp}(m/2, \mathbb{C}) \oplus \mathbb{C}$ standardly represented in $\mathbb{C}^m$, $m \geqslant 2$. The Borel-Weil data $(X, L)$ associated with these four representations are $(\mathbb{CP}_{m-1}, \mathcal{O}(1))$. Therefore, if $(X, L) \neq (\mathbb{CP}_n, \mathcal{O}(1))$, then, by Theorem 3.3, $\mathrm{H}^0(X, L \otimes \odot^{k+1} N^*) = 0$ for all $k \geqslant 3$. Hence $\mathrm{H}^0(X, \odot^k N^*) = 0$ for all $k \geqslant 3$. This proves our corollary for $k \geqslant 3$.

Asume now that $k=2$ and denote $\tilde{L} := L^2$ and $\tilde{N} := J^1\tilde{L} \simeq L\otimes N$. Then

$$\mathrm{H}^0(X, \tilde{L}\otimes \odot^2\tilde{N}^*) = \mathrm{H}^0(X, \odot^2 N^*).$$

Again, using Theorem 3.3 and Table 5 one concludes that the only irreducibly acting reductive Lie subalgebra $\mathfrak{g}\subset\mathfrak{gl}(V)$ which has $\mathfrak{g}^{(1)}\neq 0$ and whose associated pair $(X,\tilde{L})$ is such that $\tilde{L}$ is a *square* of an ample line bundle on $X$ is $\mathfrak{gl}(n,\mathbb{C})$ represented in $\odot^2\mathbb{C}^n$ with the associated Borel-Weil data $(\mathbb{CP}_{n-1}, \mathcal{O}(2))$. Therefore, $\mathrm{H}^0(X,\odot^2 N^*) = 0$ for all $(X,L)\neq(\mathbb{CP}_n,\mathcal{O}(1))$. The proof is complete. □

5.3. *The case* $\mathbb{X} = X_1\times X_2$ *with* $\dim X_i \geqslant 2$. The long exact sequence of (22) implies

$$\mathrm{H}^0\left(\mathbb{X}, \mathbb{L}\otimes\odot^2\mathbb{N}^*\right) = \mathrm{H}^0(X_1, N_1^*)\otimes \mathrm{H}^0(X_2,N_2^*) \tag{25}$$

while the long exact sequence of (23) contains the following piece

$$0\longrightarrow \begin{array}{c} \mathrm{H}^0\left(X_1,\odot^2 N_1^*\right)\otimes \mathrm{H}^1\left(X_2, TX_2\otimes L_2^{*2}\right) \\ + \\ \mathrm{H}^0\left(X_2,\odot^2 N_2^*\right)\otimes \mathrm{H}^1\left(X_1, TX_1\otimes L_1^{*2}\right)\end{array} \longrightarrow \mathrm{H}^1\left(\mathbb{X},\mathbb{L}\otimes\odot^3\mathbb{N}^*\right)\longrightarrow$$

$$\longrightarrow \mathrm{H}^1\left(X_1, TX_1\otimes L_1^{*2}\right)\otimes \mathrm{H}^1\left(X_2, TX_2\otimes L_2^{*2}\right)\longrightarrow \dots . \tag{26}$$

LEMMA 5.4. *Let* $G_i\subset \mathrm{GL}(V_i)$, $i=1,2$, *be an irreducible complex semisimple matrix Lie group such that the associated Borel-Weil data* $(X_i,L_i)$ *satisfy* $\dim X_i\geqslant 2$. *Then* $G = \mathrm{T}_{\mathbb{C}}\cdot G_1\cdot G_2\subset \mathrm{GL}(V_1\otimes V_2)$ *can have* $K(\mathfrak{g})\neq 0$ *only if each* $G_i$ *is isomorphic to one of the following representations:*

| Group: | $\mathrm{SL}(n,\mathbb{C})$ | $\mathrm{Sp}(m,\mathbb{C})$ | $\mathrm{SO}(p,\mathbb{C})$ | $\mathrm{G}_2$ | $\mathrm{Spin}(7,\mathbb{C})$ |
|---|---|---|---|---|---|
| Representation space: | $\mathbb{C}^n$ | $\mathbb{C}^{2m}$ | $\mathbb{C}^p$ | $\mathbb{C}^7$ | $\mathbb{C}^8$ |

*where* $n\geqslant 3$, $m\geqslant 2$ *and* $p\geqslant 4$.

*Proof.* $K(\mathfrak{g})\neq 0$ if and only if $\mathfrak{g}^{(1)}\neq 0$ and/or $\mathrm{H}^{1,2}(\mathfrak{g})\neq 0$. By Theorem 3.3, Corollary 5.3 for $k=1$ and (25), one has

$$\mathfrak{g}^{(1)} = \mathrm{H}^0\left(\mathbb{X},\mathbb{L}\otimes\odot^2\mathbb{N}^*\right) = \begin{cases} \mathbb{C}^{n_1+1}\otimes\mathbb{C}^{n_2+1} & \text{for } (X_i,L_i) = (\mathbb{CP}_{n_i},\mathcal{O}(1)),\ n_i\geqslant 1 \\ 0 & \text{otherwise.}\end{cases}$$

On the other hand, a glance at Table 3 shows that $\mathrm{H}^1(X, TX\otimes L^{*2})\neq 0$ only if $(X,L) = (Q_n, j^*\mathcal{O}(1))$ where $Q_n$ is the $n$-dimensional quadric and $j: Q_n\hookrightarrow\mathbb{CP}_{n+1}$ is its standard embedding. Then Theorem 3.3, Corollary 5.3 for $k=2$ and (26) imply that $\mathrm{H}^{1,2}(\mathfrak{g})$ can be nonzero only if each pair $(X_i,L_i)$ is isomorphic either to $(\mathbb{CP}_n,\mathcal{O}(1))$ or to $(Q_n, j^*\mathcal{O}(1))$.

These observations combined with Fact 3.1 imply that $K(\mathfrak{g})$ can be nonzero only for representations listed in Lemma 5.4. □

*Example* 1. Let $G$ be the representation of $\mathbb{Z}_{\mathbb{C}} \cdot \mathrm{SL}(m,\mathbb{C}) \cdot \mathrm{SL}(n,\mathbb{C})$ in the vector space $V = V_m \otimes V_n$ where $V_m$ and $V_n$ are $m$- and, respectively, $n$-dimensional complex vector spaces with $m, n \geqslant 3$. The associated Borel-Weil datum $(\mathbb{X}, \mathbb{L})$ is $(\mathbb{CP}_{m-1} \times \mathbb{CP}_{n-1}, \pi_1^*(\mathcal{O}(1)) \otimes \pi_2^*(\mathcal{O}(1)))$, implying

$$\mathfrak{g}^{(1)} = \mathrm{H}^0(\mathbb{X}, \mathbb{L} \otimes \odot^2 \mathbb{N}^*) = V^*, \qquad \mathrm{H}^{1,2}(\mathfrak{g}) = \mathrm{H}^1(\mathbb{X}, \mathbb{L} \otimes \odot^3 \mathbb{N}^*) = 0.$$

Therefore, $K(\mathfrak{g}) = \partial(\mathfrak{g}^{(1)} \otimes V^*) \simeq V^* \otimes V^*$. Denoting typical elements of $V$, $V_m$ and $V_n$ by $v^a$, $v^A$ and $v^{\dot{A}}$ respectively,[1] and identifying $v^a \in V$ with its image $v^{A\dot{A}}$ under the isomorphism $V = V_m \otimes V_n$, one may write a typical element $R_{abc}{}^d \in K(\mathfrak{g}) \subset \Lambda^2 V^* \otimes V^* \otimes V$ as

$$R_{abc}{}^d \equiv R_{A\dot{A}B\dot{B}C\dot{C}}{}^{D\dot{D}} = \left[\delta_A^D Q_{B\dot{B}C\dot{A}} - \delta_B^D Q_{A\dot{A}C\dot{B}}\right]\delta_{\dot{C}}^{\dot{D}} + \left[\delta_{\dot{A}}^{\dot{D}} Q_{B\dot{B}A\dot{C}} - \delta_{\dot{B}}^{\dot{D}} Q_{A\dot{A}B\dot{C}}\right]\delta_C^D \tag{27}$$

for some $Q_{ab} \equiv Q_{A\dot{A}B\dot{B}} \in V^* \otimes V^*$. Therefore, a torsion-free connection $\nabla$ on an $mn$-dimensional manifold $M$ with holonomy in $\mathbb{Z}_{\mathbb{C}} \cdot \mathrm{SL}(m,\mathbb{C}) \cdot \mathrm{SL}(n,\mathbb{C})$ has at an arbitrary point $x \in M$ the curvature tensor of the form (27) for some $Q_{ab}(x) \in \Omega_x M \otimes \Omega_x M$. It is not hard to show that the second Bianchi identities for $\nabla$,

$$\nabla_e R_{abc}{}^d + \nabla_b R_{eac}{}^d + \nabla_a R_{bec}{}^d = 0,$$

imply

$$\begin{aligned} 0 \;=\; & m\left(\nabla_{A\dot{A}} Q_{B\dot{B}C\dot{C}} - \nabla_{B\dot{B}} Q_{A\dot{A}C\dot{C}}\right) + n\left(\nabla_{C\dot{C}} Q_{A\dot{A}B\dot{B}} - \nabla_{A\dot{A}} Q_{C\dot{C}B\dot{B}}\right) \\ & + \left(\nabla_{B\dot{B}} Q_{A\dot{C}C\dot{A}} - \nabla_{A\dot{C}} Q_{B\dot{B}C\dot{A}}\right) + \left(\nabla_{B\dot{A}} Q_{C\dot{C}A\dot{B}} - \nabla_{C\dot{C}} Q_{B\dot{A}A\dot{B}}\right). \end{aligned} \tag{28}$$

*Example* 2. Keeping notation of the preceding paragraph, we consider a subgroup $G_o \subset G$ which is $\mathbb{Z}_{\mathbb{C}} \cdot \mathrm{SL}(m,\mathbb{C}) \cdot \mathrm{SO}(n,\mathbb{C})$ represented in $V = V_m \otimes V_n$ with $m, n \geqslant 3$. The $G_o$-module $K(\mathfrak{g}_o)$ is a subset of $K(\mathfrak{g})$ consisting of all elements $R_{abc}{}^d$ satisfying

$$R_{A\dot{A}B\dot{B}C\dot{C}}{}^{D\dot{D}}\, g_{\dot{E}\dot{D}} + R_{A\dot{A}B\dot{B}C\dot{E}}{}^{D\dot{D}}\, g_{\dot{C}\dot{D}} = 0,$$

where $g_{\dot{E}\dot{D}} \in \odot^2 V_n^*$ is the $\mathrm{SO}(n,\mathbb{C})$-invariant quadratic form. Substituting (27) into the above equation, one obtains after elementary algebraic manipulations that

$$Q_{A\dot{A}B\dot{B}} = P_{AB}\, g_{\dot{A}\dot{B}}$$

[1] One may view indices of the type $a$, $A$ or $\dot{A}$ as refering to some fixed basis in a relevant vector space or, alternatively, as abstract labels providing us with a transparent notation for such basic tensor operations as (anti)symmetrization, contraction, etc. (cf. [7], [35]).

for some symmetric tensor $P_{AB} \in \odot^2 V_m^*$. (Another way to obtain this result is to note that the Borel-Weil datum $(\mathbb{X}, \mathbb{L})$ associated to $(G_o, V)$ is $(\mathbb{CP}_{m-1} \times Q_{n-1}, \pi_1^*(\mathcal{O}(1)) \otimes \pi_2^*(j^*\mathcal{O}(1)))$ implying $\mathfrak{g}_o^{(1)} = \mathrm{H}^0(\mathbb{X}, \mathbb{L} \otimes \odot^2 \mathbb{N}^*) = 0$ and $K(\mathfrak{g}_o) = \mathrm{H}^1(\mathbb{X}, \mathbb{L} \otimes \odot^3 \mathbb{N}^*) = \odot^2 V_m^* \otimes C \subset V^* \otimes V^*$, where $C$ is the 1-dimensional subspace of $\odot^2 V_n$ spanned by $g_{\dot{A}\dot{B}}$.) Then the second Bianchi identities (28) imply $\nabla_a Q_{bc} = 0$ which in turn implies $\nabla_m R_{abc}{}^d = 0$. These arguments imply essentially the following:

Lemma 5.5. *Let $G$ be the irreducible representation of a subgroup of $\mathrm{GL}(m, \mathbb{C}) \cdot \mathrm{SO}(n, \mathbb{C})$ in the $mn$-dimensional vector space $V_m \otimes V_n$. If $m, n \geqslant 3$, then $K^1(\mathfrak{g}) = 0$.*

*Example* 3. Keeping notation of Example 1, consider a subgroup $G_s \subset G$ which is $\mathbb{Z}_{\mathbb{C}} \cdot \mathrm{SL}(m, \mathbb{C}) \cdot \mathrm{Sp}(n, \mathbb{C})$ represented in $V = V_m \otimes V_{2n}$ with $m \geqslant 3$, $n \geqslant 2$, and note that $K(\mathfrak{g}_s)$ is a subset of $K(\mathfrak{g})$ consisting of all elements $R_{abc}{}^d$ satisfying

$$R_{A\dot{A}B\dot{B}C\dot{C}}{}^{D\dot{D}}\, \varepsilon_{\dot{E}\dot{D}} - R_{A\dot{A}B\dot{B}C\dot{E}}{}^{D\dot{D}}\, \varepsilon_{\dot{C}\dot{D}} = 0,$$

where $\varepsilon_{\dot{E}\dot{D}} \in \Lambda^2 V_{2n}^*$ is the $\mathrm{Sp}(n, \mathbb{C})$-invariant symplectic form. Substituting (27) into this equation, one easily finds

$$Q_{A\dot{A}B\dot{B}} = S_{AB}\, \varepsilon_{\dot{A}\dot{B}}$$

for some antisymmetric tensor $S_{AB} \in \Lambda^2 V_m$. Then the second Bianchi identities (28) imply $\nabla_a Q_{bc} = 0$ which in turn implies $\nabla_m R_{abc}{}^d = 0$. We may summarize these arguments as follows.

Lemma 5.6. *Let $G$ be the irreducible representation of a subgroup of $\mathrm{GL}(m, \mathbb{C}) \cdot \mathrm{Sp}(n, \mathbb{C})$ in the $2mn$-dimensional vector space $V_m \otimes V_{2n}$. If $m \geqslant 3$, $n \geqslant 2$, then $K^1(\mathfrak{g}) = 0$.*

An immediate corollary of Lemmas 5.4–5.6 is the following:

Proposition 5.7. *Let $G_i \subset \mathrm{GL}(V_{n_i})$, $i = 1, 2$, be an irreducible complex semisimple matrix Lie group such that the associated Borel-Weil data $(X_i, L_i)$ satisfy $\dim X_i \geqslant 2$. Then $G = \mathbb{Z}_{\mathbb{C}} \cdot G_1 \cdot G_2 \subset \mathrm{GL}(V_{n_1} \otimes V_{n_2})$ can have $K^1(\mathfrak{g}) \neq 0$ only if $G_1 = \mathrm{SL}(n_1, \mathbb{C})$ and $G_2 = \mathrm{SL}(n_2, \mathbb{C})$.*

5.4. *The case $\mathbb{X} = X \times \mathbb{CP}_1$ with $\dim X \geqslant 2$.* Any ample line bundle on $\mathbb{X}$ is of the form $\mathbb{L} = \pi_1^*(L) \otimes \pi_2^*(\mathcal{O}(k))$ for some ample line bundle $L \to X$ and $k \geqslant 1$. We denote in this subsection $V_n := \mathrm{H}^0(X, L)$, $V_2 := \mathrm{H}^0(\mathbb{CP}_1, \mathcal{O}(1))$, $N := J^1 L$, and $\mathfrak{g}$ stands for the Lie algebra $\mathrm{H}^0(X, L \otimes N^*) + \mathfrak{sl}(2, \mathbb{C})$ represented in $V = V_n \otimes V_2$.

If $k \geqslant 2$, then the associated matrix Lie group $G = \exp(\mathfrak{g})$ is an irreducible matrix subgroup of either $\mathrm{GL}(n, \mathbb{C})\mathrm{SO}(p, \mathbb{C})$ represented in $\mathbb{C}^{np}$ for some $n, p \geqslant$ 3 or $\mathrm{GL}(n,\mathbb{C})\mathrm{Sp}(q,\mathbb{C})$ represented in $\mathbb{C}^{2nq}$ for some $n \geqslant 3, q \geqslant 2$. Then, by Lemmas 5.5 and 5.6, $K^1(\mathfrak{g}) = 0$.

So we may assume that $k = 1$.

PROPOSITION 5.8. *Let $(X, L)$ be a pair consisting of a compact complex homogeneous-rational manifold $X$ and an ample line bundle $L \to X$. If $\dim X \geqslant 2$, then*

$$\mathfrak{g}^{(1)} = \mathrm{H}^0(X, N^*) \otimes V_2^*$$

*and there is an exact sequence of $\mathfrak{g}$-modules*

$$0 \longrightarrow \mathrm{H}^{1,2}(\mathfrak{g}) \longrightarrow \begin{array}{c} \mathrm{H}^0(X, \odot^3 TX \otimes L^{*2}) \otimes \Lambda^2 V_2^* \\ + \\ \mathrm{H}^1(X, TX \otimes L^{*2}) \otimes \odot^2 V_2^* \end{array} \longrightarrow \mathrm{H}^1(X, TX \otimes L^*) \otimes V^* \otimes V_2^*.$$

*Proof.* Since $\dim X \geqslant 2$, the Kodaira vanishing theorem implies $\mathrm{H}^1(X, L^*) = 0$. Then the long exact sequence of

$$0 \longrightarrow L^* \longrightarrow N^* \longrightarrow TX \otimes L^* \longrightarrow 0$$

implies $\mathrm{H}^1(X, N^*) = \mathrm{H}^1(X, TX \otimes L^*)$, while the long exact sequence of (22) with $(X_1, L_1) = (X, L)$ and $(X_2, L_2) = (\mathbb{CP}_1, \mathcal{O}(1))$ implies

$$\mathrm{H}^0(\mathbb{X}, \mathbb{L} \otimes \odot^2 \mathbb{N}^*) = \mathrm{H}^0(X, N^*) \otimes V_2^*,$$

$$\mathrm{H}^1(\mathbb{X}, \mathbb{L} \otimes \odot^2 \mathbb{N}^*) = \mathrm{H}^1(X, N^*) \otimes V_2^* = \mathrm{H}^1(X, TX \otimes L^*) \otimes V_2^*.$$

Analogously, the long exact sequence of (23) implies

$$0 \longrightarrow \mathrm{H}^0(\mathbb{X}, \mathbb{L} \otimes \odot^3 \mathbb{N}^*) \longrightarrow \begin{array}{c} \mathrm{H}^0(X, \odot^2 N^*) \otimes \Lambda^2 V_2^* \\ + \\ 0 \end{array} \longrightarrow \begin{array}{c} \mathrm{H}^0(X, L \otimes \odot^3 N^*) \otimes \Lambda^2 V_2^* \\ + \\ 0 \end{array} \longrightarrow$$

$$\longrightarrow \mathrm{H}^1(\mathbb{X}, \mathbb{L} \otimes \odot^3 \mathbb{N}^*) \longrightarrow \begin{array}{c} \mathrm{H}^1(X, \odot^2 N^*) \otimes \Lambda^2 V_2^* \\ + \\ 0 \end{array} \longrightarrow \begin{array}{c} \mathrm{H}^1(X, L \otimes \odot^3 N^*) \otimes \Lambda^2 V_2^* \\ + \\ \mathrm{H}^2(X, L^* \otimes N^*) \otimes \odot^2 V_2^* \end{array} \longrightarrow \ldots.$$

Comparing this with the long exact sequence of

$$0 \longrightarrow \odot^2 N^* \longrightarrow L \otimes \odot^3 N^* \longrightarrow \odot^3 TX \otimes L^{*2} \longrightarrow 0 \tag{29}$$

one obtains

$$\mathrm{H}^1(\mathbb{X}, \mathbb{L} \otimes \odot^3 \mathbb{N}^*) = \mathrm{H}^0(X, \odot^3 TX \otimes L^{*2}) \otimes \Lambda^2 V_2^* + \mathrm{H}^1(X, TX \otimes L^{*2}) \otimes \odot^2 V_2^*.$$

Then Theorem 3.3 implies the desired result. □

LEMMA 5.9. *Let* $G \subset \mathrm{GL}(V)$ *be an irreducible complex representation of a complex reductive Lie group such that the associated Borel-Weil data are of the form* $(\mathbb{X} = X \times \mathbb{CP}_1, \mathbb{L})$ *with* $\dim X \geqslant 2$ *and* $\mathrm{H}^1(X, TX \otimes L^{*2}) = 0$. *Then* $K^1(\mathfrak{g}) \neq 0$ *if only if* $G$ *is* $\mathrm{T}_{\mathbb{C}} \cdot \mathrm{SL}(n, \mathbb{C}) \cdot \mathrm{SL}(2, \mathbb{C})$ *represented in* $V = V_n \otimes V_2$.

*Proof.* If $(X, L) = (\mathbb{CP}_n, \mathcal{O}(1))$, then $G$ has both modules $K(\mathfrak{g})$ and $K^1(\mathfrak{g})$ nonzero only if it is $\mathrm{T}_{\mathbb{C}} \cdot \mathrm{SL}(n, \mathbb{C}) \cdot \mathrm{SL}(2, \mathbb{C})$ represented in $V_n \otimes V_2$.

Assume now that $(X, L) \neq (\mathbb{CP}_n, \mathcal{O}(1))$. Then, by Proposition 5.9,

$$K(\mathfrak{g}) \subset \mathrm{H}^0(X, \odot^3 TX \otimes L^{*2}) \otimes \Lambda^2 V_2^*.$$

Since $\mathrm{H}^0(X, L \otimes \odot^3 N^*) = 0$, the long exact sequence of (29) implies $\mathrm{H}^0(X, \odot^3 TX \otimes L^{*2}) \subset H^1(X, \odot^2 N^*)$, while the exact sequence

$$0 \longrightarrow \odot^2 N^* \longrightarrow N^* \otimes V_n^* \longrightarrow \Lambda^2 V_n^*$$

implies $H^1(X, \odot^2 N^*) \subset \Lambda^2 V_n^*$. Thus $K(\mathfrak{g}) \subset \Lambda^2 V_n^* \otimes \Lambda^2 V_2^*$ which easily implies $K^1(\mathfrak{g}) = 0$. □

LEMMA 5.10. *Let* $G \subset \mathrm{GL}(V)$ *be an irreducible complex representation of a complex reductive Lie group such that the associated Borel-Weil data are of the form* $(\mathbb{X} = X \times \mathbb{CP}_1, \mathbb{L})$ *with* $\dim X \geqslant 2$ *and* $\mathrm{H}^1(X, TX \otimes L^{*2}) \neq 0$. *Then* $K^1(\mathfrak{g}) \neq 0$ *only if* $G$ *is either* $\mathrm{T}_{\mathbb{C}} \cdot \mathrm{SL}(n, \mathbb{C}) \cdot \mathrm{SL}(2, \mathbb{C})$ *or* $\mathrm{T}_{\mathbb{C}} \cdot \mathrm{SO}(n, \mathbb{C}) \cdot \mathrm{SL}(2, \mathbb{C})$, *both represented in* $V_n \otimes V_2$.

*Proof.* It follows from Table 3 that $\mathrm{H}^1(X, TX \otimes L^{*2}) \neq 0$ only if $(X, L) = (Q_n, j^*\mathcal{O}(1))$ where $Q_n$ is the $n$-dimensional quadric and $j : Q_n \hookrightarrow \mathbb{CP}_{n+1}$ is its standard embedding. This together with Fact 12 imply that $G$ must be of the form $\mathbb{Z}_{\mathbb{C}} \cdot H \cdot \mathrm{SL}(2, \mathbb{C}) \subset \mathfrak{gl}(V_n \otimes V_2)$ where $H$ is one of the following representations:

| Group $H$: | $\mathrm{SO}(n, \mathbb{C})$ | $\mathrm{G}_2$ | $\mathrm{Spin}(7, \mathbb{C})$ |
|---|---|---|---|
| Representation space $V_n$: | $\mathbb{C}^n$ | $\mathbb{C}^7$ | $\mathbb{C}^8$ |

Since the Borel-Weil data associated to $G = \mathbb{Z}_{\mathbb{C}} \cdot \mathrm{SO}(n, \mathbb{C}) \cdot \mathrm{SL}(2, \mathbb{C})$ represented in $V_n \otimes V_2$ are $(Q_{n-1} \times \mathbb{CP}_1, \pi_1^*(j^*\mathcal{O}(1)) \otimes \pi_2^*(\mathcal{O}(1))$, one has $\mathfrak{g}^{(1)} = \mathrm{H}^0(\mathbb{X}, \mathbb{L} \otimes \odot^2 \mathbb{N}^*) = 0$ and

$$K(\mathfrak{g}) = \mathrm{H}^1(\mathbb{X}, \mathbb{L} \otimes \odot^3 \mathbb{N}^*) = \Lambda^2 V_n^* \otimes \Lambda^2 V_2^* + C \otimes \odot^2 V_2^*$$

where $C$ is the 1-dimensional subspace of $\odot^2 V_n$ spanned by the $\mathrm{SO}(n, \mathbb{C})$-invariant metric $g_{AB}$. Then a generic element of $K(\mathfrak{g})$ must be of the form (cf.

[20])

$$R_{abc}{}^{d} = [\varepsilon_{\dot{A}\dot{B}} g^{DE}(g_{AB}S_{CE} + g_{AC}S_{BE} + g_{BC}S_{AE} - g_{AE}S_{BC} - g_{BE}S_{AC}) + \Phi_{\dot{A}\dot{B}}(g_{BC}\delta_A^D - g_{AC}\delta_B^D]\delta_{\dot{C}}^{\dot{D}} + [g_{AB}\varepsilon_{\dot{A}\dot{B}}\varepsilon^{\dot{D}\dot{E}}\Phi_{\dot{C}\dot{E}} - S_{AB}(\varepsilon_{\dot{B}\dot{C}}\delta_{\dot{A}}^{\dot{D}} + \varepsilon_{\dot{A}\dot{C}}\delta_{\dot{A}}^{\dot{B}})]\delta_C^D \tag{30}$$

for some $S_{AB} \in \Lambda^2 V_n^*$ and $\Phi_{\dot{A}\dot{B}} \in \odot^2 V_2^*$.

Let $g \subset \mathfrak{gl}(V)$ be the Lie algebra of the representation of $\mathbb{Z}_{\mathbb{C}} \cdot \mathrm{G}_2 \cdot \mathrm{SL}(2,\mathbb{C})$ (resp. $\mathbb{Z}_{\mathbb{C}} \cdot \mathrm{Spin}(7,\mathbb{C}) \cdot \mathrm{SL}(2,\mathbb{C})$) in $V = \mathbb{C}^7 \otimes \mathbb{C}^2$ (resp. in $V = \mathbb{C}^8 \otimes \mathbb{C}^2$). It is a proper matrix subalgebra of the Lie algebra $\mathfrak{g}$ of the representation of $\mathbb{Z}_{\mathbb{C}} \cdot \mathrm{SO}(7,\mathbb{C}) \cdot \mathrm{SL}(2,\mathbb{C})$ (resp. $\mathbb{Z}_{\mathbb{C}} \cdot \mathrm{SO}(8,\mathbb{C}) \cdot \mathrm{SL}(2,\mathbb{C})$) in $V = \mathbb{C}^7 \otimes \mathbb{C}^2$ (resp. in $V = \mathbb{C}^8 \otimes \mathbb{C}^2$). Then

$$K(g) \subset K(\mathfrak{g}) = \Lambda^2 V_n^* \otimes \Lambda^2 V_2^* + C \otimes \odot^2 V_2^*$$

and $K^1(g) \subset K^1(\mathfrak{g})$. We claim

$$K(g) \subset \Lambda^2 V_n^* \otimes \Lambda^2 V_2^*. \tag{31}$$

If not, then a typical element $R_{abc}{}^{d} \in K(g)$ contains a nonzero term $\Phi_{\dot{A}\dot{B}}(g_{BC}\delta_A^D - g_{AC}\delta_B^D)\delta_{\dot{C}}^{\dot{D}}$ which easily implies that the image of the map

$$\Lambda^2 V \longrightarrow g$$

defined by $R_{abc}{}^{d} \in g \otimes \Lambda^2 V^*$ contains $\Lambda^2 V_n^* \simeq \mathfrak{so}(n,\mathbb{C})$. This contradicts the fact that $g$ is a proper subalgebra of $\mathfrak{g}$.

Finally, it is straightforward to check that the inclusion (31) implies that $K^1(g) = 0$. □

PROPOSITION 5.11. *Let $G \subset \mathrm{GL}(V)$ be an irreducible complex representation of a complex reductive Lie group such that the associated Borel-Weil data are of the form $(\mathbb{X} = X \times \mathbb{CP}_1, \mathbb{L})$ with $\dim X \geqslant 2$. Then $K^1(\mathfrak{g}) \neq 0$ if and only if $G$ is either $\mathrm{T}_{\mathbb{C}} \cdot \mathrm{SL}(n,\mathbb{C}) \cdot \mathrm{SL}(2,\mathbb{C})$ or $\mathrm{SO}(n,\mathbb{C}) \cdot \mathrm{SL}(2,\mathbb{C})$, both represented in $V_n \otimes V_2$.*

*Proof.* By Lemmas 5.9 and 5.10, one has only to rule out the case $\mathbb{C}^* \cdot \mathrm{SO}(n,\mathbb{C}) \cdot \mathrm{SL}(2,\mathbb{C})$. But this follows from $R_{abc}{}^{c} = 0$ which itself follows from (30). □

5.5. *The case* $\mathbb{X} = \mathbb{CP}_1 \times \mathbb{CP}_1$. This is the case of $\mathbb{Z}_{\mathbb{C}} \cdot \mathrm{SL}(2,\mathbb{C}) \cdot \mathrm{SL}(2,\mathbb{C})$ represented in $\odot^m V_2 \otimes \odot^n V_2$. In the context of the holonomy classification, this class of representations was studied in [22] and [32] where the following result was established by two different methods.

PROPOSITION 5.12. *Let $G \subset \mathrm{GL}(V)$ be an irreducible complex representation of a complex reductive Lie group such that the associated Borel-Weil* *data* $(\mathbb{X}, \mathbb{L})$ *has* $\mathbb{X} = \mathbb{CP}_1 \times \mathbb{CP}_1$. *Then* $K^1(\mathfrak{g}) \neq 0$ *only if* $G$ *is either the representation of* $\mathrm{T}_{\mathbb{C}} \cdot \mathrm{SO}(4, \mathbb{C})$ *in* $\mathbb{C}^4$ *or the representation of* $\mathrm{T}_{\mathbb{C}} \cdot \mathrm{SO}(3, \mathbb{C}) \cdot \mathrm{SL}(2, \mathbb{C})$ *in* $\mathbb{C}^6$.

In fact, for the above representation, $K(\mathbb{C} + \mathfrak{so}(3, \mathbb{C}) + \mathfrak{sl}(2, \mathbb{C})) = K(\mathfrak{so}(3, \mathbb{C}) + \mathfrak{sl}(2, \mathbb{C}))$ which means that $\mathbb{C}^* \cdot \mathrm{SO}(3, \mathbb{C}) \cdot \mathrm{SL}(2, \mathbb{C})$ can not occur as the holonomy of a torsion-free affine connection.

5.6. *Proof of Theorem* 5.1. Let $G \subset \mathrm{GL}(V)$ be an irreducible complex representation of a complex reductive Lie group which can be represented as a tensor product of two or more non-Abelian complex representations. Then, by Propositions 5.7, 5.11 and 5.12, $G$ may occur as the holonomy of a nonmetric torsion-free affine connection only if it is either $\mathbb{Z}_{\mathbb{C}} \cdot \mathrm{SL}(m, \mathbb{C}) \cdot \mathrm{SL}(n, \mathbb{C})$ represented in $\mathbb{C}^m \otimes \mathbb{C}^n$ for $m, n \geqslant 2$, or $\mathrm{SO}(l, \mathbb{C}) \cdot \mathrm{SL}(2, \mathbb{C})$ represented in $\mathbb{C}^l \otimes \mathbb{C}^2$ for $l \geqslant 3$.

Let $\rho : G \to \mathrm{GL}(V)$ be an irreducible representation of a real reductive Lie group $G$ in a real vector space $V$ and let $\rho : \mathfrak{g} \to \mathfrak{gl}(V)$ be the associated real irreducible representation of the Lie algebra $\mathfrak{g}$ of $G$. The latter defines naturally a complex representation $\rho_{\mathbb{C}} : \mathfrak{g}_{\mathbb{C}} \to \mathfrak{gl}(V_{\mathbb{C}})$, where $\mathfrak{g}_{\mathbb{C}} = \mathfrak{g} \otimes \mathbb{C}$ and $V_{\mathbb{C}} = V \otimes \mathbb{C}$. Then two situations may arise [25]:

(i) The complex representation $\rho_{\mathbb{C}} : \mathfrak{g}_{\mathbb{C}} \to \mathfrak{gl}(V_{\mathbb{C}})$ is irreducible; in this case we denote $\rho_{\mathbb{C}}$ by $\tilde{\rho}$;

(ii) There is a complex vector space $W_{\mathbb{C}}$ and an irreducible complex representation $\rho' : \mathfrak{g} \to \mathfrak{gl}(W_{\mathbb{C}})$ such that $V$ is the underlying real vector space of $W_{\mathbb{C}}$ and $\rho$ is the composition $\rho : \mathfrak{g} \xrightarrow{\rho'} \mathfrak{gl}(W_{\mathbb{C}}) \longrightarrow \mathfrak{gl}(V)$, where the second arrow is the natural inclusion of the algebra of all complex automorphisms of $V$ into the algebra of all real automorphisms of $V$. Then the $\mathfrak{g}_{\mathbb{C}}$-module $V_{\mathbb{C}}$ splits as a direct sum of two irreducible $\mathfrak{g}_{\mathbb{C}}$-submodules $W_{\mathbb{C}} + W_{\mathbb{C}}$ and we denote by $\tilde{\rho} : \mathfrak{g}_{\mathbb{C}} \to \mathfrak{gl}(W_{\mathbb{C}})$ the restriction of $\rho_{\mathbb{C}}$ to one of these.

In both cases, the $\mathfrak{g}$-modules $K(\rho(\mathfrak{g}))$ and $K^1(\rho(\mathfrak{g}))$ are subsets of $K(\tilde{\rho}(\mathfrak{g}_{\mathbb{C}}))$ and $K^1(\tilde{\rho}(\mathfrak{g}_{\mathbb{C}}))$ respectively. In particular, if $K(\rho(\mathfrak{g}))$ and $K^1(\rho(\mathfrak{g}))$ are nonzero, then $K(\tilde{\rho}(\mathfrak{g}_{\mathbb{C}}))$ and $K^1(\tilde{\rho}(\mathfrak{g}_{\mathbb{C}}))$ are nonzero as well.

Assume now that the semisimple part of $\mathfrak{g}$ has at least two non-Abelian ideals. Then the Borel-Weil data associated to the irreducible matrix subalgebra $\tilde{\rho}(\mathfrak{g}_{\mathbb{C}})$ must be of the form $(\mathbb{X}, \mathbb{L}) = (X_1 \times X_2, \pi_1^*(L_1) \otimes \pi_2^*(L_2))$ for some compact complex homogeneous-rational manifolds $X_1$ and $X_2$ and ample line bundles $L_1 \to X_1$ and $L_2 \to X_2$.

Then, by Propositions 5.7, 5.11 and 5.12, $\rho(G) \subset \mathrm{GL}(V)$ must be of the form $\mathbb{Z}_{\mathbb{C}} \cdot G_1 \cdot G_2$, where $\mathbb{Z}_{\mathbb{C}}$ is a connected real Lie subgroup of $\mathbb{C}^*$ and $G_i \subset \mathrm{GL}(V_i)$, $i = 1, 2$, is one of the following real matrix groups

| Group $G_i$: | $\mathrm{SL}(n,\mathbb{C})$ | $\mathrm{SL}(n,\mathbb{R})$ | $\mathrm{SU}(n)$ | $\mathrm{SL}(m,\mathrm{H})$ |
|---|---|---|---|---|
| Representation space $V_i$: | $\mathbb{R}^{2n}$ | $\mathbb{R}^{n}$ | $\mathbb{R}^{2n}$ | $\mathbb{R}^{4m}$ |
| Group $G_i$: | $\mathrm{SO}(l,\mathbb{C})$ | $\mathrm{SO}(p,q)$ | $\mathrm{SO}(n,\mathrm{H})$ | |
| Representation space $V_i$: | $\mathbb{R}^{2l}$ | $\mathbb{R}^{p+q}$ | $\mathbb{R}^{4n}$ | |

with $n \geqslant 2$, $m \geqslant 1$, $l \geqslant 3$, $p+q \geqslant 3$.

Since we know $K(\tilde{\rho}(\mathfrak{g}_{\mathbb{C}}))$ explicitly, it is straightforward to check that the only combinations $\rho(G) = \mathbb{Z}_{\mathbb{C}} \cdot G_1 \cdot G_2$ which (i) have $K^1(\rho(\mathfrak{g})) \neq 0$ and (ii) have no proper subgroup $G' \subset G$ with $K(\rho(\mathfrak{g}')) = K(\rho(\mathfrak{g}))$ are the ones given in the Table 7. □

## 6. Proof of the Main Theorem

6.1. *Adjoint representations.* In this subsection, we shall rule out the (central extensions of the) adjoint representations of simple Lie groups as candidates for holonomies. In the Borel-Weil language, this is a family of representations which have contact complex homogeneous manifolds as their skies.

Proposition 6.1. *Let $\mathfrak{g}_s$ be a complex simple Lie algebra with* $\mathrm{rank}(\mathfrak{g}_s) \geqslant 3$, *and let* $\mathfrak{g} := \mathbb{C} \oplus \mathfrak{g}_s$. *Consider the embedding $\imath : \mathfrak{g} \hookrightarrow \mathfrak{gl}(\mathfrak{g}_s)$ via the central extension of the adjoint representation. Then $\imath(\mathfrak{g}) \subset \mathfrak{gl}(\mathfrak{g}_s)$ is not a Berger algebra.*

*Proof.* Fix a Cartan decomposition $\mathfrak{g}_s = \mathfrak{t} \oplus \bigoplus_{\alpha\in\Phi} \mathfrak{g}_\alpha$, and fix elements $0 \neq A_\alpha \in \mathfrak{g}_\alpha$. Moreover, denote elements of $\mathfrak{t}$ by $A_0, B_0, \dots$.

Let us suppose that $\imath(\mathfrak{g})$ *is* a Berger group. Then, by Lemma 2.3, $K(\imath(\mathfrak{g}))$ is a nontrivial $\mathfrak{g}$-module. Also, it is evident that the weights of $K(\imath(\mathfrak{g})) \subset \Lambda^2\mathfrak{g}_s^* \otimes \mathfrak{g}$ all lie in the root lattice of $\mathfrak{t}^*$. It follows that there must be a weigth $\rho$ of $K(\imath(\mathfrak{g}))$ which is a root. Let $0 \neq R \in K(\imath(\mathfrak{g}))$ be a corresponding weight element.

We write $R(x,y) = \{x,y\} - \omega(x,y)\mathrm{Id}$ with $\{x,y\} \in \mathfrak{g}_s$ and a 2-form $\omega$. Then the Bianchi identity reads

$$[\{x,y\},z] + [\{y,z\},x] + [\{z,x\},y] = \omega(x,y)z + \omega(y,z)x + \omega(z,x)y. \tag{32}$$

*Step* 1. For all $A_0, B_0 \in \mathfrak{t}$, $\{A_0,B_0\} = 0$ and $\omega(A_0, _) = 0$.

*Proof.* For $(x,y,z) = (A_0, B_0, A_{-\rho})$, (32) reads

$$[\underbrace{\{A_0,B_0\}}_{\in\mathfrak{g}_\rho}, A_{-\rho}] + [\underbrace{\{B_0,A_{-\rho}\}}_{\in\mathfrak{t}}, A_0] + [\underbrace{\{A_{-\rho},A_0\}}_{\in\mathfrak{t}}, B_0]$$

$$= \underbrace{\omega(A_0, B_0)}_{=0} A_{-\rho} + \omega(B_0, A_{-\rho})A_0 + \omega(A_{-\rho}, A_0)B_0,$$

and thus, if $[A_\rho, A_{-\rho}] \notin \mathrm{span}(A_0, B_0)$, then $\{A_0, B_0\} = 0$ and $\omega(A_0, A_{-\rho}) = 0$. Since $\mathrm{rank}(\mathfrak{g}) \geqslant 3$ and $\omega$ has weight $\rho$, the assertion follows. □

*Step* 2. For each root $\alpha \neq -\rho$, there is a constant $c_\alpha$ such that $\{A_0, A_\alpha\} = c_\alpha(\alpha + \rho)(A_0)A_{\alpha+\rho}$ for all $A_0 \in \mathfrak{t}$. We set $c_\alpha = 0$ if $\alpha + \rho$ is *not* a root.

*Proof.* The assertion is obvious if $\alpha + \rho$ is not a root. If $\alpha + \rho$ *is* a root, then from Step 1 and (32) for $(x, y, z) = (A_0, B_0, A_\alpha)$ we get

$$\begin{aligned} 0 &= [\underbrace{\{B_0, A_\alpha\}}_{\in \mathfrak{g}_{\alpha+\rho}}, A_0] + [\underbrace{\{A_\alpha, A_0\}}_{\in \mathfrak{g}_{\alpha+\rho}}, B_0] \\ &= -(\alpha + \rho)(A_0)\{B_0, A_\alpha\} + (\alpha + \rho)(B_0)\{A_0, A_\alpha\}. \end{aligned}$$

It follows that for all $A_0 \in (\alpha + \rho)^\perp$ we have $\{A_0, A_\alpha\} = 0$, and from there the claim follows.

*Step* 3. For each root $\alpha \neq -\rho$ we have $\{A_0, A_\alpha\} = 0$; moreover, $\omega = 0$.

*Proof.* The first claim is obvious if $\alpha + \rho$ is not a root; thus we assume the contrary and let $\beta := -(\alpha + \rho)$. Then, using Step 2 and $\{A_\alpha, A_\beta\} \in \mathfrak{t}$, we get from (32) for $(x, y, z) = (A_0, A_\alpha, A_\beta)$ that

$$-c_\alpha \beta(A_0)[A_{-\beta}, A_\beta] + c_\beta \alpha(A_0)[A_{-\alpha}, A_\alpha] = \omega(A_\alpha, A_\beta)A_0.$$

Since $\alpha \neq \pm\rho$, it follows that $\alpha$ and $\beta$ are independent, and then we prove the claim using the fact that $\mathrm{rank}(\mathfrak{g}) \geqslant 3$.

*Step* 4. $\{A_\alpha, A_\beta\} = 0$ for all roots $\alpha, \beta \neq -\rho$.

*Proof.* If $\alpha + \beta + \rho \neq 0$ then $\{A_\alpha, A_\beta\} = 0$ follows from the previous steps by application of (32) to $(x, y, z) = (A_0, A_\alpha, A_\beta)$. Thus, we assume now that $\alpha + \beta + \rho = 0$.

Let $\gamma \notin \{\alpha, \beta, \alpha + \beta\}$ be a root. Again, applying (32) to $(x, y, z) = (A_\alpha, A_\beta, A_\gamma)$ and using the previous steps and the preceding observation, we get $\gamma(\{A_\alpha, A_\beta\}) = 0$ for all such roots $\gamma$. The claim follows.

*Step* 5. For each root $\alpha$, there is a constant $d_\alpha$ such that $\{A_{-\rho}, A_\alpha\} = d_\alpha A_\alpha$. Moreover, if $\alpha, \beta \neq -\rho$ are roots such that $\alpha + \beta$ is a root, then $d_\alpha = -d_\beta$. If $\alpha \neq \pm\rho$, then $d_{-\alpha} = -d_\alpha$.

*Proof.* Since $\{A_{-\rho}, A_\alpha\} \in \mathfrak{g}_\alpha$, the existence of $d_\alpha$ follows. Applying (32) to $(x, y, z) = (A_{-\rho}, A_\alpha, A_\beta)$ yields $(d_\alpha + d_\beta)[A_\alpha, A_\beta] = 0$. The second assertion follows, while the third is obtained by setting $\beta = -\alpha$.

*Step* 6. $d_\alpha = 0$ for all roots $\alpha$; i.e. $\{A_\alpha, A_{-\rho}\} = 0$ for all roots $\alpha$.

*Proof.* Let $\alpha, \beta$ be roots such that $\alpha + \beta$ is a root and $\rho \notin \operatorname{span}(\alpha, \beta)$. Then applying Step 5 several times, we get $d_\alpha = -d_\beta = d_{-(\alpha+\beta)} = -d_{\alpha+\beta} = d_{-\alpha} = -d_\alpha$, and thus, $d_\alpha = 0$ for all $\alpha \neq \pm\rho$. Next, let $\alpha \neq \pm\rho$ be a root such that $\alpha + \rho$ is a root. Then Step 5 implies that $d_\rho = -d_\alpha = 0$. Finally, $d_{-\rho} = 0$ follows from the definition.

*Step* 7. $R = 0$, which yields a contradiction.

*Proof.* Let $\alpha \neq -\rho$ be a root. Applying (32) to $(x, y, z) = (A_0, A_\alpha, A_{-\rho})$ and using the previous steps yields $\alpha(\{A_0, A_{-\rho}\}) = 0$ for all $\alpha \neq -\rho$, and hence $\{A_0, A_{-\rho}\} = 0$. But now it follows that $R = 0$. □

6.2. *Classification of complex Berger algebras.* Let $\mathfrak{g} \subset \mathfrak{gl}(V)$ be an irreducible complex representation, and let $\mathfrak{g}_s$ denote the semisimple part of $\mathfrak{g}$. For each root $\alpha$ of $\mathfrak{g}_s$, we fix $0 \neq A_\alpha \in \mathfrak{g}_\alpha$, and let

$$W_\alpha := \{\text{weights of } A_\alpha V\}.$$

LEMMA 6.2. *Let $\mathfrak{g} \subset \mathfrak{gl}(V)$ be a Berger algebra. Then for every root $\alpha$ of $\mathfrak{g}_s$, there are weights $\lambda_0$,$\lambda_1$ of $V$ such that*

$$W_\alpha \subset \{\lambda_0 + \beta, \lambda_1 + \beta \mid \beta \text{ a root or} = 0\}. \tag{33}$$

*In fact, if $R \in K(\mathfrak{g})$ is a weight element and if there are weight vectors $x_i \in V$ of weights $\lambda_i$ for $i = 0, 1$ such that $R(x_0, x_1) = A_\alpha$, then* (33) *holds for these weights.*

*Proof.* We first show the second assertion. Let $R \in K(\mathfrak{g})$ and $x_i \in V$ as required. Then, for any $y \in V$, the first Bianchi identity of $R \in K(\mathfrak{g})$ reads

$$A_\alpha y = R(y, x_1)x_0 - R(y, x_0)x_1 \in \operatorname{span}\{\mathfrak{g}x_0, \mathfrak{g}x_1\},$$

i.e. $A_\alpha V \subset \operatorname{span}\{\mathfrak{g}x_0, \mathfrak{g}x_1\}$. Then (33) holds since both $A_\alpha V$ and $\operatorname{span}\{\mathfrak{g}x_0, \mathfrak{g}x_1\}$ are a direct sum of weight spaces, and the weights of the latter are contained in the right-hand side of (33).

To show that such an $R$ exists for all roots, let

$$S := \left\{ \alpha \in \Phi \;\middle|\; \begin{array}{l} \text{there are weight elements } R \in K(\mathfrak{g}),\ x_0, x_1 \in V \\ \text{such that } R(x_0, x_1) = A_\alpha \end{array} \right\}.$$

Since $K(\mathfrak{g})$ and $V$ are spanned by their weight vectors, it follows that

$$\underline{\mathfrak{g}} \subset \mathfrak{t}_0 \oplus \bigoplus_{\alpha \in S} \mathfrak{g}_\alpha.$$

Then, since $\mathfrak{g}$ is Berger, it follows that $S = \Phi$. □

Given an irreducible $\mathfrak{g} \subset \mathfrak{gl}(V)$, we call a weight $\lambda$ of $V$ *extremal* if it lies in the orbit of the maximal weight under the Weyl group $W$. We then get the following slight generalization of Lemma 6.2.

LEMMA 6.3. *Let $\mathfrak{g} \subset \mathfrak{gl}(V)$ be a Berger algebra. Then at least one of the following holds.*

1. *There are extremal weights $\lambda_0, \lambda_1$ and a root $\alpha$ for which* (33) *holds.*
2. *There is an affine hyperplane in $\mathfrak{t}^*$ which contains all but two extremal weights.*
3. *There are an extremal weight $\lambda_0$ and a root $\alpha$ such that for all extremal weights $\lambda_1$ with $\langle \lambda_1, \alpha \rangle < 0$, $\lambda_1 - \lambda_0 + \alpha$ is a root or $= 0$.*

*Proof.* Suppose there is a Berger algebra $\mathfrak{g} \subset \mathfrak{gl}(V)$ for which none of these conditions holds. Let $R \in K(\mathfrak{g})$ be a weight element. By Lemma 6.2, the negation of the first condition implies that if $x_i$ is an extremal weight vector, i.e. if $x_i \in V_{\lambda_i}$ with extremal weights $\lambda_i$, then $R(x_0, x_1) \notin \mathfrak{g}_\alpha$ for any root $\alpha$; hence $R(x_0, x_1) \in \mathfrak{t}_0$. Now let $\lambda_2 \neq \lambda_0, \lambda_1$ be any other extremal weight, and let $x_2 \in V_{\lambda_2}$. Then the first Bianchi identity for $(x_0, x_1, x_2)$ implies that $R(x_0, x_1)x_2 = 0$ for all such $x_2$. But if the second condition from above is violated, then this is impossible unless $R(x_0, x_1) = 0$.

Thus, we have that $R(x_0, x_1) = 0$ whenever both $x_i$ are extremal weight vectors. Now let $\lambda$ be an arbitrary weight and let $x \in V_\lambda$. If $R(x_0, x) = A_\alpha$ then the first Bianchi identity for $(x_0, x_1, x)$ would imply that $A_\alpha x_1 \in \mathfrak{g} \cdot x_0$. But if the third condition from above is violated, we can choose $\lambda_1$ such that this is not the case, contradicting $R(x_0, x) = A_\alpha$. On the other hand, if $R(x_0, x) \in \mathfrak{t}_0$ for all $x \in V$, then by the Bianchi identity for $(x_0, x_1, x)$ with an extremal weight vector $x_1$, the violation of the second condition above implies that $R(x_0, x) = 0$.

Thus, it follows that $R(x_0, _) = 0$ for all extremal weight vectors $x_0$. Now, let $x, y \in V$ be arbitrary and apply the Bianchi identity to the triple $(x_0, x, y)$. It follows that $R(x, y)x_0 = 0$ for all extremal weight vectors $x_0$. From here, it follows that $R = 0$, i.e. $K(\mathfrak{g}) = 0$, contradicting that $\mathfrak{g} \subset \mathfrak{gl}(V)$ is Berger. $\square$

For the remainder of this section, we assume that $\mathfrak{g}_s$ is simple. If in addition $\mathfrak{g}$ is not isomorphic to the exceptional Lie algebra $\mathfrak{g}_2$, then we have the following conditions:

(34) If $\alpha, \beta$ are roots then $\alpha + 3\beta$ is not a root.

(35) If $\alpha, \beta$ are roots, then $|\langle \beta, \alpha \rangle| \leqslant 2$; if $\alpha$ is a long root then equality holds if and only if $\alpha = \pm\beta$.

(36) If $\alpha$ is a long root and $\beta$ is any root, then $2\alpha + \beta$ is a root if and only if $\beta = -\alpha$.

Proposition 6.4. *Let $\mathfrak{g} \subset \mathfrak{gl}(V)$ be an irreducible Berger algebra such that $\mathfrak{g}_s$ is simple and not isomorphic to $\mathfrak{g}_2$. If $\alpha$ is a root such that there is a root $\beta$ with $\langle \beta, \alpha \rangle = 1$, then $|\langle \lambda, \alpha \rangle| \leqslant 3$ for all weights $\lambda$ of $V$.*

*Proof. Step* 1. If $\alpha, \beta$ are roots with $\langle \beta, \alpha \rangle = 1$ and if there is a weight $\lambda$ with $\langle \lambda, \alpha \rangle \leqslant -4$, then $\langle \lambda, \beta \rangle \geqslant -2$.

*Proof.* Let $\alpha, \beta, \lambda$ be as required, and let $\lambda_0, \lambda_1$ be as in (33). Then $\lambda + k\alpha \in W_\alpha$ for $k = 1, \ldots, 4$, and thus $\lambda + \alpha = \lambda_0 + \beta_0$, where $\beta_0$ is a root or $= 0$. Also, since by (34) $\beta_0 + 3\alpha$ is *not* a root,

$$(37) \qquad \begin{aligned} \lambda + \alpha &= \lambda_0 + \beta_0 \\ \lambda + 4\alpha &= \lambda_1 + \beta_1 \end{aligned} \qquad \text{where } \beta_i \text{ is a root or} = 0 \text{ for } i = 1, 2.$$

Suppose that $\langle \lambda, \beta \rangle \leqslant -3$. Then $\lambda + 3\beta$ is a weight, and $\langle \lambda + 3\beta, \alpha \rangle \leqslant -1$. Therefore, $W_\alpha \ni \lambda + 3\beta + \alpha = \lambda_0 + \beta_0 + 3\beta = \lambda_1 + \beta_1 + 3(\beta - \alpha)$. Now since $\langle \beta, \alpha \rangle = 1$, $\beta - \alpha$ is a root, hence, by (34), neither $\beta_0 + 3\beta$ nor $\beta_1 + 3(\beta - \alpha)$ are roots, contradicting (33).

*Step* 2. Suppose that $\alpha$ is a root and $\lambda$ a weight with $\langle \lambda, \alpha \rangle \leqslant -4$. Then for any root $\beta$ with $\langle \beta, \alpha \rangle = 1$ we have $\langle \alpha, \beta \rangle = 1$ and $\langle \lambda, \beta \rangle = -2$.

*Proof.* We rescale the inner product such that $(\alpha, \alpha) = 2$. Let $\alpha, \lambda$ be as requested, and let $\beta$ be a root with $\langle \beta, \alpha \rangle = 1$. Then $\alpha - \beta$ is another root with $\langle \alpha - \beta, \alpha \rangle = 1$. Thus, we can apply Step 1 to both $\beta$ and $\alpha - \beta$:

$$\begin{aligned} (\lambda, \beta) &\geqslant -(\beta, \beta), \\ (\lambda, \alpha - \beta) &\geqslant -(\alpha - \beta, \alpha - \beta) = -(\beta, \beta). \end{aligned}$$

Adding these, we have $(\lambda, \alpha) \geqslant -2(\beta, \beta)$. But $\langle \beta, \alpha \rangle = 1$ implies that $(\beta, \beta) \leqslant 2$, and, by hypothesis, $(\lambda, \alpha) \leqslant -4$; hence all inequalities are equalities. In particular, $(\beta, \beta) = (\alpha, \alpha)$; thus $\langle \alpha, \beta \rangle = \langle \beta, \alpha \rangle = 1$. □

*Step* 3. With the hypothesis from the proposition, $\langle \lambda, \alpha \rangle \geqslant -3$.

*Proof.* Let $\alpha, \beta$ be roots with $\langle \beta, \alpha \rangle = 1$, and suppose there is a weight $\lambda$ with $\langle \lambda, \alpha \rangle \leqslant -4$. Choose $\beta_0, \beta_1$ as in (37). By Step 2, we have that $\langle \lambda, \beta \rangle = -2$, hence $\lambda + 2\beta$ is a weight, and $\langle \lambda + 2\beta, \alpha \rangle \leqslant -2$. Therefore, $W_\alpha \ni \lambda + 2\beta + 2\alpha = \lambda_0 + \beta_0 + \alpha + 2\beta$. If $\beta_0 + \alpha + 2\beta$ were a root or $= 0$, then applying (34) twice and using $\langle \alpha, \beta \rangle = 1$ from Step 2, we would get

$$2 \geqslant \langle \beta_0 + \alpha + 2\beta, \beta \rangle = \langle \beta_0, \beta \rangle + 5 \geqslant 3,$$

which is impossible. Thus, since $W_\alpha \ni \lambda + 2\beta + 2\alpha = \lambda_1 + \beta_1 + 2(\beta - \alpha)$, (33) implies that $\beta_1 + 2(\beta - \alpha)$ must be a root. This means that $\beta_1$ is a root with $\langle \alpha - \beta, \beta_1 \rangle = 2$ for every root $\beta$ with $\langle \beta, \alpha \rangle = 1$. But then, replacing $\beta$ by $\alpha - \beta$, we also have $\langle \beta, \beta_1 \rangle = 2$, and thus $\langle \alpha, \beta_1 \rangle = \langle \alpha - \beta, \beta_1 \rangle + \langle \beta, \beta_1 \rangle = 4$ which is impossible.

Application of Step 3 to both $\alpha$ and $-\alpha$ concludes the proof of the proposition.

PROPOSITION 6.5. *Let $\mathfrak{g} \subset \mathfrak{gl}(V)$ be an irreducible Berger algebra such that $\mathfrak{g}_s$ is simple, not isomorphic to $\mathfrak{g}_2$ and* $\mathrm{rank}(\mathfrak{g}_s) > 1$. *Moreover, assume that the root system of long roots of $\mathfrak{g}_s$ is irreducible. Then for every weight $\lambda$ and every long root $\alpha$, $|\langle \lambda, \alpha \rangle| \leqslant 2$.*

Note that the additional hypothesis only excludes those $\mathfrak{g}_s$ whose Dynkin diagram is of type $C_n$; i.e., we assume $\mathfrak{g}_s \not\cong \mathfrak{sp}(n, \mathbb{C})$.

*Proof.* Since for a long root $\alpha$ there exists a root $\beta$ with $\langle \beta, \alpha \rangle = 1$ if $\mathrm{rank}(\mathfrak{g}_s) > 1$, Proposition 6.4 implies that $|\langle \lambda, \alpha \rangle| \leqslant 3$ if $\alpha$ is a long root. Suppose that $|\langle \lambda, \alpha \rangle| = 3$. After replacing $\alpha$ by its negative if necessary, we assume without loss of generality that $\langle \lambda, \alpha \rangle = -3$. Let $\lambda_0, \lambda_1$ be weights satisfying (33).

*Case* 1. $\lambda_i \neq \lambda + 2\alpha$ for $i = 0, 1$. By hypothesis, $\lambda + k\alpha \in W_\alpha$ for $k = 1, 2, 3$. Then $\lambda + \alpha = \lambda_0 + \beta_0$ with $\beta_0 \neq -\alpha$. Thus, by (36), $\beta_0 + 2\alpha$ is not a root, and we have

$$\begin{aligned} \lambda + \alpha &= \lambda_0 + \beta_0 \\ \lambda + 3\alpha &= \lambda_1 + \beta_1 \end{aligned} \quad \text{where } \beta_i \text{ is a root or} = 0 \text{ for } i = 1, 2. \tag{38}$$

Let $\beta$ be a long root with $\langle \beta, \alpha \rangle = 1$, which exists since the root system of long roots is irreducible. After replacing $\beta$ by $\alpha - \beta$ if necessary, we may assume that $\langle \lambda, \beta \rangle \leqslant -2$. Then $\lambda + 2\beta$ is a weight with $\langle \lambda + 2\beta, \alpha \rangle = -1$, and hence $W_\alpha \ni \lambda + 2\beta + \alpha = \lambda_0 + \beta_0 + 2\beta = \lambda_1 + \beta_1 + 2(\beta - \alpha)$. Since $\langle \beta, \alpha \rangle = 1$, $\beta - \alpha$ is a long root. Thus, (33) and (36) imply that either $\beta = -\beta_0$ or $\beta = \alpha - \beta_1$. In other words, the only possible long roots $\beta$ with $\langle \beta, \alpha \rangle = 1$ are $-\beta_0$, $\alpha + \beta_0$, $\alpha - \beta_1$ and $\beta_1$. If $\mathrm{rank}(\mathfrak{g}_s) \geqslant 4$ then there are *more than four* long roots $\beta$ with $\langle \beta, \alpha \rangle = 1$, which yields a contradiction.

Next, suppose that $\mathrm{rank}(\mathfrak{g}_s) = 3$. Then, under the above hypothesis, the root system of long roots of $\mathfrak{g}_s$ is $A_3$. Thus for any long root $\alpha$, there are exactly four long roots $\beta$ with $\langle \beta, \alpha \rangle = 1$, and therefore, $\langle -\beta_0, \alpha \rangle = \langle \beta_1, \alpha \rangle = 1$. Now let $R \in K(\mathfrak{g})$ be a weight element of weight $\rho$ such that $R(x_0, x_1) = A_\alpha$ for weight vectors $x_i \in V_{\lambda_i}$. Then $\rho = \lambda_0 + \lambda_1 - \alpha$, and hence, $\langle \rho, \alpha \rangle = 0$.

Let $x \in V_\lambda$ and $x' \in V_{\lambda + 2\alpha}$. Then $R(x_0, x) \in \mathfrak{g}_{\alpha + \lambda - \lambda_1} = \mathfrak{g}_{2\alpha - \beta_1} = 0$, and $R(x_1, x) \in \mathfrak{g}_{\alpha + \lambda - \lambda_0} = \mathfrak{g}_{-\beta_0}$. Then the Bianchi identity for $(x_0, x_1, x)$ and the fact that $A_\alpha x \neq 0$ imply that $0 \neq R(x_1, x) \in g_{-\beta_0}$. Similarly, the Bianchi identity for $(x_0, x_1, x')$ implies that $0 \neq R(x_0, x') \in \mathfrak{g}_{\beta_1}$.

*Case* 2. $\lambda_0 = \lambda + 2\alpha$. Let $\beta$ as before, and again assume without loss of generality that $\langle \lambda, \beta \rangle \leqslant -2$. Then $W_\alpha \ni \lambda + 2\beta + \alpha = \lambda_0 - \alpha + 2\beta$, and since $-\alpha + 2\beta$ is not a root, it follows from (33) that

$$\lambda + \alpha + 2\beta = \lambda_1 + \beta_1, \quad \beta_1 \text{ a root or } = 0.$$

Let $\gamma$ be a long root with $\langle \gamma, \alpha \rangle = \langle \gamma, \beta \rangle = 1$.

Suppose that $\langle \lambda, \beta - \gamma \rangle \leqslant -2$. Then $\lambda + 2(\beta - \gamma) + k\alpha \in W_\alpha$ for $k = 1, 2, 3$. Now

$$\begin{aligned} \lambda + 2(\beta - \gamma) + \alpha &= \lambda_0 - \alpha + 2(\beta - \gamma) = \lambda_1 + \beta_1 - 2\gamma, \\ \lambda + 2(\beta - \gamma) + 3\alpha &= \lambda_0 + \alpha + 2(\beta - \gamma) = \lambda_1 + \beta_1 + 2(\alpha - \gamma). \end{aligned}$$

But by (36), $\pm\alpha + 2(\beta - \gamma)$ is not a root or $= 0$, and not both $\beta_1 - 2\gamma$ and $\beta_1 + 2(\alpha - \gamma)$ can be a root or $= 0$, violating (33). Thus, $\langle \lambda, \beta - \gamma \rangle \geqslant -1$. Since $\langle \lambda, \beta \rangle \leqslant -2$, this implies, in particular, that $\langle \lambda, \gamma \rangle \leqslant -1$.

Suppose that $\langle \lambda, \gamma \rangle = -1$. Then, by the above, $\langle \lambda, \beta \rangle = -2$, and hence $\lambda + (\beta - \gamma) + k\alpha \in W_\alpha$ for $k = 1, 2, 3$. Now,

$$\begin{aligned} \lambda + (\beta - \gamma) + \alpha &= \lambda_0 - \alpha + (\beta - \gamma) = \lambda_1 + \beta_1 - \beta - \gamma, \\ \lambda + (\beta - \gamma) + 3\alpha &= \lambda_0 + \alpha + (\beta - \gamma) = \lambda_1 + \beta_1 - \beta - \gamma + 2\alpha. \end{aligned}$$

Since $\langle \beta - \gamma, \alpha \rangle = 0$ and these are long roots, $\pm\alpha + (\beta - \gamma)$ is not a root. Therefore, by (33), both $\beta_1 - \beta - \gamma$ and $\beta_1 - \beta - \gamma + 2\alpha$ must be roots, implying, by (36), that $\beta_1 - \beta - \gamma = -\alpha$. But then $1 = \langle \alpha, \gamma \rangle = -\langle \beta_1, \gamma \rangle + 3$, implying by (35) that $\beta_1 = \gamma$; i.e., $\alpha = \beta$, which is impossible.

Thus, $\langle \lambda, \gamma \rangle \leqslant -2$, and therefore, $W_\alpha \ni \lambda + 2\gamma + \alpha = \lambda_0 - \alpha + 2\gamma = \lambda_1 + \beta_1 - 2(\beta - \gamma)$. Now, since $-\alpha + 2\gamma$ is neither a root nor $= 0$, $\beta_1 - 2(\beta - \gamma)$ must be a root or $= 0$ by (33), and (36) implies that $\beta_1 = \beta - \gamma$; i.e., $\lambda_1 = \lambda + \alpha + \beta + \gamma$. But this means that there is at most one long root $\gamma$ with $\langle \gamma, \alpha \rangle = \langle \gamma, \beta \rangle = 1$, namely $\gamma = \beta - \beta_1$. If $\operatorname{rank}(\mathfrak{g}_s) \geqslant 4$ then there is more than one such $\gamma$, which yields a contradiction.

Next, suppose that $\operatorname{rank}(\mathfrak{g}_s) = 3$. Let $R \in K(\mathfrak{g})$ be a weight element of weight $\rho$ such that $R(x_0, x_1) = A_\alpha$ for weight vectors $x_i \in V_{\lambda_i}$. Then again, $\rho = \lambda_0 + \lambda_1 - \alpha$, and hence, $\langle \rho, \alpha \rangle = 0$ follows.

Let $x \in V_{\lambda+\beta}$ and $x' \in V_{\lambda+2\beta}$ such that $A_\alpha x, A_\alpha x'$ are both $\neq 0$. Then an analysis similar to Case 1 implies that $0 \neq R(x_0, x') \in \mathfrak{g}_{\beta-\gamma}$, and either $0 \neq R(x_0, x) \in \mathfrak{g}_{-\gamma}$ or $0 \neq R(x_1, x) \in \mathfrak{g}_{\beta-\alpha}$.

From these two cases, the proposition follows for $\operatorname{rank}(\mathfrak{g}) \geqslant 4$. If $\operatorname{rank}(\mathfrak{g}) = 3$, then in either case, if $R \in K(\mathfrak{g})$ is a weight element of weight $\rho$ such that $R(x_0, x_1) = A_\alpha$ for weight vectors $x_i \in V_{\lambda_i}$ and a long root $\alpha$ then $\langle \rho, \alpha \rangle = 0$. But then, in the first case above, it follows that $\langle \rho, \alpha \rangle = \langle \rho, \beta_0 \rangle =$ $\langle \rho, \beta_1 \rangle = 0$, and hence $\rho = 0$. Treating the second case similarly we conclude the following:

For every weight element $R \in K(\mathfrak{g})$ of weight $\rho \neq 0$ and any $x_i \in V_{\lambda_i}$, either $R(x_0, x_1) \in \mathfrak{t}_0$ or $R(x_0, x_1)$ is a root vector of a short root.

But now a straightforward analysis shows that this is impossible; i.e., the only weight of $K(\mathfrak{g})$ is $\rho = 0$; i.e., $K(\mathfrak{g})$ is a trivial $\mathfrak{g}$-module, and we get a contradiction from Lemma 2.3.

In the case $\mathrm{rank}(\mathfrak{g}) = 2$, i.e., if the root system is of type $A_2$, the statement is shown by a direct and straightforward analysis. □

PROPOSITION 6.6. *Let $\mathfrak{g} \subset \mathfrak{gl}(V)$ be as in Proposition* 6.4, *such that $\mathfrak{g}_s \not\cong \mathfrak{so}(7, \mathbb{C})$ and* $\mathrm{rank}(\mathfrak{g}_s) > 2$. *Suppose there is a weight $\lambda$ and a (not necessarily long) root $\alpha$ with $|\langle \lambda, \alpha \rangle| \geqslant 2$. Then one of the following must hold.*

1. $|\langle \lambda, \beta \rangle| \leqslant 1$ *for all long roots $\beta$ with $\langle \alpha, \beta \rangle = 0$, or*
2. *There is a unique long root $\gamma$ with $\langle \alpha, \gamma \rangle = 0$ and $\langle \lambda, \gamma \rangle \geqslant 2$, and the root system*
   $\{$*long roots* $\beta \mid \langle \alpha, \beta \rangle = 0\}$ *contains* $\{\pm\gamma\}$ *as a direct summand.*

We shall split the proof of this proposition into two lemmas.

LEMMA 6.7. *Let $\mathfrak{g}$, $V$, $\lambda$ and $\alpha$ be as in Proposition* 6.6. *Then the set*

$$\{\textit{long roots } \gamma \mid \langle \alpha, \gamma \rangle = 0,\ \langle \lambda, \gamma \rangle \geqslant 2\}$$

*contains at most one element.*

*Proof.* After replacing $\alpha$ by $-\alpha$ if necessary, we assume without loss of generality that $\langle \lambda, \alpha \rangle \leqslant -2$. Let $\lambda_0, \lambda_1$ be as in (33). Suppose there are two distinct long roots $\gamma_1, \gamma_2$ with $\langle \alpha, \gamma_i \rangle = 0$ and $\langle \lambda, \gamma_i \rangle \geqslant 2$ for $i = 1, 2$. If $\gamma_1 + \gamma_2$ were a root then $\langle \lambda, \gamma_1 + \gamma_2 \rangle \geqslant 4$, contradicting Proposition 6.4. It follows that $\langle \gamma_1, \gamma_2 \rangle \in \{0, 1\}$.

*Case* 1. $\langle \gamma_1, \gamma_2 \rangle = 0$. In this case, we have that $\lambda + k\alpha - r\gamma_1 - s\gamma_2 \in W_\alpha$ for $k = 1, 2$, $r, s = 0, 1, 2$. Thus, by (33),

$$\lambda + \alpha = \lambda_0 + \beta_0, \quad \text{where } \beta_0 \text{ is a root or } = 0.$$

After possibly switching $\gamma_1$ and $\gamma_2$, we assume without loss of generality that $\beta_0 \neq \gamma_2$, i.e. $\beta_0 - 2\gamma_2$ is not a root. Therefore, by (33),

$$\lambda + \alpha - 2\gamma_2 = \lambda_1 + \beta_1, \quad \text{where } \beta_1 \text{ is a root or } = 0.$$

Next, we consider $W_\alpha \ni \lambda + \alpha - 2\gamma_1 = \lambda_0 + \beta_0 - 2\gamma_1 = \lambda_1 + \beta_1 - 2\gamma_1 + 2\gamma_2$. Since $\langle \beta_1 - 2\gamma_1 + 2\gamma_2, \gamma_2 \rangle = \langle \beta_1, \gamma_2 \rangle + 4 \geqslant 2$ with equality if and only if $\beta_1 = -\gamma_2$ by (35), it follows that $\beta_1 - 2\gamma_1 + 2\gamma_2$ is not a root. Thus, by (33) and (36), $\beta_0 = \gamma_1$, i.e. $\lambda_0 = \lambda + \alpha - \gamma_1$.

Now consider $W_\alpha \ni \lambda + 2\alpha - 2\gamma_1 = \lambda_0 + \alpha - \gamma_1 = \lambda_1 + \beta_1 + \alpha - 2\gamma_1 + 2\gamma_2$. Since $\langle \gamma_1, \alpha \rangle = 0$, $\alpha - \gamma_1$ is not a root, hence $\beta_1 + \alpha - 2\gamma_1 + 2\gamma_2$ must be a root or $= 0$ by (33), and this yields a contradiction in a similar way.

*Case* 2. $\langle \gamma_1, \gamma_2 \rangle = 1$. In this case, we have that $\lambda + k\alpha - r\gamma_1 - s\gamma_2 \in W_{-\alpha}$ for $k = 1, 2$, $r, s = 0, 1, 2$ with $r + s \leqslant 2$. Thus, by (33),

$$\lambda + \alpha = \lambda_0 + \beta_0, \quad \text{where } \beta_0 \text{ is a root or } = 0.$$

Again, without loss of generality we may assume $\beta_0 \neq \gamma_2$, and therefore,

$$\lambda + \alpha - 2\gamma_2 = \lambda_1 + \beta_1, \quad \text{where } \beta_1 \text{ is a root or } = 0.$$

Next, we consider $W_\alpha \ni \lambda + \alpha - 2\gamma_1 = \lambda_0 + \beta_0 - 2\gamma_1 = \lambda_1 + \beta_1 + 2(\gamma_2 - \gamma_1)$. By (33) and (36), either $\beta_0 = \gamma_1$ or $\beta_1 = \gamma_1 - \gamma_2$. If $\beta_0 = \gamma_1$, i.e. $\lambda_0 = \lambda + \alpha - \gamma_1$, then $W_\alpha \ni \lambda + 2\alpha = \lambda_0 + \alpha + \gamma_1 = \lambda_1 + \beta_1 + \alpha + 2\gamma_2$. But $\alpha + \gamma_1$ is not a root as $\langle \gamma_1, \alpha \rangle = 0$, and $\beta_1 + \alpha + 2\gamma_2$ is seen not to be a root or $= 0$ in a similar way as before. Hence, we get a contradiction from (33). The case $\beta_1 = \gamma_1 - \gamma_2$ is ruled out in a similar way, and this proves the lemma. □

Lemma 6.8. *Let $\mathfrak{g}$, $V$, $\lambda$ and $\alpha$ as before, and suppose the set in Lemma* 6.7 *is nonempty. Then there is no long root $\delta$ other than $\pm\gamma$ with $\langle \alpha, \delta \rangle = 0$ and $\langle \gamma, \delta \rangle \neq 0$, unless $\mathfrak{g}_s \cong \mathfrak{so}(7, \mathbb{C})$, and $\alpha$ is a short root.*

*Proof.* The claim is obvious if the root system of long roots is reducible (i.e. $\mathfrak{g}_s \cong \mathfrak{sp}(n, \mathbb{C})$) or if $\operatorname{rank}(\mathfrak{g}) \leqslant 3$. Thus, without loss of generality we may assume that the premises of Proposition 6.5 are satisfied, and therefore we may assume that $\langle \lambda, \alpha \rangle = -\langle \lambda, \gamma \rangle = -2$. As usual, let $\lambda_0, \lambda_1$ be the weights satisfying (33).

Suppose there is a long root $\delta \neq \pm\gamma$ with $\langle \alpha, \delta \rangle = 0$ and $\langle \gamma, \delta \rangle \neq 0$. After replacing $\delta$ by $-\delta$ if necessary, we may assume that $\langle \gamma, \delta \rangle = 1$. Then $\gamma - \delta$ is another such root, and because of Lemma 6.7, we have both $\langle \lambda, \delta \rangle \leqslant 1$ and $\langle \lambda, \gamma - \delta \rangle \leqslant 1$, i.e. $\langle \lambda, \delta \rangle = 1$. It follows that

$$\{\lambda + k\alpha - r\gamma,\ \lambda + k\alpha - \delta \mid k = 1, 2,\ r = 0, 1, 2\} \subset W_\alpha.$$

*Case* 1. $\lambda_i \neq \lambda + \alpha - \gamma$ for $i = 0, 1$. As in previous cases, (33) yields

$$\begin{aligned} \lambda + \alpha &= \lambda_0 + \beta_0, \\ \lambda + \alpha - 2\gamma &= \lambda_1 + \beta_1, \end{aligned} \quad \text{where } \beta_0 \neq \gamma \text{ and } \beta_1 \neq -\gamma \text{ are roots or } = 0.$$

Now, $W_\alpha \ni \lambda + 2\alpha - 2\gamma = \lambda_0 + \beta_0 + \alpha - 2\gamma = \lambda_1 + \beta_1 + \alpha$. But $\langle \beta_0 + \alpha - 2\gamma, \gamma \rangle = \langle \beta_0, \gamma \rangle - 4 \leqslant -3$, since $\beta_0 \neq \gamma$, and thus $\beta_0 + \alpha - 2\gamma$ is not a root. Therefore, $\beta_1 + \alpha$ must be a root or $= 0$.

Next, for $k = 1, 2$, $W_\alpha \ni \lambda + k\alpha - \delta = \lambda_0 + \beta_0 + (k-1)\alpha - \delta = \lambda_1 + \beta_1 + (k-1)\alpha + 2\gamma - \delta$. But $\langle \beta_1 + (k-1)\alpha + 2\gamma - \delta, \gamma \rangle = \langle \beta_1, \gamma \rangle + 3 \geqslant 2$, since $\beta_1 \neq -\gamma$. Hence, if $\beta_1 + (k-1)\alpha + 2\gamma - \delta$ was a root it would have to equal $\gamma$; i.e., $\beta_1 = -(k-1)\alpha + (\delta - \gamma)$. If $k = 2$, this is impossible since $\langle \delta - \gamma, \alpha \rangle = 0$. Also, if $k = 1$ this would contradict $\beta_1 + \alpha$ being a root or $= 0$. Therefore, both $\beta_0 - \delta$ and $\beta_0 + \alpha - \delta$ must be roots or $= 0$.

Finally, $W_\alpha \ni \lambda + 2\alpha = \lambda_0 + \beta_0 + \alpha = \lambda_1 + \beta_1 + \alpha + 2\gamma$. We exclude that $\beta_1 + \alpha + 2\gamma$ is a root or $= 0$ as before; hence $\beta_0 + \alpha$ is a root or $= 0$.

Since $\beta_0, \beta_0 + \alpha, \beta_0 - \delta, \beta_0 + \alpha - \delta$ are all roots or $= 0$, we conclude that $\beta_0 \notin \{0, -\alpha, \delta\}$, and from there it follows that $\langle \beta_0, \delta \rangle = 1$ for all long roots $\delta$ with $\langle \alpha, \delta \rangle = 0$, $\langle \gamma, \delta \rangle = 1$. However, if such a $\delta$ exists, then $\gamma - \delta$ is another such root; hence $1 = \langle \beta_0, \gamma - \delta \rangle = \langle \beta_0, \gamma \rangle - \langle \beta_0, \delta \rangle$, and hence $\langle \beta_0, \gamma \rangle = 2$, contradicting that $\beta_0 \neq \gamma$.

*Case* 2. $\lambda_0 = \lambda + \alpha - \gamma$. In this case, since $\gamma + \alpha$ is not a root, it follows that

$$\lambda + 2\alpha = \lambda_1 + \beta_1, \qquad \text{where } \beta_1 \text{ is a root or } = 0.$$

Moreover, $\lambda + 2\alpha - 2\gamma = \lambda_0 + \alpha - \gamma = \lambda_1 + \beta_1 - 2\gamma$. Since $\alpha - \gamma$ is not a root, (36) implies that $\beta_1 = \gamma$, i.e.

$$\lambda_1 = \lambda_0 + \alpha.$$

Thus, using Lemma 6.2 and the fact that the Weyl group acts transitively on the set of long roots, we conclude:

(39) If $x, y \in V$ are weight vectors of weight $\lambda, \mu$ respectively, and if $R \in K(\mathfrak{g})$ is a weight element such that $R(x, y) = A_\beta$ for some long root $\beta$, then $\lambda - \mu$ must be a long root.

Now, by Lemma 6.2, there is a weight element $R \in K(\mathfrak{g})$ such that $R(x_0, x_1) = A_\alpha$ where $x_i \in V$ are weight vectors of weights $\lambda_i$ for $i = 0, 1$. Let $x_2$ be a weight vector of weight $\lambda - \delta$. Then $R(x_0, x_2) \in \mathfrak{g}_{\alpha+\gamma-\delta} = 0$, since $\alpha + \gamma - \delta$ is not a root, and $R(x_1, x_2) \in \mathfrak{g}_{\gamma-\delta}$; however, since $\gamma - \delta$ is a long root, and $\lambda_1 - (\lambda - \delta) = 2\alpha - \gamma + \delta$ is not a root, (39) implies that $R(x_1, x_2) = 0$. But now the Bianchi identity for $x_0, x_1, x_2$ reads

$$A_\alpha x_2 = 0 \qquad \text{for all } x_2 \in V_{\lambda-\delta},$$

contradicting that $\lambda - \delta + \alpha$ is a weight. □

It is now clear that Proposition 6.6 follows immediately from Lemmas 6.7 and 6.8.

PROPOSITION 6.9. *Let $\mathfrak{g} \subset \mathfrak{gl}(V)$ be an irreducible complex Lie subalgebra whose semisimple part is simple and has rank at least three. Then $\mathfrak{g}$ is a*

*Berger algebra only if either its sky is biholomorphic to a compact hermitean-symmetric manifold, or if* $\mathfrak{g}$ *is* (*the central extension of*) *the image of one of the following representations*:

(i) 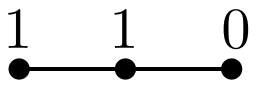

(ii) 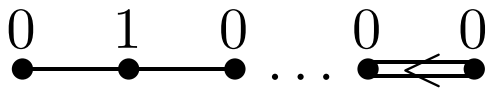

(iii) 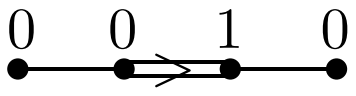

(iv) 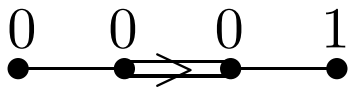

(v) 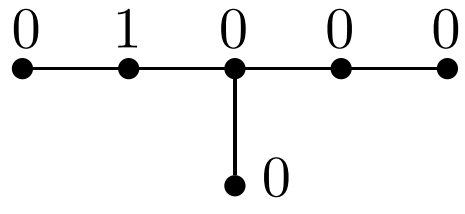

(vi) 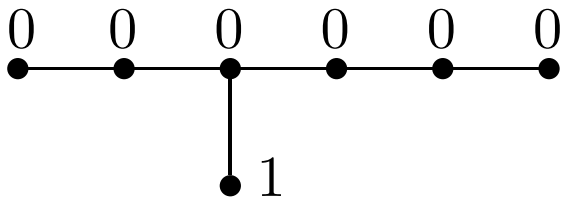

.

*Proof.* The list of irreducible hermitean symmetric spaces is given at the beginning of Section 4. We give the proof for each type of the Dynkin diagram of $G$.

1. *Type* $A_n$. In this case, the root system is $\Phi = \{\alpha_{i,j} := \theta_i - \theta_j \mid i \neq j \in \{1, \dots, n+1\}\ \}$, and the positive roots are $\Phi^+ = \{\alpha_{i,j} \mid i < j\}$. The dominant weight of $G$ can be represented in a unique way as $\lambda_0 = c_1\theta_1 + \ldots + c_n\theta_n$ with integers $c_1 \geqslant \ldots c_n \geqslant 0$. For convenience, we set $c_{n+1} = 0$. Note that due to the symmetry of the root system $A_n$ we may assume without loss of generality that $c_1 - c_2 \geqslant c_n$. From Proposition 6.5 it follows that $\langle \lambda_0, \alpha_{1,n+1} \rangle = c_1 \leqslant 2$, i.e. $c_1 = 1$ or 2. If $c_1 = 1$ then the representation of $G$ on $V$ is a fundamental representation, and its sky is hermitean symmetric. Thus, we may assume that $c_1 = 2$. Since the system of roots perpendicular to $\alpha_{1,n+1}$ is $A_{n-2}$ which is irreducible and does not contain $A_1$ as a summand for $n \geqslant 4$, it follows from Proposition 6.6 that $\langle \lambda_0, \alpha_{2,n} \rangle = c_2 - c_n \leqslant 1$ if $n \geqslant 4$. If $n = 3$ and $c_2 - c_n = 2$ then $\lambda_0 = 2(\theta_1 + \theta_2)$, and then the sky of $G$ is hermitean symmetric. If $c_2 - c_n = 0$ then, since $c_1 - c_2 \geqslant c_n$, either $\lambda_0 = 2\theta_1$ or $\lambda_0 = 2\theta_1 + \theta_2 + \ldots + \theta_n$. In the first case, the sky of $G$ is hermitean symmetric, while in the second case $G$ acts via the adjoint representation and hence is not Berger by Proposition 6.1. Thus, we may assume that $c_2 - c_n = 1$ and, because of $c_1 - c_2 \geqslant c_n$, we have $c_2 = 1$, $c_n = 0$. In other words,

$$\lambda_0 = 2\theta_1 + \theta_2 + \ldots + \theta_k, \quad \text{for some } k \in \{2, \dots, n-1\}.$$

The Weyl group of $A_n$ is the permutation group $S_{n+1}$ which acts by permutation of the indices of $\theta_1, \dots, \theta_{n+1}$. Thus, the extremal weights are all elements of the form $\lambda = 2\theta_{m_0} + \theta_{m_1} + \theta_{m_k}$ for distinct indices $m_0, \dots, m_k$. From here, it is now straightforward to show that the last two conditions of Lemma 6.3 are violated, while the first one holds only if $n = 3$ and $\lambda_0 = 2\theta_1 + \theta_2$. This corresponds precisely to (i) above.

2. *Type* $B_n$. The root system is $\Phi = \{\pm\theta_i,\ \pm\theta_i \pm \theta_j \mid i < j, i = 1, \dots, n\}$, and the positive roots are $\Phi^+ = \{\theta_i,\ \theta_i \pm \theta_j \mid i < j, i = 1, \dots, n\}$. The maximal weight is given by $\lambda_0 = c_1\theta_1 + \ldots + c_n\theta_n$ with $c_1 \geqslant \ldots \geqslant c_n \geqslant 0$, where either all $c_k$ are integers, or all $c_k$ are half-integers. From Proposition 6.5 it follows that $\langle \lambda_0, \theta_1 + \theta_2 \rangle = c_1 + c_2 \leqslant 2$. If $c_2 = 0$, then obviously the sky of $G$ is hermitean symmetric. Thus, we may assume that $c_2 > 0$, and hence $c_1 \leqslant \frac{3}{2}$.

Suppose that $c_1 = \frac{3}{2}$. Then $\lambda_0 = \theta_1 + \frac{1}{2}(\theta_1 + \ldots + \theta_n)$. In this case, the extremal weights are of the form $\pm\theta_k + \frac{1}{2}(\pm\theta_1 \pm \ldots \pm \theta_n)$. From here, one sees easily that the last two conditions of Lemma 6.3 are violated, while the first one holds only if $n \leqslant 3$. However, a more careful analysis shows that the case $n = 3$ can be ruled out for Berger groups as well.

Next, suppose that $c_1 = 1$. If $c_2 = 0$ then the corresponding sky is hermitean symmetric; thus we assume that $c_2 = 1$, i.e.,

$$\lambda_0 = \theta_1 + \ldots + \theta_k.$$

If $k = 2$ then this is the adjoint representation, which is not Berger; if $k = n$ then the corresponding sky is hermitean symmetric. Thus, we assume that $3 \leqslant k \leqslant n - 1$ and, in particular, that $n \geqslant 4$. It is now straightforward to show that none of the conditions in Lemma 6.3 hold; thus these groups are not Berger.

Finally, suppose that $c_1 = \frac{1}{2}$. Then the corresponding representation is the spinor representation whose sky is hermitean symmetric.

3. *Type* $C_n$. The root system is $\Phi = \{\pm 2\theta_i,\ \pm\theta_i \pm \theta_j \mid i < j, i = 1, \ldots, n\}$, and the positive roots are $\Phi^+ = \{2\theta_i,\ \theta_i \pm \theta_j \mid i < j, i = 1, \ldots, n\}$. The maximal weight is given by $\lambda_0 = c_1\theta_1 + \ldots + c_n\theta_n$ with integers $c_1 \geqslant \ldots \geqslant c_n \geqslant 0$. From Proposition 6.4 it follows that $\langle \lambda_0, \theta_1 + \theta_2 \rangle = c_1 + c_2 \leqslant 3$. If $c_2 = 0$, then obviously the sky of $G$ is hermitean symmetric. Thus, we may assume that $c_2 > 0$, and hence $c_1 \leqslant 2$.

Suppose that $c_1 = 2$. Then $\lambda_0 = 2\theta_1 + \theta_2 + \ldots + \theta_k$, with $2 \leqslant k \leqslant n$. The extremal weights are then all weights of the form $\pm 2\theta_{m_1} \pm \theta_{m_2} \pm \ldots \pm \theta_{m_k}$ for some pairwise distinct indices $m_1, \ldots, m_k$. It is now easy to see that the conditions in Lemma 6.3 are violated if $n \geqslant 3$; hence these subgroups cannot be Berger.

Next, suppose that $c_1 = 1$; hence $\lambda_0 = \theta_1 + \ldots + \theta_k$ for some $k \in \{1, \ldots, n\}$. The sky of this representation is hermitean symmetric if $k = 1$ or $k = n$; hence we assume that $2 \leqslant k \leqslant n - 1$. Once again, the extremal weights are of the form $\pm\theta_{m_1} \pm \ldots \pm \theta_{m_k}$ for some pairwise distinct indices $m_1, \ldots, m_k$, and then straightforward calculation shows that the conditions in Lemma 6.3 are violated if $k \geqslant 3$.

Thus, we assume that $k = 2$, and this is precisely (ii) from above.

4. *Type* $D_n$. The root system is $\Phi = \{\pm\theta_i \pm \theta_j \mid i < j, i = 1, \ldots, n\}$, and is hence contained in the root system $B_n$. Thus, the investigation of this case is very similar to that of the case $B_n$. The details are omitted.

5. *Type $F_4$.* Suppose that the maximal weight $\lambda$ is given by $\overset{c_1}{\bullet}\!-\!\overset{c_2}{\bullet}\!\Rightarrow\!\overset{c_3}{\bullet}\!-\!\overset{c_4}{\bullet}$. The maximal short root of $F_4$ is given by $\alpha = \overset{1}{\bullet}\!-\!\overset{2}{\bullet}\!\Rightarrow\!\overset{3}{\bullet}\!-\!\overset{2}{\bullet}$, and for $\beta = \overset{0}{\bullet}\!-\!\overset{0}{\bullet}\!\Rightarrow\!\overset{0}{\bullet}\!-\!\overset{1}{\bullet}$ we have $\langle \beta, \alpha \rangle = 1$. Thus, Proposition 6.4 implies that $\langle \lambda, \alpha \rangle = 2c_1 + 4c_2 + 3c_3 + 2c_4 \leqslant 3$. Thus, $c_2 = 0$ and exactly one of $c_1, c_3, c_4$ is equal to 1. However, the case $c_1 = 1$ corresponds to the adjoint representation which is not Berger by Proposition 6.1, and the two remaining cases are precisely (iii) and (iv).

6. *Type $E_6$.* Let $\lambda$ be given by $\begin{array}{c}\overset{c_1}{\bullet}\!-\!\overset{c_2}{\bullet}\!-\!\overset{c_3}{\bullet}\!-\!\overset{c_4}{\bullet}\!-\!\overset{c_5}{\bullet}\\ \quad\ \ |\\ \quad\ \ \bullet\, c_6\end{array}$. The maximal root of $E_6$ is given by $\begin{array}{c}\overset{1}{\bullet}\!-\!\overset{2}{\bullet}\!-\!\overset{3}{\bullet}\!-\!\overset{2}{\bullet}\!-\!\overset{1}{\bullet}\\ \quad\ \ |\\ \quad\ \ \bullet\, 2\end{array}$, and the root system orthogonal to the maximal root is $A_5$ with maximal root $\begin{array}{c}\overset{1}{\bullet}\!-\!\overset{1}{\bullet}\!-\!\overset{1}{\bullet}\!-\!\overset{1}{\bullet}\!-\!\overset{1}{\bullet}\\ \quad\ \ |\\ \quad\ \ \bullet\, 0\end{array}$. It follows from Propositions 6.5 and 6.6 that $c_1 + 2c_2 + 3c_3 + 2c_4 + c_5 + 2c_6 \leqslant 2$ and $c_1 + \ldots + c_5 \leqslant 1$. If $c_6 = 1$ then all other $c_i = 0$, and this describes the adjoint representation which is not Berger by Proposition 6.1. If $c_1 = 1$ or $c_5 = 1$, then all other $c_i = 0$, hence the sky of $G$ is hermitean symmetric. Thus, the only possibility left is (v).

7. *Type $E_7$ or $E_8$.* These cases are dealt with in complete analogy to $E_6$. □

6.3. *Skies of complex Berger algebras.*

THEOREM 6.10. *Let $\mathfrak{g} \subset \mathfrak{gl}(V)$ be an irreducible complex Lie subalgebra whose semisimple part is simple. Then $\mathfrak{g}$ is a Berger algebra only if its sky is biholomorphic to a compact hermitean-symmetric manifold.*

*Proof.* By Proposition 6.9, it suffices to show that the representations given there are not Berger algebras, and to consider the case where $\mathrm{rank}(\mathfrak{g}_s) = 2$.

*Type $A_3$* (i). We have the associated Borel-Weil data $(X, L) = \overset{1}{\times}\!-\!\overset{1}{\times}\!-\!\overset{0}{\bullet}$ There is a natural projection $X \longrightarrow X_0 = \times\!-\!\bullet\!-\!\bullet$ whose fibre $X_v$ is $\mathbb{P}_2$. The short exact sequence (6) takes the form

$$TX = \overset{-1}{\times}\!-\!\overset{1}{\times}\!-\!\overset{1}{\bullet} \ + \ \nu^*(\overset{1}{\times}\!-\!\overset{0}{\bullet}\!-\!\overset{1}{\bullet}).$$

By Lemma 2.7(i),

$$\mathrm{H}^1(X, \odot^2 TX \otimes L^{*2}) \subset \mathrm{H}^1(\overset{-2}{\times}\!-\!\overset{0}{\times}\!-\!\overset{2}{\bullet}) = 0$$

and

$$\mathrm{H}^1(X, \odot^3 TX \otimes L^{*2}) \subset \mathrm{H}^1(\overset{-5}{\times}\!-\!\overset{1}{\times}\!-\!\overset{3}{\bullet}) + \mathrm{H}^1(\overset{-4}{\times}\!-\!\overset{0}{\bullet}\!-\!\overset{2}{\bullet} \otimes \overset{-1}{\times}\!-\!\overset{0}{\bullet}\!-\!\overset{1}{\bullet}) = 0.$$

Then, by the estimation (19), $\mathrm{H}^1(X, L \otimes \odot^3 N^*) = 0$ implying $K(\mathfrak{g} \oplus \mathbb{C}) = 0$.

*Type $A_2$.* Let $(X, L) = \overset{m}{\times}\!-\!\overset{n}{\times}$ with $m \geqslant 2$ and $n \geqslant 1$. There is a fibration $\nu : X \to X_0 = \times\!-\!\bullet$. The short exact sequence (6) takes the form

$$TX = \overset{-1}{\times}\!-\!\overset{2}{\times} \; + \; \nu^*(\overset{1}{\times}\!-\!\overset{1}{\bullet}),$$

implying

$$\odot^2 TX \otimes L^{*2} = \overset{-2-2m}{\times}\!-\!\overset{4-2n}{\times} + \overset{-1-2m}{\times}\!-\!\overset{2-2n}{\times} \otimes \nu^*(\overset{1}{\times}\!-\!\overset{1}{\bullet}) \; + \overset{-2m}{\times}\!-\!\overset{-2n}{\times} \otimes \nu^*(\overset{2}{\times}\!-\!\overset{2}{\bullet})$$

and

$$\begin{aligned} \odot^3 TX \otimes L^{*2} \;\; &= \overset{-3-2m}{\times}\!-\!\overset{6-2n}{\times} + \overset{-2-2m}{\times}\!-\!\overset{4-2n}{\times} \otimes \nu^*(\overset{1}{\times}\!-\!\overset{1}{\bullet}) \\ &\quad + \overset{-1-2m}{\times}\!-\!\overset{2-2n}{\times} \otimes \nu^*(\overset{2}{\times}\!-\!\overset{2}{\bullet}) \\ &\quad + \overset{-2m}{\times}\!-\!\overset{-2n}{\times} \otimes \nu^*(\overset{3}{\times}\!-\!\overset{3}{\bullet}). \end{aligned}$$

Applying the direct image functor $\nu^i_*$ and using the formula

$$\overset{a}{\times}\!-\!\overset{b}{\bullet} \otimes \overset{c}{\times}\!-\!\overset{d}{\bullet} = \bigoplus_{0 \leqslant i \leqslant b} \overset{a+c+i-1}{\times}\!-\!\overset{c+d-2i}{\bullet} \quad , \qquad b \leqslant d,$$

one easily finds that $\mathrm{H}^1(X, \odot^2 TX \otimes L^{*2}) = 0$ for all $m \geqslant 2, n \geqslant 1$ and

$$\mathrm{H}^1(X, \odot^3 TX \otimes L^{*2}) \subset \begin{cases} \mathrm{H}^1(\overset{-3}{\times}\!-\!\overset{2}{\bullet}) & \text{for } m = 2, n = 1 \\ 0 & \text{otherwise.} \end{cases}$$

Hence, by (19), we get an estimation

$$K(\mathfrak{g} \oplus \mathbb{C}) \subset \begin{cases} \overset{1}{\bullet}\!-\!\overset{0}{\bullet} & \text{for } m = 2, n = 1 \\ 0 & \text{otherwise.} \end{cases}$$

This leaves the cases where $m = 2, n = 1$ and $m = n = 1$. However, a direct analysis using arguments similar to those from subsection 7.1. shows that these representations do not yield Berger algebras. We omit further details.

*Type $B_2$.* The Borel-Weil data associated with a representation $\overset{s}{\bullet}\!\Longrightarrow\!\overset{t}{\bullet}$, $s, t \geqslant 1$, are $(X, L) = \overset{s}{\times}\!\Longrightarrow\!\overset{t}{\times}$. There is a projection $\nu : X \to X_0 = \times\!\Longrightarrow\!\bullet$ with the associated short exact sequence (6) given by

$$TX = \overset{-1}{\times}\!\Longrightarrow\!\overset{2}{\times} \; + \; \nu^*(\overset{0}{\times}\!\Longrightarrow\!\overset{2}{\bullet}).$$

Then

$$\odot^2 TX \otimes L^{*2} = \overset{-2-2s}{\times}\!\Longrightarrow\!\overset{4-2t}{\times} + \overset{-1-2s}{\times}\!\Longrightarrow\!\overset{2-2t}{\times} \otimes \nu^*(\overset{0}{\times}\!\Longrightarrow\!\overset{2}{\bullet}) + \overset{-2s}{\times}\!\Longrightarrow\!\overset{-2t}{\times} \otimes \nu^*(\odot^2 \overset{0}{\times}\!\Longrightarrow\!\overset{2}{\bullet}),$$

implying $\mathrm{H}^1(\odot^2 TX \otimes L^{*2}) = 0$ for all $s, t \geqslant 1$. A similar analysis of $\odot^3 TX \otimes L^{*2}$ produces an estimation

$$\mathrm{H}^1(X, \odot^3 TX \otimes L^{*2}) \subset \begin{cases} \mathrm{H}^1(\overset{-3}{\times}\!\Longrightarrow\!\overset{4}{\bullet}) & \text{for } m = n = 1, \\ 0 & \text{otherwise.} \end{cases}$$

Hence, by (19), we get an estimation

$$K(\mathfrak{g} \oplus \mathbb{C}) \subset \begin{cases} \overset{1}{\bullet}\!\Longrightarrow\!\overset{0}{\bullet} & \text{for } m = n = 1 \\ 0 & \text{otherwise.} \end{cases}$$

Finally, the case $m = n = 1$ can be analyzed by methods similar to those of of subsection 7.1. It turns out that $\mathfrak{g} \oplus \mathbb{C}$ is not a Berger algebra in this case either. We omit the details.

*Type $C_n$* (ii), $n \geqslant 3$. The Borel-Weil data associated with the fundamental representation of the second node are $(X, L) = \overset{0}{\bullet}\!-\!\overset{1}{\times}\!-\!\overset{0}{\bullet}\!-\!\overset{0}{\bullet}\cdots\overset{0}{\bullet}\!\Longleftarrow\!\overset{0}{\bullet}$. By (5),

$$TX = \overset{1}{\bullet}\!-\!\overset{-1}{\times}\!-\!\overset{1}{\bullet}\!-\!\overset{0}{\bullet}\cdots\overset{0}{\bullet}\!\Longleftarrow\!\overset{0}{\bullet} + \overset{2}{\bullet}\!-\!\overset{0}{\times}\!-\!\overset{0}{\bullet}\!-\!\overset{0}{\bullet}\cdots\overset{0}{\bullet}\!\Longleftarrow\!\overset{0}{\bullet},$$

implying

$$\odot^2 TX \otimes L^{*2} = \begin{array}{c} \overset{2}{\bullet}\!-\!\overset{-2}{\times}\!-\!\overset{2}{\bullet}\!-\!\overset{0}{\bullet}\cdots\overset{0}{\bullet}\!\Longleftarrow\!\overset{0}{\bullet} \\ \oplus \\ \overset{0}{\bullet}\!-\!\overset{-1}{\times}\!-\!\overset{0}{\bullet}\!-\!\overset{0}{\bullet}\cdots\overset{0}{\bullet}\!\Longleftarrow\!\overset{0}{\bullet} \\ \oplus \\ \overset{0}{\bullet}\!-\!\overset{-2}{\times}\!-\!\overset{0}{\bullet}\!-\!\overset{1}{\bullet}\cdots\overset{0}{\bullet}\!\Longleftarrow\!\overset{0}{\bullet} \end{array} + \begin{array}{c} \overset{3}{\bullet}\!-\!\overset{-3}{\times}\!-\!\overset{1}{\bullet}\!-\!\overset{0}{\bullet}\cdots\overset{0}{\bullet}\!\Longleftarrow\!\overset{0}{\bullet} \\ \oplus \\ \overset{1}{\bullet}\!-\!\overset{-2}{\times}\!-\!\overset{1}{\bullet}\!-\!\overset{0}{\bullet}\cdots\overset{0}{\bullet}\!\Longleftarrow\!\overset{0}{\bullet} \end{array} + \begin{array}{c} \overset{4}{\bullet}\!-\!\overset{-2}{\times}\!-\!\overset{0}{\bullet}\!-\!\overset{0}{\bullet}\cdots\overset{0}{\bullet}\!\Longleftarrow\!\overset{0}{\bullet} \\ \oplus \\ \overset{0}{\bullet}\!-\!\overset{0}{\times}\!-\!\overset{0}{\bullet}\!-\!\overset{0}{\bullet}\cdots\overset{0}{\bullet}\!\Longleftarrow\!\overset{0}{\bullet} \end{array}$$

where the summand $\overset{0}{\bullet}\!-\!\overset{-2}{\times}\!-\!\overset{0}{\bullet}\!-\!\overset{1}{\bullet}\cdots\overset{0}{\bullet}\!\Longleftarrow\!\overset{0}{\bullet}$ must be set to zero in the case $n = 3$. All the irreducible bundles on the right-hand side have the first cohomology vanishing except $\mathrm{H}^1(\overset{1}{\bullet}\!-\!\overset{-2}{\times}\!-\!\overset{1}{\bullet}\!-\!\overset{0}{\bullet}\cdots\overset{0}{\bullet}\!\Longleftarrow\!\overset{0}{\bullet}) = \mathbb{C}$ which implies an estimation $\mathrm{H}^1(X, \odot^2 TX \otimes L^{*2}) \subset \mathbb{C}$. In fact, the long exact sequence

$$\begin{aligned} 0 \;&\rightarrow\; \mathrm{H}^0(\overset{0}{\bullet}\!-\!\overset{0}{\times}\!-\!\overset{0}{\bullet}\!-\!\overset{0}{\bullet}\cdots\overset{0}{\bullet}\!\Longleftarrow\!\overset{0}{\bullet}) \rightarrow \mathrm{H}^1(\overset{1}{\bullet}\!-\!\overset{-2}{\times}\!-\!\overset{1}{\bullet}\!-\!\overset{0}{\bullet}\cdots\overset{0}{\bullet}\!\Longleftarrow\!\overset{0}{\bullet}) \\ &\rightarrow\; \mathrm{H}^1(X, \odot^2 TX \otimes L^{*2}) \rightarrow 0 \end{aligned}$$

implies the vanishing of $\mathrm{H}^1(X, \odot^2 TX \otimes L^{*2})$. A similar analysis of $\odot^3 TX \otimes L^{*2}$ produces an estimation

$$\mathrm{H}^1(X, \odot^3 TX \otimes L^{*2}) \subset \mathrm{H}^1(\overset{1}{\bullet}\!-\!\overset{-2}{\times}\!-\!\overset{1}{\bullet}\!-\!\overset{0}{\bullet}\cdots\overset{0}{\bullet}\!\Longleftarrow\!\overset{0}{\bullet}) = \mathbb{C}$$

which, by (19), implies $K(\mathfrak{g}\oplus\mathbb{C}) \subset \mathrm{H}^1(X, L\otimes\odot^3 N^*) \subset \mathbb{C}$ which in turn implies $K^1(\mathfrak{g}\oplus\mathbb{C}) = 0$.

*Type* $F_4$ (iv). The associated Borel-Weil data are $(X,L) = \overset{0}{\bullet}\!-\!\overset{0}{\bullet}\!\Rightarrow\!\overset{0}{\bullet}\!-\!\overset{1}{\times}$. By (5), $TX = \overset{0}{\bullet}\!-\!\overset{0}{\bullet}\!\Rightarrow\!\overset{1}{\bullet}\!-\!\overset{-1}{\times} + \overset{1}{\bullet}\!-\!\overset{0}{\bullet}\!\Rightarrow\!\overset{0}{\bullet}\!-\!\overset{0}{\times}$, implying

$$\odot^2 TX\otimes L^{*2} = \begin{matrix}\overset{0}{\bullet}\!-\!\overset{0}{\bullet}\!\Rightarrow\!\overset{2}{\bullet}\!-\!\overset{-4}{\times}\\ \oplus \\ \overset{0}{\bullet}\!-\!\overset{0}{\bullet}\!\Rightarrow\!\overset{0}{\bullet}\!-\!\overset{-1}{\times}\end{matrix} + \begin{matrix}\overset{1}{\bullet}\!-\!\overset{0}{\bullet}\!\Rightarrow\!\overset{1}{\bullet}\!-\!\overset{-3}{\times}\\ \oplus \\ \overset{0}{\bullet}\!-\!\overset{0}{\bullet}\!\Rightarrow\!\overset{1}{\bullet}\!-\!\overset{-2}{\times}\end{matrix} + \begin{matrix}\overset{2}{\bullet}\!-\!\overset{0}{\bullet}\!\Rightarrow\!\overset{0}{\bullet}\!-\!\overset{-2}{\times}\\ \oplus \\ \overset{0}{\bullet}\!-\!\overset{0}{\bullet}\!\Rightarrow\!\overset{0}{\bullet}\!-\!\overset{0}{\times}\end{matrix}.$$

All the irreducible bundles on the right-hand side have the first cohomology vanishing except $\overset{0}{\bullet}\!-\!\overset{0}{\bullet}\!\Rightarrow\!\overset{1}{\bullet}\!-\!\overset{-2}{\times}$ whose first cohomology is $\mathbb{C}$. Therefore, $H^1(X,\odot^2 TX\otimes L^{*2}) \subset \mathbb{C}$. In fact, the long exact sequence

$$0 \to \mathrm{H}^0(\overset{0}{\bullet}\!-\!\overset{0}{\bullet}\!\Rightarrow\!\overset{0}{\bullet}\!-\!\overset{0}{\times}) \to \mathrm{H}^1(\overset{0}{\bullet}\!-\!\overset{0}{\bullet}\!\Rightarrow\!\overset{1}{\bullet}\!-\!\overset{-2}{\times}) \to \mathrm{H}^1(X,\odot^2 TX\otimes L^{*2}) \to 0$$

implies the vanishing of $\mathrm{H}^1(X,\odot^2 TX\otimes L^{*2})$. A similar analysis of $\odot^3 TX\otimes L^{*2}$ produces an estimation

$$\mathrm{H}^1(X,\odot^3 TX\otimes L^{*2}) \subset \mathrm{H}^1(\overset{0}{\bullet}\!-\!\overset{0}{\bullet}\!\Rightarrow\!\overset{1}{\bullet}\!-\!\overset{-2}{\times}) = \mathbb{C}$$

which, by (19), implies $K(\mathfrak{g}\oplus\mathbb{C}) \subset \mathrm{H}^1(X, L\otimes\odot^3 N^*) \subset \mathbb{C}$ which in turn implies $K^1(\mathfrak{g}\oplus\mathbb{C}) = 0$.

*Type* $F_4$ (iii). The associated Borel-Weil data are $(X,L) = \overset{0}{\bullet}\!-\!\overset{0}{\bullet}\!\Rightarrow\!\overset{1}{\times}\!-\!\overset{0}{\bullet}$. By (5),

$$TX = \overset{1}{\bullet}\!-\!\overset{0}{\bullet}\!\Rightarrow\!\overset{-1}{\times}\!-\!\overset{1}{\bullet} + \overset{0}{\bullet}\!-\!\overset{1}{\bullet}\!\Rightarrow\!\overset{-2}{\times}\!-\!\overset{2}{\bullet} + \overset{0}{\bullet}\!-\!\overset{0}{\bullet}\!\Rightarrow\!\overset{0}{\times}\!-\!\overset{1}{\bullet} + \overset{1}{\bullet}\!-\!\overset{0}{\bullet}\!\Rightarrow\!\overset{0}{\times}\!-\!\overset{0}{\bullet}.$$

By this formula for $TX$, it is straighforward but very tedious to obtain the associated (non-direct) decompositions into irreducible homogeneous bundles $E_i$ for $\odot^2 TX\otimes L^{*2}$ and $\odot^3 TX\otimes L^{*2}$, and to observe that all irreducible summonds $E_i$ have the first cohomology zero implying $\mathrm{H}^1(X,\odot^2 TX\otimes L^{*2}) = \mathrm{H}^1(X,\odot^2 TX\otimes L^{*2}) = 0$. By (19), this implies $\mathrm{H}^1(X, L\otimes\odot^3 N^*) = 0$ which in turn implies $K(\mathfrak{g}\oplus\mathbb{C}) = 0$.

*Type* $E_6$ (v). The associated Borel-Weil data are $(X,L) = \overset{0}{\bullet}\!-\!\overset{1}{\times}\!-\!\underset{\bullet\,0}{\overset{0}{\bullet}}\!-\!\overset{0}{\bullet}\!-\!\overset{0}{\bullet}$. By (5),

$$TX = \overset{1}{\bullet}\!-\!\overset{-1}{\times}\!-\!\underset{\bullet\,0}{\overset{0}{\bullet}}\!-\!\overset{1}{\bullet}\!-\!\overset{0}{\bullet} + \overset{0}{\bullet}\!-\!\overset{0}{\times}\!-\!\underset{\bullet\,1}{\overset{0}{\bullet}}\!-\!\overset{0}{\bullet}\!-\!\overset{0}{\bullet},$$

implying

$\odot^2 TX \otimes L^{*2} =$ 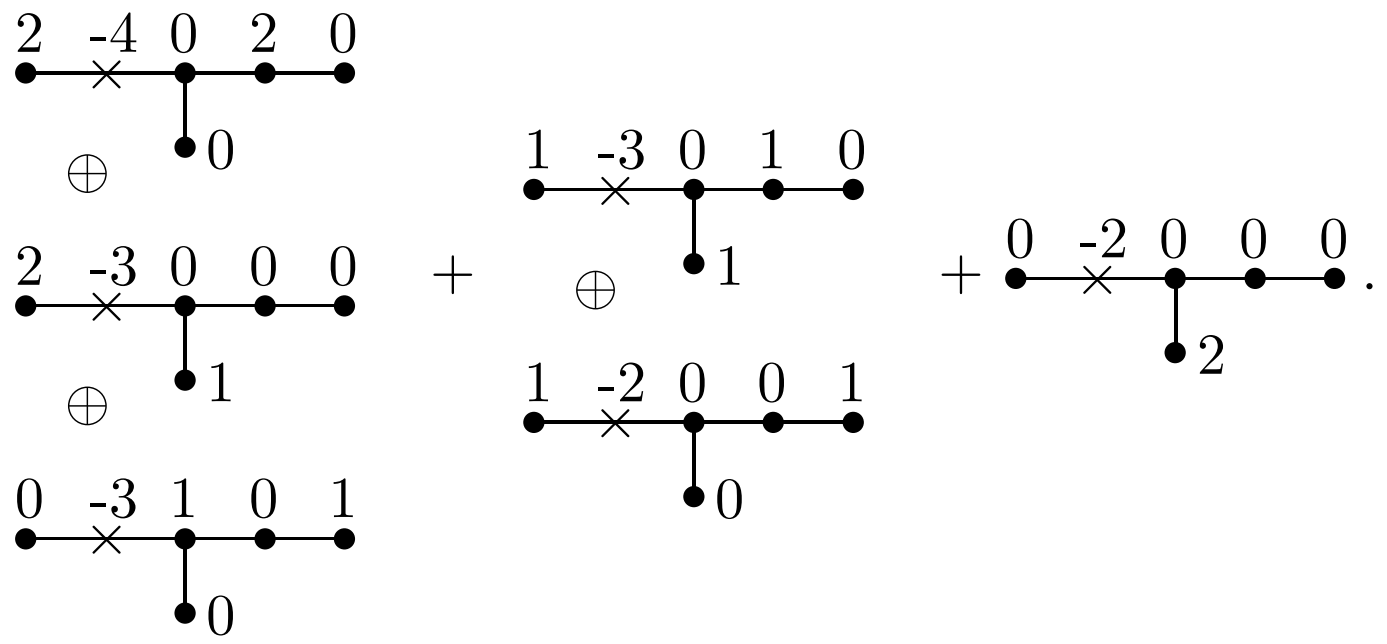

All irreducible bundles on the right-hand side of the above equality have the first cohomology vanishing, implying $\mathrm{H}^1(X, \odot^2 TX \otimes L^{*2}) = 0$. Analogously, one obtains $\mathrm{H}^1(X, \odot^3 TX \otimes L^{*2}) = 0$. By (19), $\mathrm{H}^1(X, L \otimes \odot^3 N^*) = 0$ implying $K(\mathfrak{g} \oplus \mathbb{C}) = 0$.

*Type $E_7$* (vi). The associated Borel-Weil data are

$(X, L) =$ 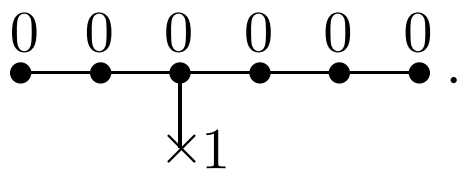

By (5),

$TX =$ 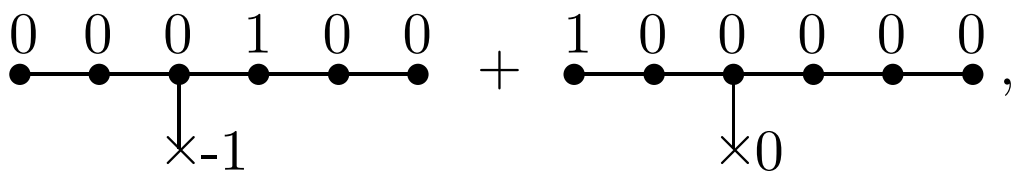

implying

$\odot^2 TX \otimes L^{*2} =$ 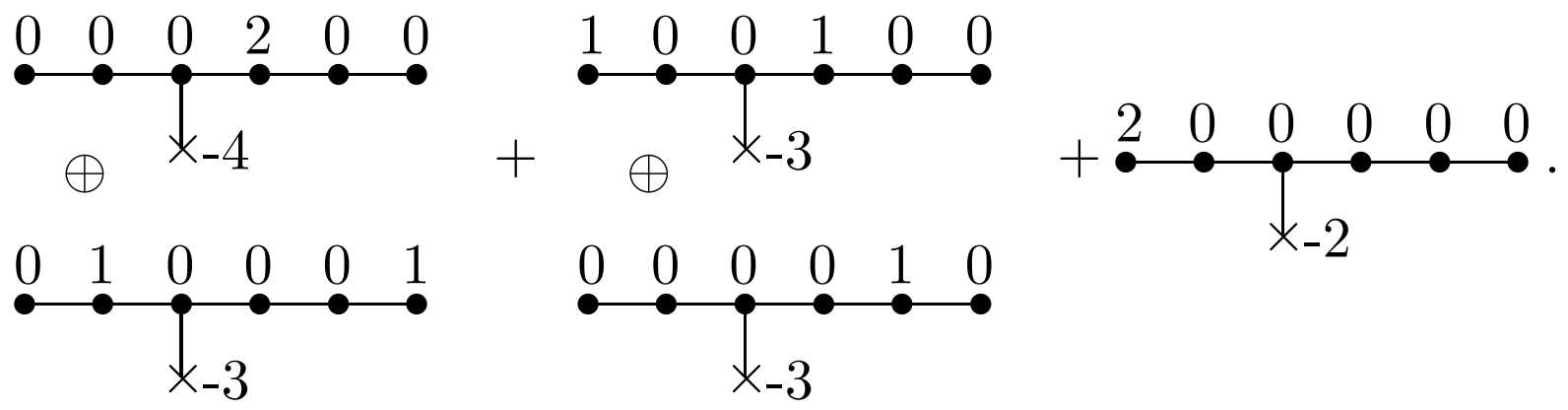

Since every irreducible bundle on the right-hand side of the above equality has the first cohomology vanishing, one obtains $\mathrm{H}^1(X, \odot^2 TX \otimes L^{*2}) = 0$. A similar analysis of $\odot^2 TX \otimes L^{*2}$ results in $\mathrm{H}^1(X, \odot^3 TX \otimes L^{*2}) = 0$. Then, by (19), $\mathrm{H}^1(X, L \otimes \odot^3 N^*) = 0$ implying $K(\mathfrak{g} \oplus \mathbb{C}) = 0$.

*Type $G_2$*. The Borel-Weil data associated with the representation 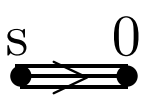 are $(X, L) =$ 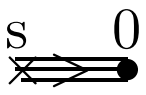. By (5), $TX =$ 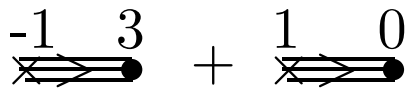, implying

$\odot^2 TX =$ 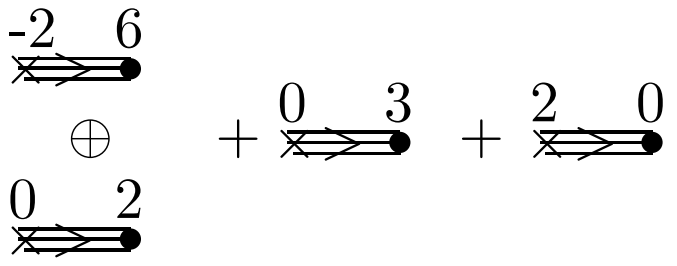

and

$$\odot^3 TX = \begin{array}{c} \overset{-3\quad 9}{\times\!\!\Rrightarrow\!\!\bullet} \\ \oplus \\ \overset{-1\quad 5}{\times\!\!\Rrightarrow\!\!\bullet} \\ \oplus \\ \overset{0\quad 3}{\times\!\!\Rrightarrow\!\!\bullet} \end{array} + \begin{array}{c} \overset{-1\quad 6}{\times\!\!\Rrightarrow\!\!\bullet} \\ \oplus \\ \overset{1\quad 2}{\times\!\!\Rrightarrow\!\!\bullet} \end{array} + \overset{1\quad 3}{\times\!\!\Rrightarrow\!\!\bullet} + \overset{3\quad 0}{\times\!\!\Rrightarrow\!\!\bullet}.$$

Thus $\mathrm{H}^1(X, \odot^2 TX \otimes L^{*2}) = 0$ for any $s \geqslant 1$ and

$$\mathrm{H}^1(X, \odot^3 TX \otimes L^{*2}) \subset \begin{cases} \mathbb{C} & \text{for } s = 1, \\ 0 & \text{otherwise.} \end{cases}$$

Now, $K(\mathfrak{g} \oplus \mathbb{C}) \subset \mathrm{H}^1(X, L \otimes \odot^3 N^*) \subset \mathbb{C}$ implying $K^1(\mathfrak{g} \oplus \mathbb{C}) = 0$.

Next we consider the representation $\overset{s\quad t}{\bullet\!\!\Rrightarrow\!\!\bullet}$, $s, t \geqslant 1$, which has the associated Borel-Weil data $(X, L) = \overset{s\quad t}{\times\!\!\Rrightarrow\!\!\times}$. There is a projection $\nu : X \to X_0 = \times\!\!\Rrightarrow\!\!\bullet$. One has $TX = \overset{-1\quad 2}{\times\!\!\Rrightarrow\!\!\times} + \nu^*(TX_0)$ with $TX_0 = \overset{-1\quad 3}{\times\!\!\Rrightarrow\!\!\bullet} + \overset{1\quad 0}{\times\!\!\Rrightarrow\!\!\bullet}$. Hence

$$\odot^2 TX = \overset{-2\quad 4}{\times\!\!\Rrightarrow\!\!\times} + \overset{-1\quad 2}{\times\!\!\Rrightarrow\!\!\times} \otimes \nu^*(TX_0) + \nu^*(TX_0).$$

Computing $\nu^i_*(\odot^2 TX \otimes \otimes L^{*2})$, $i = 0, 1$, one easily finds

$$\mathrm{H}^1(X, \odot^2 TX \otimes L^{*2}) \subset \begin{cases} \mathrm{H}^1(\overset{-1-2s\quad 0}{\times\!\!\Rrightarrow\!\!\bullet} \otimes TX_0) & \text{for any } s \text{ and } t = 1, \\ 0 & \text{otherwise.} \end{cases}$$

Since $\mathrm{H}^1(\overset{-1-2s\quad 0}{\times\!\!\Rrightarrow\!\!\bullet} \otimes TX_0) = 0$, we obtain $\mathrm{H}^1(X, \odot^3 TX \otimes L^{*2}) = 0$ for any $s, t \geqslant 1$. Analogous computation shows that $\mathrm{H}^1(X, \odot^3 TX \otimes L^{*2})$ vanishes for all $s, t \geqslant 1$ as well. Hence, by (19), $\mathrm{H}^1(X, L \otimes \odot^3 N^*) = 0$ implying $K(\mathfrak{g} \oplus C) = 0$.

The proof of Theorem 6.10 is completed.

6.4. *Proof of the Main Theorem.* Suppose $\mathfrak{g} \subset \mathfrak{gl}(V)$ is a real Lie subalgebra which occurs as the irreducible holonomy of a torsion-free affine connection which is not locally symmetric. Then $\mathfrak{g}$ is a Berger algebra by Proposition 2.2. By Theorem 5.1, we may assume that the semisimple part of $\mathfrak{g}$ is simple. Moreover, by Proposition 2.4, it follows that $K(\mathfrak{g}_{\mathbb{C}}) \neq 0$. Therefore, $\mathfrak{g}_{\mathbb{C}}$ must be one of the entries of Table 6. So, in order to complete the proof of the Main Theorem, one first determines all possible real forms of real type of the entries of Table 6, and then investigates which real forms of the complex Lie algebras with $(\mathfrak{g}_{\mathbb{C}})^{(1)} \neq 0$ (cf. Table 5) satisfy Berger's criteria.

For the entries 1–3, 6, 7, 9 and 12, this has been carried out by Berger [8], Bryant [14], [16] and McClean [16]; all admissible real forms of entries 4, 10 and 13 have also been studied by Bryant [15], [16]. Finally, the entries 5, 8, 11 and 14 satisfy $(\mathfrak{g}_{\mathbb{C}})^{(1)} = 0$, and thus Table 3 contains all real forms of real type for these entries. That all of these do occur as holonomies will be shown in Section 7. □

This argument immediately implies the following corollary, completing the classification of complex Berger algebras:

COROLLARY 6.11. *Let* $\mathfrak{g} \subset \mathfrak{gl}(V)$ *be an irreducible complex Lie subalgebra whose semisimple part is simple. Then* $\mathfrak{g}$ *is a Berger algebra if and only it is an entry of Table* 7.

Since all our arguments above relied essentially on the assumption $\dim K(\mathfrak{g}) > 1$, we obtained actually a stronger result than the one formulated in the above corollary:

THEOREM 6.12. *Let* $\mathfrak{g} \subset \mathfrak{gl}(V)$ *be an irreducible complex semisimple Lie subalgebra. Then* $\dim K(\mathfrak{g} \oplus \mathbb{C}) > 1$ *if and only if* $\mathfrak{g}$ *is* (*the semisimple part of*) *an entry of Table* 6 *or a complex entry of Table* 7.

## 7. Existence and supersymmetry of exotic holonomies

7.1. *Existence theorem.* Let $V$ be a $2n$-dimensional vector space over $\mathbb{F} = \mathbb{R}$ or $\mathbb{C}$ equipped with a non-degenerate symplectic form $\omega$. Denote by $\mathrm{Sp}(n, \mathbb{F})$ the subgroup of $\mathrm{GL}(V)$ which preserves $\omega$ and by $\mathfrak{sp}(n, \mathbb{F})$ the Lie algebra of $\mathrm{Sp}(n, \mathbb{F})$.

Since $\mathfrak{sp}(n, \mathbb{F}) \simeq \odot^2 V^*$, for any irreducible subalgebra $\mathfrak{g} \subset \mathfrak{sp}(n, \mathbb{F})$ one has a nonzero $\mathfrak{g}$-invariant symmetric map

$$\begin{array}{ccc} \odot^2 V & \longrightarrow & \mathfrak{g} \\ u \otimes v & \longrightarrow & u \circ v, \end{array}$$

which is unique up to nonzero scalar factor.

Let $B : \odot^2 \mathfrak{g} \to \mathbb{F}$ be the Killing form of $\mathfrak{g}$.

LEMMA 7.1. *The products* $\omega : \Lambda^2 V \to \mathbb{F}$ *and* $\circ : \odot^2 V \to \mathfrak{g}$ *satisfy, for some nonzero constants* $\lambda, \mu \in \mathbb{F}$ *, the identities*

$$\omega(Au, v) = \lambda B(A, u \circ v), \tag{40}$$

$$B(u \circ v, s \circ t) - B(u \circ t, s \circ v) = \mu(2\omega(u,s)\omega(v,t) + \omega(u,t)\omega(v,s) + \omega(u,v)\omega(s,t)), \tag{41}$$

*for all* $A \in \mathfrak{g}$ *and all* $u, v, s, t \in V$, *if and only if the pair* $(G, V)$ *is an entry of the table of Theorem* C.

*Proof.* Both identities are trivially satisfied for $\mathfrak{g} = \mathfrak{sp}(n, \mathbb{F})$.

The first identity holds for an arbitrary $\mathfrak{g} \subset \mathfrak{sp}(n, \mathbb{F})$, since the Killing form $B$ on $\mathfrak{g}$ coincides with the restriction of the Killing form on $\mathfrak{sp}(n, \mathbb{F})$ which is proportional to the metric $\omega \otimes \omega$ induced on $\odot^2 V^*$ from $V^* \otimes V^*$.

A straightforward calculation shows that the second identity (41) is equivalent to saying that, for any fixed $A \in \mathfrak{g}$, the map

$$\tag{42} \begin{array}{rcl} \rho_A : \ \Lambda^2 V & \longrightarrow & \mathfrak{g} \\ u \otimes v & \longrightarrow & \mu\, \omega(u, v)\, A - u \circ (Av) + v \circ (Au) \end{array}$$

defines an element of $K(\mathfrak{g})$, i.e. lies in the kernel of the composition

$$\mathfrak{g} \otimes \Lambda^2 V^* \longrightarrow V \otimes V^* \otimes \Lambda^2 V^* \longrightarrow V \otimes \Lambda^3 V^*.$$

Therefore, if (41) is satisfied, then $K(\mathfrak{g})$ contains an irreducible summand isomorphic to $\mathrm{Ad}(\mathfrak{g})$. So $\mathfrak{g}$ can only be a representation listed in the table of Theorem C.

The $\mathfrak{sp}(n, \mathbb{F})$-module $K(\mathfrak{sp}(n, \mathbb{F}))$ contains an irreducible submodule $\mathrm{Ad}(\mathfrak{sp}(n, \mathbb{F}))$. Since the above identities do hold for $\mathfrak{g} = \mathfrak{sp}(n, \mathbb{F})$, the elements of this submodule are given explicitly by the formula (42). As follows from Theorem 6.12, the only irreducible symplectic subgroups $\mathfrak{g} \subset \mathfrak{sp}(n, \mathbb{F})$ which have $\dim K(\mathfrak{g}) > 1$ are the ones given in the table of Theorem C. They all have $K(\mathfrak{g}) \simeq \mathrm{Ad}(\mathfrak{g})$ so that their elements must be given explicitly by the same formula (42). Then the Bianchi identities for maps (42) imply (41). □

For $\mathfrak{g} = \mathfrak{e}_7^{\mathbb{C}}$ represented in $\mathbb{C}^{56}$ the identities (40) and (41) were first established by Adams [1] via a reduction from $\mathfrak{e}_8^{\mathbb{C}}$. We refer to [21] for details.

Having obtained in the proof of Lemma 7.1 an explicit structure of $K(\mathfrak{g})$ for all proper subgroups $\mathfrak{g}$ of $\mathfrak{sp}(n, \mathbb{F})$ listed in the table of Theorem C, we can show easily that a generic element of $K^1(\mathfrak{g}) \subset \mathfrak{g} \otimes V^* \otimes \Lambda^2 V^*$ is of the form

$$\begin{array}{rcl} V \otimes \Lambda^2 V & \longrightarrow & \mathfrak{g} \\ s \otimes u \otimes v & \longrightarrow & \mu\, \omega(u, v)\, (s \circ w) - (u \circ (s \circ w))v + (v \circ (s \circ w))u \end{array}$$

for some fixed $w \in V \simeq V^*$. This establishes the isomorphism $K^1(\mathfrak{g}) = V^*$.

Lemma 7.2. *The $\mathfrak{g}$-module $S(\mathfrak{g}) = \odot^2 \mathfrak{g} \otimes \Lambda^2 V^* \cap \mathfrak{g} \otimes K(\mathfrak{g})$ contains a one-dimensional submodule if and only if $G$ is one of the entries of the table of Theorem* C.

*Proof.* First note that the $\mathfrak{g}$-module $K(\mathfrak{g}) \otimes \mathfrak{g} \simeq \mathfrak{g} \otimes \mathfrak{g}$ has only one 1-dimensional $\mathfrak{g}$-submodule. So, if there is an invariant element $\phi$ in $S(\mathfrak{g})$, it is unique up to a nonzero scalar factor. Since $S(\mathfrak{g}) \subset V \otimes K^1(\mathfrak{g})$ such an element gives rise to $\mathfrak{g}$-equivariant maps

$$\tag{43} \phi' : \mathfrak{g}^* \longrightarrow K(\mathfrak{g}),$$

$$\tag{44} \phi'' : V^* \longrightarrow K^1(\mathfrak{g}),$$

implying, in particular, that $\mathrm{Ad}(\mathfrak{g}) \subset K(\mathfrak{g})$. Thus, if $S(\mathfrak{g})$ has a one-dimensional submodule, then $\mathfrak{g}$ must be in the table of Theorem C.

Let us now construct a nonzero $\mathfrak{g}$-invariant element $\phi \in S(\mathfrak{g})$ for every entry of that table. Formula (42) gives no choice but to consider the following element of $\mathfrak{g} \otimes K(\mathfrak{g})$ as a candidate to $\phi$:

$$\phi(C, D, u, v) = \mu\, \omega(u, v)\, B(C, D) + B(u \circ (Cv), D) - B(v \circ (Cu), D),$$

where $C, D \in \mathfrak{g}$, $u, v \in V$ and we identify $\mathfrak{g} = \mathfrak{g}^*$ via the Killing form. Clearly, this element is $\mathfrak{g}$-invariant. Since $B(u \circ (Cv), D) = \lambda\, \omega(Du, Cv)$, we have

$$\begin{aligned} \phi(C, D, u, v) &= \mu\, \omega(u, v)\, B(C, D) + \lambda\, \omega(Du, Cv) - \lambda\, \omega(Dv, Cu) \\ &= \mu\, \omega(u, v)\, B(C, D) - \lambda\, \omega(Cv, Du) - \lambda\, \omega(Dv, Cu), \end{aligned}$$

which makes it evident that $\phi \in (\odot^2\mathfrak{g} \otimes \Lambda^2 V^*) \cap (\mathfrak{g} \otimes K(\mathfrak{g}))$. That $\phi' : \mathfrak{g}^* \longrightarrow K(\mathfrak{g})$ is an isomorphism follows from the very definition of $\phi$. Since $\phi''_2 : V^* \longrightarrow K^1(\mathfrak{g})$ is evidently nonzero, it is an isomorphism as well. □

*Proof of Theorem* A. It was shown in [20] that all statements of Theorem A hold for those $\mathfrak{g} \subset \mathfrak{gl}(V)$ whose associated module $S(\mathfrak{g})$ contains an invariant element $\phi$ for which the corresponding maps (43) and (44) are isomorphisms.

Since such an element $\phi$ was given above for all representations listed in the table of Theorem C, and thus, in particular, for the new series of exotic holonomies, Theorem A follows. □

This argument also shows that the statements of Theorem A are true for all torsion-free affine connections $\nabla$ with "generic" symplectic holonomy $\mathfrak{sp}(n, \mathbb{F})$ satisfying the condition that the curvature tensor $R^\nabla$ takes values in the submodule $\mathrm{Ad}(\mathfrak{sp}(n, \mathbb{F})) \subset K(\mathfrak{sp}(n, \mathbb{F}))$.

Finally, we show the following result.

LEMMA 7.3. *If $\mathfrak{spin}(11, \mathbb{C}) \subset \mathfrak{gl}(V)$ with $V \cong \mathbb{C}^{16}$ is the image of the spinor representation, then* $K(\mathfrak{spin}(11, \mathbb{C})) = K(\mathfrak{spin}(11, \mathbb{C}) \oplus \mathbb{C}\,\mathrm{Id}) = 0$.

*Proof.* Since $\mathfrak{spin}(11, \mathbb{C}) \subset \mathfrak{spin}(12, \mathbb{C}) \subset \mathfrak{gl}(V)$, it follows that $K(\mathfrak{spin}(11, \mathbb{C}) \oplus \mathbb{C}\,\mathrm{Id}) \subset K(\mathfrak{spin}(12, \mathbb{C}) \oplus \mathbb{C}\,\mathrm{Id}) = K(\mathfrak{spin}(12, \mathbb{C})) \cong \mathrm{Ad}(\mathfrak{spin}(12))$ by the above result. In particular, $K(\mathfrak{spin}(11, \mathbb{C}) \oplus \mathbb{C}\,\mathrm{Id}) = K(\mathfrak{spin}(11, \mathbb{C}))$.

But now, a calculation shows that there is no nonzero element $A \in \mathfrak{spin}(12, \mathbb{C})$ such that $\rho_A(u, v) \in \mathfrak{spin}(11, \mathbb{C})$ for all $u, v \in V$ with $\rho_A$ as in (42). Thus, $K(\mathfrak{spin}(11, \mathbb{C})) = 0$. □

7.2. *Supersymmetry.* A (smooth, analytic, etc.) *supermanifold* of dimension $(m|n)$ is a locally ringed space $(M, \mathcal{O}_M)$ with the following properties [29]:

(i) the structure sheaf $\mathcal{O}_M = \mathcal{O}_{M,0} \oplus \mathcal{O}_{M,1}$ is a sheaf of $\mathbb{Z}_2$-graded supercommutative rings; (ii) $M_{\text{red}} = (M, \mathcal{O}_{M,\text{red}} := \mathcal{O}_M/[\mathcal{O}_{M,1} + \mathcal{O}^2_{M,1}])$ is a (smooth, analytic, etc.) classical manifold of dimension $m$; (iii) $\mathcal{O}_M$ is locally isomorphic to the exterior algebra $\Lambda(E)$ of a locally free $\mathcal{O}_{M,\text{red}}$-module $E$ of rank $n$. If $\phi : \Lambda(E) \to \mathcal{O}_M$ is any such local isomorphism, $\bar{x}^1, \ldots, \bar{x}^m$ are local coordinates on $M_{\text{red}}$ and $\bar{\theta}^1, \ldots, \bar{\theta}^n$ are free local generators of $E$, then the set of $m+n$ sections

$$x^1 = \phi(\bar{x}^1), \ldots, x^m = \phi(\bar{x}^m),\ \theta^1 = \phi(\bar{\theta}^1), \ldots, \theta^n = \phi(\bar{\theta}^n)$$

of the structure sheaf $\mathcal{O}_M$ form a local coordinate system on $M$. Any local function $f$ on $M$ can be expressed as a polynomial in anticommuting odd coordinates $\theta^\alpha$, $\alpha = 1, \ldots, n$,

$$f(x, \theta) = \sum_{k=0}^{n} \sum_{\alpha_1, \ldots, \alpha_k = 1}^{n} f_{\alpha_1 \ldots \alpha_k}(x) \theta^{\alpha_1} \cdots \theta^{\alpha_k}$$

whose coefficients $f_{\alpha_1 \ldots \alpha_k}(x)$ are classical (smooth, analytic, etc.) functions of the commuting variables $x^a$, $a = 1, \ldots, m$.

When a need arises to use odd constants in the structure sheaf of a supermanifold $M$, one simply replaces $M$ by its relative version, i.e. by a submersion of supermanifolds $\pi : \mathcal{M} \to \mathcal{S}$ whose typical fibre is $M$. Then (odd) constants are just (odd) elements of $\pi^{-1}(\mathcal{O}_{\mathcal{S}})$. Such a replacement is made often tacitly — the necessary changes are routine; cf. [29], [30].

A SUSY-*structure* on an $(m|n)$-dimensional supermanifold $M$ is, by definition, a rank $0|n$ locally split subsheaf $\mathcal{T}_1 \subset TM$ such that the associated Frobenius form

$$\begin{aligned} \Phi : \quad \Lambda^2 \mathcal{T}_1 &\longrightarrow \mathcal{T}_0 := TM/\mathcal{T}_1 \\ X \otimes Y &\longrightarrow [X, Y] \bmod \mathcal{T}_1 \end{aligned}$$

is surjective.

*Examples.* (i) A $(1|1)$-dimensional supermanifold $M$ with a SUSY-structure is the same as a $\text{SUSY}_1$-curve [30]. (ii) A SUSY-structure on a $(3|2)$-dimensional supermanifold is equivalent to $N = 1$ superconformal supergravity in 3 dimensions [33]. (iii) A SUSY-structure on a $(4|4)$-dimensional supermanifold with $\mathcal{T}_1$ being a direct sum of two integrable rank $(0|2)$ distributions $\mathcal{T}_l$ and $\mathcal{T}_r$ is the same as an $N = 1$ superconformal supergravity in 4 dimensions [29].

Let $V$ be a $2n$-dimensional vector space over $\mathbb{F} = \mathbb{R}$ or $\mathbb{C}$ equipped with a nondegenerate symplectic form $\omega$ and $\mathfrak{g} \subseteq \mathfrak{sp}(n, \mathbb{F})$ be an irreducible Lie subalgebra. Define $M_{\mathfrak{g}}$ as a linear supermanifold $\mathfrak{g} \oplus \Pi(V^*)$, where $\Pi$ is the parity change functor [29]. Then, as usual for linear (super)manifolds, there is an isomorphism $i : [\mathfrak{g} \oplus \Pi(V^*)] \otimes \mathcal{O}_M \to TM_{\mathfrak{g}}$ so that $TM$ comes equipped with a free subsheaf $i(\Pi(V^*) \otimes \mathcal{O}_M)$ of rank $0|n$. Let $(e_a)$, $a = 1, \ldots, m = \dim \mathfrak{g}$, be a basis in the vector space $\mathfrak{g}$ and $(e_\alpha)$, $\alpha = 1, \dots, 2n$, be a basis in the vector space $\Pi(V^*)$, and let $(x^a)$ and $(\theta^\alpha)$ be the associated commuting and, respectively anticommuting linear coordinates. Then the system $(x^a, \theta^\alpha)$ is a global coordinate system on $M_{\mathfrak{g}}$. Since $i(\Pi(V^*) \otimes \mathcal{O}_M)$ is freely generated by $\partial_\alpha = \partial/\partial\theta^\alpha$, the associated Frobenius form is zero, implying that this subsheaf is *not* a SUSY-structure on $M$. This can be corrected by the following natural map

$$\begin{array}{rclcl} \Phi:\ \Lambda^2(\Pi(V^*)) & \xrightarrow{\simeq} & \odot^2 V^* & \longrightarrow & \mathfrak{g} \\ e_\alpha \wedge e_\beta & \longrightarrow & \Pi(e_\alpha)\odot\Pi(e_\beta) & \longrightarrow & \phi^a_{\alpha\beta} e_a, \end{array}$$

which gives rise to a monomorphism

$$\begin{array}{rcl} j: \Pi(V^*)\otimes\mathcal{O}_M & \longrightarrow & TM \\ e_\alpha & \longrightarrow & \partial_\alpha + \phi^a_{\alpha\beta}\theta^\beta\partial_a, \end{array}$$

where $\partial_a = \partial/\partial x^a$ and summation over repeated indices is assumed.

LEMMA 7.4. *The subsheaf* $\mathcal{T}_1 = j(\Pi(V^*)) \subset TM_{\mathfrak{g}}$ *is a* SUSY-*structure on* $M_{\mathfrak{g}}$.

*Proof.* The subsheaf $\mathcal{T}_1$ is freely generated by the vector fields $D_\alpha = \partial_\alpha + \phi^b_{\alpha\beta}\theta^\beta\partial_b$. A simple calculation shows $[D_\alpha, D_\beta] = 2\phi^d_{\alpha\beta}\partial_d$ implying that the Frobenius form of $\mathcal{T}_1$ is equal to $\Phi$ and is, therefore, nondegenerate. □

The SUSY-structure on $M_{\mathfrak{g}}$ gives rise to a group $\mathcal{G}$ of SUSY-translations $x^a \longrightarrow x^a + \phi^a_{\alpha\beta}\epsilon^\alpha\theta^\beta$ which mix even and odd coordinates, $\epsilon^\alpha$ being some odd constants. A function $f$ on $M_{\mathfrak{g}}$ is said to be SUSY-*invariant* if it is infinitesimally invariant under $\mathcal{G}$, i.e. if it satisfies the differential equation $D_\alpha f = \partial_\alpha f$.

*Proof of Theorem C.* Let $f$ be a nonzero $\mathfrak{g}$- and SUSY-invariant polynomial of order two in $\theta^\alpha$ and of order at least two in $x^a$; i.e.,

$$f = \sum_{k\geqslant 2} \phi_{a_1\cdots a_k\,\mu\nu} x^{a_1}\dots x^{a_k}\theta^\mu\theta^\nu$$

for some $\phi_{a_1\cdots a_k\,\mu\nu} \in \odot^k\mathfrak{g}\otimes\Lambda^2(V^*)$. An easy calculation shows that $f$ satisfies the equation $D_\alpha f = \partial_\alpha f$ if and only if $\phi_{a_1\cdots a_k\,\mu\nu} \in \odot^k\mathfrak{g}\otimes\Lambda^2(V^*)\cap\odot^{k-1}\mathfrak{g}\otimes K(\mathfrak{g})$. Then, since $f$ is $\mathfrak{g}$-invariant, $K(\mathfrak{g})$ must contain $\mathfrak{g}$-submodules $\odot^{k-1}\mathfrak{g}$, $k \geqslant 2$. By Theorem 6.12, $\mathfrak{g}$ can only be an entry of the table of Theorem C. As Lemma 7.2 shows, every such entry $\mathfrak{g}$ does admit a nonzero invariant element $\phi_{a_1 a_2\,\mu\nu} \subset S(\mathfrak{g})$ so that the set of $\mathfrak{g}$- and SUSY-invariant functions of the required type on the associated supermanifold $M_{\mathfrak{g}}$ is nonempty. □

## 8. Classification of complex homogeneous-rational manifolds $X$

**and ample line bundles $L \to X$ such that $\mathrm{H}^1(X, TX \otimes L^*) \neq 0$**

*Proof of Theorem* B. Statement (i) of Theorem B follows from Proposition 5.2. Let us prove statement (ii).

If $X$ is reducible, say $X = X_1 \times X_2$ and $L = \pi_1^*(L_1) \otimes \pi_2^*(L_2)$, then

$$\begin{aligned} \mathrm{H}^1(X, TX \otimes L^*) &= \mathrm{H}^0(X_1, TX_1 \otimes L_1^*) \otimes \mathrm{H}^1(X_2, L_2^*) \\ &\quad + \mathrm{H}^0(X_2, TX_2 \otimes L_2^*) \otimes \mathrm{H}^1(X_1, L_1^*). \end{aligned}$$

This together with statement (i) implies that in the class of reducible $X$ only two bottom lines in Table 3 contribute to the list of all $(X, L)$ with $\mathrm{H}^1(X, TX \otimes L^*) \neq 0$.

Assume from now on that $X$ is irreducible. Though the tangent bundle $TX$ is homogeneous, it is not irreducible in general; even worse, since the parabolic $P$ is not reductive, $TX$ is *not* a direct sum of irreducible homogeneous vector bundles. This makes a naive idea of computing $\mathrm{H}^*(X, TX \otimes L^*)$ by the straightforward application of the Bott theorem impractical.

Consider the Atiyah exact sequence

$$0 \longrightarrow Q \longrightarrow \mathfrak{g} \otimes \mathcal{O}_X \longrightarrow TX \longrightarrow 0, \tag{45}$$

where $Q = G \times_{\mathrm{Ad}} \mathfrak{p}$. Since the central term of this extension is a trivial vector bundle and $\mathrm{H}^i(X, L^*) = 0$ for $0 \leqslant i \leqslant \dim X - 1$, we have, in case $\dim X \geqslant 3$,

$$\mathrm{H}^1(X, TX \otimes L^*) = \mathrm{H}^2(X, Q \otimes L^*).$$

An exact sequence of $\mathfrak{p}$-modules

$$0 \longrightarrow \tilde{\mathfrak{n}} \longrightarrow \mathfrak{p} \longrightarrow \mathfrak{p}/\tilde{\mathfrak{n}} \longrightarrow 0,$$

where $\tilde{\mathfrak{n}} = \mathfrak{n} \setminus \sum_{\alpha \in \Phi_{\mathfrak{p}}^+} \mathfrak{g}_\alpha$, gives rise to an exact sequence

$$0 \longrightarrow \Omega^1 X \longrightarrow Q \longrightarrow S \longrightarrow 0$$

of homogeneous vector bundles, where $S = G \times_{\mathrm{Ad}} \mathfrak{p}/\tilde{\mathfrak{n}}$ and where the isomorphism $G \times_{\mathrm{Ad}} \tilde{\mathfrak{n}} \simeq \Omega^1 X$ was used. According to Nakano [34], for any compact complex manifold $X$ and any positive line bundle $L$ on $X$ the groups $\mathrm{H}^i(X, \Omega^1 X \otimes L^*) = 0$ vanish for all $i \leqslant \dim X - 2$. Then, in the case $\dim X \geqslant 5$, the long exact sequence of the latter extension implies

$$\mathrm{H}^2(X, Q \otimes L^*) = \mathrm{H}^2(X, S \otimes L^*)$$

which in turn implies

$$\mathrm{H}^1(X, TX \otimes L^*) = \mathrm{H}^2(X, S \otimes L^*).$$

The advantage of working with $S$ instead of $TX$ is that $S$ can always be decomposed into a *direct* sum of irreducible homogeneous subbundles.

The Lie algebra $\mathfrak{s} = \mathfrak{p}/\tilde{\mathfrak{n}}$ is reductive and the adjoint representation of $\mathfrak{p}$ on $\mathfrak{s}$ is semisimple. Under the adjoint representation $\mathfrak{s} \to \mathfrak{gl}(\mathfrak{s})$ the Lie algebra $\mathfrak{s}$ decomposes into a direct sum of its ideals

$$\mathfrak{s} = \xi_1 + \ldots \xi_k + \mathfrak{s}_1 + \ldots + \mathfrak{s}_m$$

where $\xi_j$, $j = 1, \ldots, \mathrm{rank} X$, lie in the center of $\mathfrak{s}$ and the non-Abelian ideals $\mathfrak{s}_i$, $i = 1, \ldots, m$, are simple. Then, by Bott's theorem,

$$\mathrm{H}^2(X, S \otimes L^*) = \bigoplus_{j=1}^{k} J^2(-\lambda) + \bigoplus_{i=1}^{m} J^2(\mu_i - \lambda),$$

where $\lambda$ is the weight of $L$ and $\mu_1, \ldots, \mu_m$ are the maximal roots of the simple ideals $\mathfrak{s}_1, \ldots, \mathfrak{s}_m$ [11], [39]. Since, for $\dim X \geqslant 2$, $J^2(-\lambda) = \mathrm{H}^2(X, L^*) = 0$, we obtain the following:

Lemma 8.1. *If* $\dim X \geqslant 5$, *then* $\mathrm{H}^1(X, TX \otimes L^*) = \bigoplus_{i=1}^{m} J^2(\mu_i - \lambda)$.

There are seven irreducible compact complex homogeneous-rational manifolds $X$ with $\dim X \leqslant 4$: projective spaces $\mathbb{CP}_k$ for $k = 1, 2, 3, 4$, quadrics $Q_3$, $Q_4$ and the complete flag manifold $F(1, 2; \mathbb{C}^3)$. It is elementary to check that Theorem B is true for this family.

We assume from now on that $X$ is an irreducible complex homogeneous-rational manifold with $\dim X \geqslant 5$.

Let us introduce the following notation: If $\Gamma$ is a connected subgraph of the Dynkin diagram for $X$, then the number of simple roots $\{\alpha_j \mid j \in I\}$ in this graph is denoted by $|\Gamma|$; if $\omega$ is an integral weight such that $\omega(H_{\alpha_j}) \leqslant 0$ for all $\alpha_j \in \Gamma$, then we write $\omega|_\Gamma \leqslant 0$.

Lemma 8.2. $\mathrm{H}^1(X, TX \otimes L^*) = 0$ *for any ample line bundle* $L$ *on* $X$ *if at least one of the following conditions is satisfied*:

(i) $\mathrm{rank} X \geqslant 3$;

(ii) $\mathrm{rank} X = 2$ *and the crossed nodes are adjacent*;

(iii) $\mathrm{rank} X = 2$, *the crossed nodes are not adjacent and, for each maximal root* $\mu_i$ *of the simple ideal* $\mathfrak{s}_i$, $i = 1, \ldots, m$, *at least one crossed node is contained in a connected subgraph* $\Gamma$ *of the Dynkin diagram for* $X$ *such that* $|\Gamma| \geqslant 2$ *and* $\mu_i|_\Gamma \leqslant 0$.

(iv) $\mathrm{rank} X = 1$ *and, for each maximal root* $\mu_i$ *of the simple ideal* $\mathfrak{s}_i$, $i = 1, \ldots, m$, *the crossed node is contained in a connected subgraph* $\Gamma$ *of the Dynkin diagram for* $X$ *such that* $|\Gamma| \geqslant 3$ *and* $\mu_i|_\Gamma \leqslant 0$.

*Proof.* (i) Let $\lambda$ be the weight of $L$. Since $L$ is ample, the coefficient of $\lambda$ over each crossed node is a negative integer (its coefficient over each uncrossed node is, of course, zero). Then $(-\lambda+\mu_i+\eta)(H_{\alpha_j}) \leqslant -\lambda(H_{\alpha_j})+1 \leqslant 0$ for all crossed nodes $\alpha_j$ and all $i \in \{1,\dots,m\}$. If the number of crossed nodes is greater than or equal to 3, then either $-\lambda+\mu_i+\eta$ is singular or $\mathrm{ind}(-\lambda+\mu_i+\eta) \geqslant 3$. Whence $\bigoplus_{i=1}^m J^2(\mu_i-\lambda)=0$ and the statement follows from Lemma 8.1.

(ii) If $\alpha_j$ and $\alpha_{j+1}$ are adjacent crossed nodes, then $\alpha_j+\alpha_{j+1}$ is a positive root and one has $(-\lambda+\mu_i+\eta)(H_{\alpha_j}) \leqslant 0$, $(-\lambda+\mu_i+\eta)(H_{\alpha_{j+1}}) \leqslant 0$ and hence $(-\lambda+\mu_i+\eta)(H_{\alpha_j+\alpha_{j+1}}) \leqslant 0$. Thus either $-\lambda+\mu_i+\eta$ is singular or $\mathrm{ind}(-\lambda+\mu_i+\eta) \geqslant 3$ for all $i$ and the statement follows from Lemma 8.1.

(iii) and (iv) If $\Gamma'$ is a connected subgraph of the Dynkin diagram for $X$, then the sum of all simple roots in $\Gamma'$ is a positive root [12]. Under the conditions stated in (iii) and (iv), one easily finds at least three positive roots $\alpha_j$ such that $(-\lambda+\mu_i+\eta)(H_{\alpha_j}) \leqslant 0$ for all $i \in \{1,\dots,m\}$. Then, again, either $-\lambda+\mu_i+\eta$ is singular or $\mathrm{ind}(-\lambda+\mu_i+\eta) \geqslant 3$ implying $J^2(-\lambda+\mu_i)=0$. Thus the statement follows from Lemma 8.1. □

Therefore, we can restrict our attention to the cases $\mathrm{rank}X = 1, 2$.

*The case* rank $X=2$. It clear from items (ii) and (iii) of Lemma 8.2 that $\mathrm{H}^1(X, TX\otimes L^*) \neq 0$ only for those $X$ which have the semisimple part $\mathfrak{s}'$ of the parabolic algebra $\mathfrak{p}$ simple, i.e. the number $m$ of simple ideals of $\mathfrak{s}'$ is 1. Therefore, the crossed nodes must be located at the ends of the Dynkin diagram for $X$. Inspecting the list of maximal roots (4) leaves one with the following three candidates to the list of all $(X, L)$ with $\mathrm{H}^1(X, TX\otimes L^*) \neq 0$:

(1) $(X, L) =$ 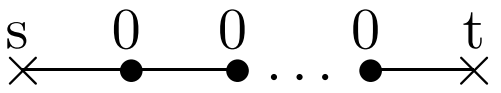

for some s,t$\geqslant$ 1. Computing $\sigma_1\circ\sigma_n(-\lambda+\mu_1+\eta)-\eta$ as shown in the following diagram

$-\lambda+\mu_1 =$ 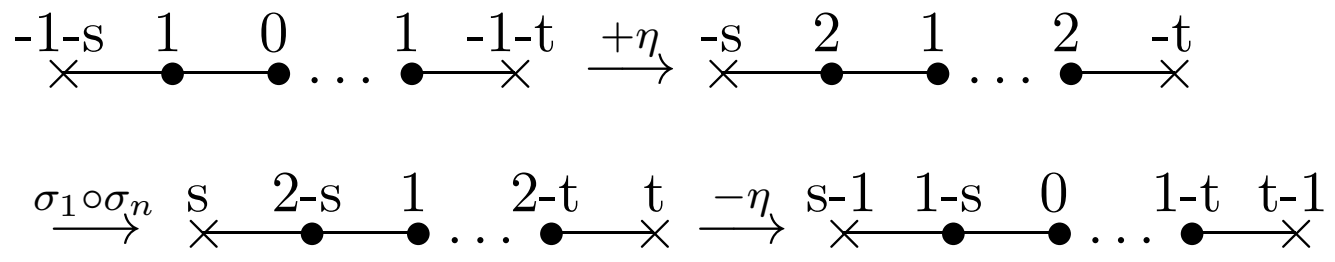

one concludes

$$\mathrm{H}^1(X, TX\otimes L^*) = J^2(-\lambda+\mu_1) = \begin{cases} \mathbb{C} & \text{s=t=1} \\ 0 & \text{otherwise.} \end{cases}$$

(2) $(X, L) =$ 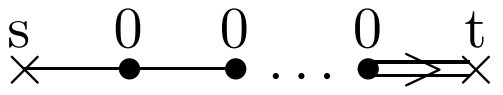

for some s,t$\geqslant$ 1. Then the graph

$-\lambda+\mu_1 =$ 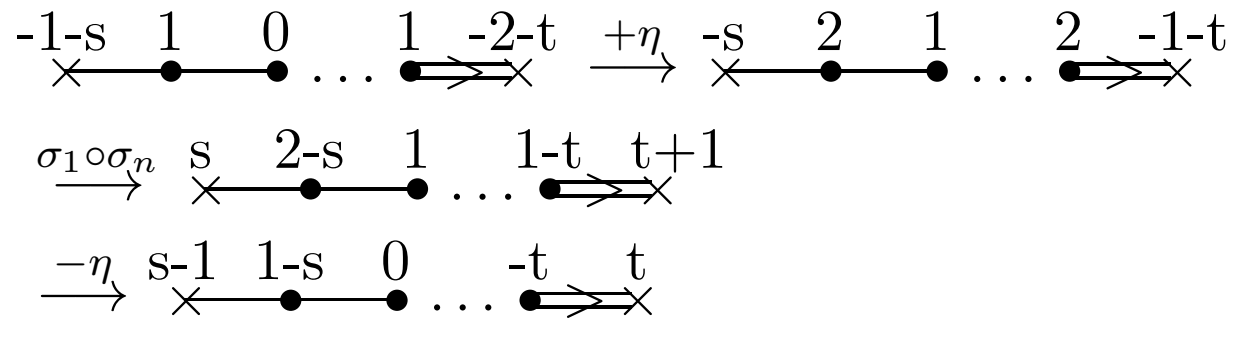

implies $\mathrm{H}^1(X, TX\otimes L^*) = J^2(-\lambda+\mu_1) = 0$ for all s,t$\geqslant$ 1.

(3) $(X, L) =$ 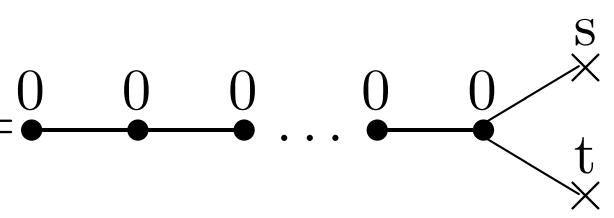

for some s,t$\geqslant 1$. In this case

$\mu_1 =$ 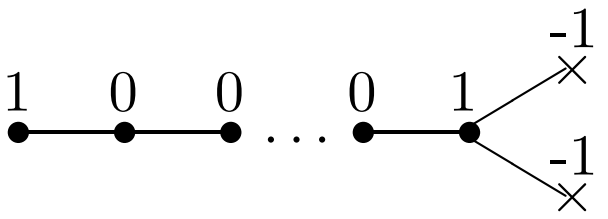

.

The only element of the Weyl group $W$ of length 2 which can, in principle, map $-\lambda + \mu_1 + \eta$ to a strictly dominant weight in $\Lambda^{++}$ is $\sigma_{n-1} \circ \sigma_n$. However, a computation as above shows that

$\sigma_{n-1} \circ \sigma_n(-\lambda + \mu_1 + \eta) - \eta =$ 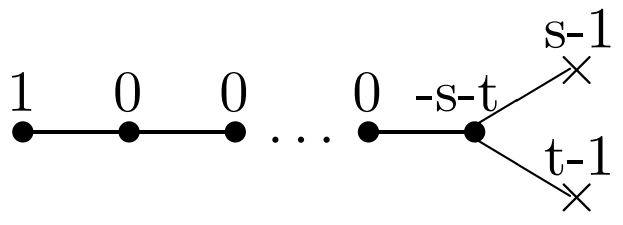

which implies $\mathrm{H}^1(X, TX \otimes L^*) = J^2(-\lambda + \mu_1) = 0$ for all s,t$\geqslant 1$.

*The case* rank$X = 1$. The number $m$ of simple ideals of the semisimple part of the parabolic algebra $\mathfrak{p}$ can, in principle, be equal to 1, 2 or 3. The case $m = 3$, however, is ruled out by Lemma 8.2(iv). If $m = 2$, then, by Lemma 8.2(iv), at least one of the ideals must be isomorphic to $\mathfrak{sl}(2, \mathbb{C})$. Therefore, the crossed node must be either an end node (for $m = 1$) or the node adjacent to an end node (for $m = 2$) of the Dynkin diagram for $X$. Inspecting the list of maximal roots (4) excludes all but the following candidates to $(X, L)$ with $\mathrm{H}^1(X, TX \otimes L^*) \neq 0$:

(1) $(X, L) = (\mathbb{CP}_n, \mathcal{O}(\mathrm{s})) =$ 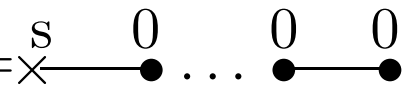

($n \geqslant 5$ nodes) for some s $\geqslant 1$. The odd dimensional projective space has another representation as

$(\mathbb{CP}_{2n-1}, \mathcal{O}(s)) =$ 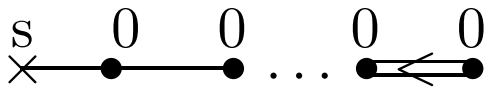

($n \geqslant 3$ nodes).

The long exact sequence

$$0 \longrightarrow \mathcal{O}(-\mathrm{s}) \longrightarrow \mathbb{C}^{n+1} \otimes \mathcal{O}(1 - \mathrm{s}) \longrightarrow TX \otimes L^* \longrightarrow 0$$

implies $\mathrm{H}^1(X, TX \otimes L^*) = 0$ for all s $\geqslant 1$.

(2) $(X, L) =$ 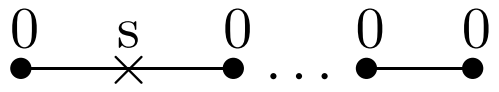

for some s $\geqslant 1$. There are two maximal roots

$\mu_1 =$ 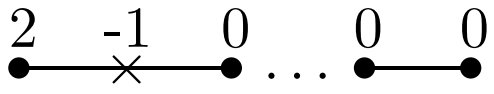

, $\mu_2 =$ 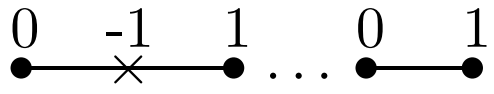

.

That $J^2(-\lambda + \mu_1) = 0$ for all s $\geqslant 1$ follows from the proof of Lemma 8.2(iv), while

$\sigma_1 \circ \sigma_2(-\lambda + \mu_2 + \eta) - \eta =$ 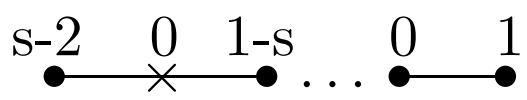

implies $J^2(-\lambda + \mu_2) = 0$ for all s $\geqslant 1$ as well. Thus $\mathrm{H}^1(X, TX \otimes L^*) = 0$ for all s $\geqslant 1$.

(3) $(X, L) =$ 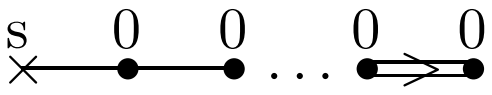 ($n \geqslant 3$ nodes) for some s $\geqslant 1$. Note that for $n = 3$ this pair is biholomorphic to 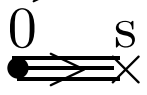 [39] .

The maximal root is

$\mu_1 = \begin{cases} \text{(diagram)} & n \geqslant 4 \text{ nodes} \\ \text{(diagram)} & n = 3 \text{ nodes.} \end{cases}$

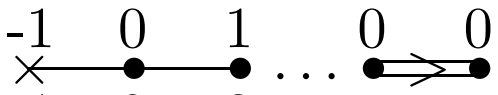
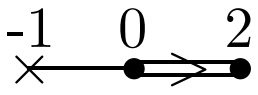

Then an easy computation shows

$\sigma_2 \circ \sigma_1(-\lambda + \mu_1 + \eta) - \eta = \begin{cases} \text{(diagram)} & n \geqslant 4 \text{ nodes} \\ \text{(diagram)} & n = 3 \text{ nodes} \end{cases}$

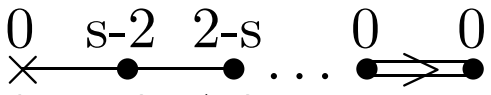
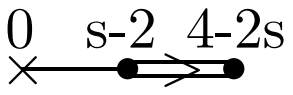

which implies

$$\mathrm{H}^1(X, TX \otimes L^*) = J^2(-\lambda + \mu_1) = \begin{cases} \mathbb{C} & \mathrm{s} = 2 \\ 0 & \text{otherwise.} \end{cases}$$

(4) $(X, L) =$ 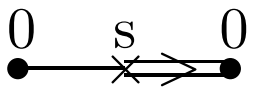 for some s $\geqslant 1$. The maximal roots are

$\mu_1 =$ 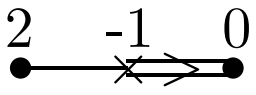 , $\mu_2 =$ 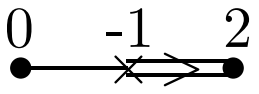 .

Then

$\sigma_3 \circ \sigma_2(-\lambda + \mu_1 + \eta) - \eta =$ 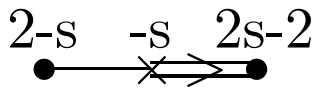 ,

$\sigma_1 \circ \sigma_2(-\lambda + \mu_1 + \eta) - \eta =$ 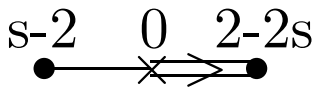 ,

implying $\mathrm{H}^1(X, TX \otimes L^*) = J^2(-\lambda + \mu_1) + J^2(-\lambda + \mu_2) = 0$ for all s $\geqslant 1$.

(5) $(X, L) =$ 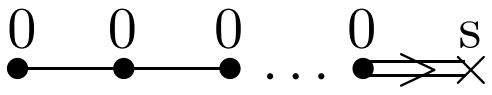 ($n$ nodes) for some s $\geqslant 1$. This pair is biholomorphic to the following one [39]:

$(X, L) =$ 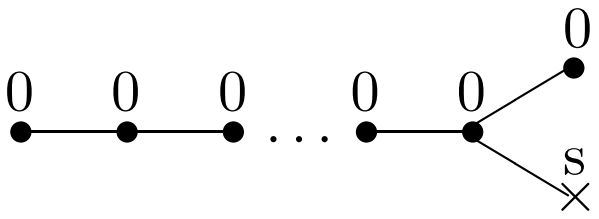 ($n + 1$ nodes).

Then, by Lemma 8.2(iv), $\mathrm{H}^1(X, TX \otimes L^*) = 0$ for all s $\geqslant 1$.

(6) $(X, L) =$ 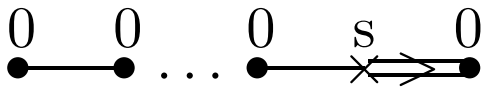 for some s $\geqslant 1$. The maximal roots are

$\mu_1 =$ 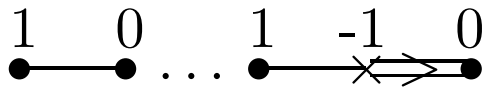 , $\mu_2 =$ 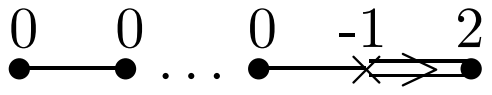 .

The proof of Lemma 8.2(iv) implies $J^2(-\lambda + \mu_2 + \eta) = 0$ for all s $\geqslant 1$, while

$\sigma_n \circ \sigma_{n-1}(-\lambda + \mu_1 + \eta) - \eta =$ 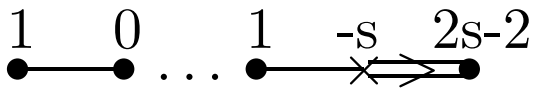

implies $J^2(-\lambda + \mu_1 + \eta) = 0$ for all s $\geqslant 1$. Whence, by Lemma 8.1, $\mathrm{H}^1(X, TX \otimes L^*) = 0$ for all s $\geqslant 1$.

(7) $(X, L) =$ 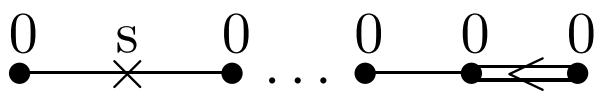 ($n \geqslant 3$ nodes) for some s $\geqslant 1$. The maximal roots are

$\mu_1 =$ 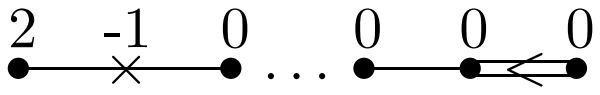 , $\mu_2 =$ 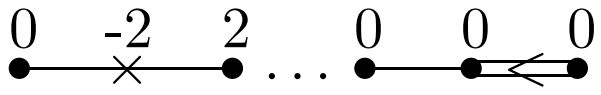 .

The vanishing of $J^2(-\lambda+\mu_1)$ for all $n \geqslant 4$, s $\geqslant 1$ follows from the proof of Lemma 8.2(iv). That this module vanishes for $n = 3$, s $\geqslant 1$ follows from a simple calculation:

$$\sigma_3 \circ \sigma_2(-\lambda+\mu_1+\eta)-\eta \;=$$
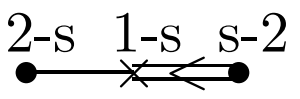

.

Analogously, one finds

$$\sigma_1 \circ \sigma_2(-\lambda+\mu_2+\eta)-\eta \;=$$
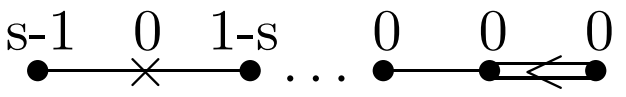

implying $J^2(-\lambda+\mu_2) = 0$ for all s $\geqslant 2$ and $J^2(-\lambda+\mu_2) = \mathbb{C}$ for s $= 1$. Therefore, by Lemma 8.1,

$$\mathrm{H}^1(X, TX \otimes L^*) = \begin{cases} \mathbb{C} & \mathrm{s}=1 \\ 0 & \text{otherwise.} \end{cases}$$

(8) $(X, L) =$
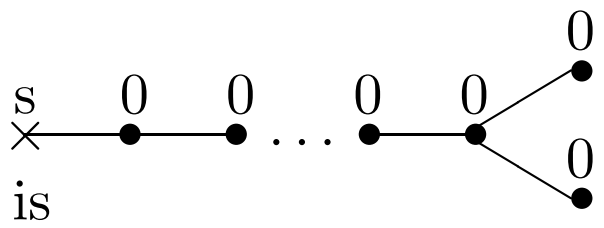

($n \geqslant 5$ nodes) for some s $\geqslant 1$. The maximal root is

$$\mu_1 =$$
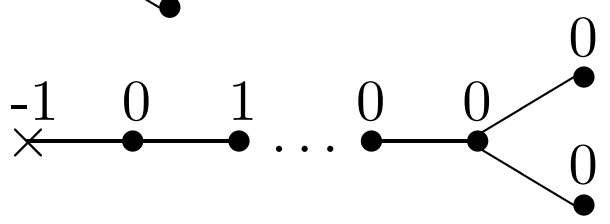

and an easy calculation shows that

$$\sigma_2 \circ \sigma_1(-\lambda+\mu_1+\eta)-\eta \;=$$
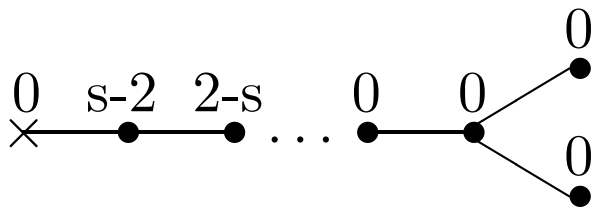

.

Therefore,

$$\mathrm{H}^1(X, TX \otimes L^*) = J^2(-\lambda+\mu_1) = \begin{cases} \mathbb{C} & \mathrm{s}=2 \\ 0 & \text{otherwise.} \end{cases}$$

(9) $(X, L) =$
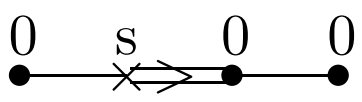

for some s $\geqslant 1$. The maximal roots are

$$\mu_1 =$$
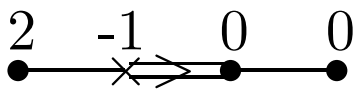

, $\mu_2 =$
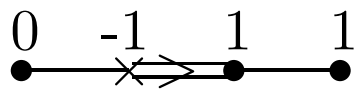

.

From the proof of Lemma 8.2(iv) it follows that $J^2(-\lambda+\mu_1) = 0$ for all s $\geqslant 1$. The only element of the Weyl group of length 2 which can, in principle, make $-\lambda+\mu_2+\eta$ strictly dominant is $\sigma_2 \circ \sigma_1$. Since

$$\sigma_2 \circ \sigma_1(-\lambda+\mu_2+\eta)-\eta \;=$$
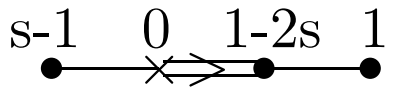

the module $J^2(-\lambda+\mu_2)$ vanishes for all s $\geqslant 1$. Therefore, $\mathrm{H}^1(X, TX \otimes L^*) = 0$ for all s $\geqslant 1$.

(10) $(X, L) =$ 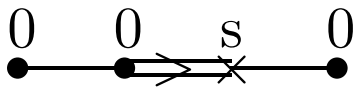

for some s $\geqslant 1$. The maximal roots are

$\mu_1 =$ 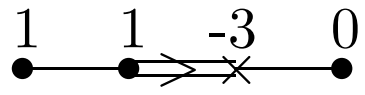

, $\mu_2 =$ 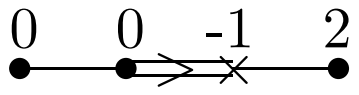

.

From the proof of Lemma 8.2(iv) it follows that $J^2(-\lambda + \mu_2) = 0$ for all s $\geqslant 1$. Since

$\sigma_4 \circ \sigma_3(-\lambda + \mu_1 + \eta) - \eta =$ 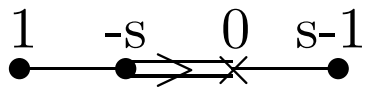

the module $J^2(-\lambda+\mu_2)$ vanishes for all s $\geqslant 1$. Therefore, $\mathrm{H}^1(X, TX \otimes L^*) = 0$ for all s $\geqslant 1$.

(11) $(X, L) =$ 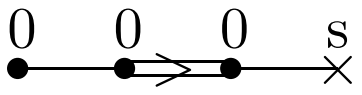

for some s $\geqslant 1$. The maximal root is

$\mu_1 =$ 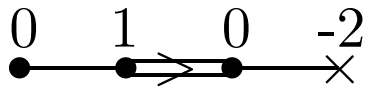

.

Since

$\sigma_3 \circ \sigma_4(-\lambda + \mu_1 + \eta) - \eta =$ 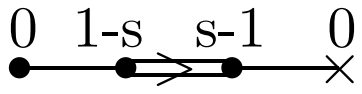

we obtain

$$\mathrm{H}^1(X, TX \otimes L^*) = J^2(-\lambda + \mu_1) = \begin{cases} \mathbb{C} & \mathrm{s} = 1 \\ 0 & \text{otherwise.} \end{cases}$$

(12) $(X, L) =$ 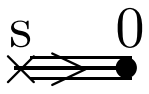

for some s $\geqslant 1$. The maximal root is $\mu_1 =$ 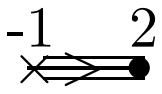

.

Hence

$\sigma_2 \circ \sigma_1(-\lambda + \mu_1 + \eta) - \eta =$ 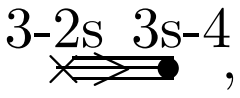

,

implying $\mathrm{H}^1(X, TX \otimes L^*) = J^2(-\lambda + \mu_1) = 0$ for all s $\geqslant 1$.

Theorem B is proved. □

*Acknowledgements.* It is a pleasure to thank D. Alekseevski, L. Berard-Bergery, N. Hitchin, Yu.I. Manin, H. Pedersen, H.-B. Rademacher and K.P. Tod for valuable discussions and comments. Hospitality of MPI in Bonn and ESI in Vienna (SM), University of Glasgow (LS), ICMS in Edinburgh (both of us) during the work on different stages of the project is gratefully acknowledged. The authors are also grateful to the British-German Council for partial financial support (grant No. 877).

Glasgow University , Glasgow, UK
*E-mail address*: sm@maths.gla.ac.uk

Mathematisches Institut, Universität Leipzig, Leipzig, Germany
*E-mail address*: schwachh@mathematik.uni-leipzig.de